%
%

\magnification=1200

\font\titfont=cmr10 at 12 pt

\font\headfont=cmr10 at 12 pt

\font \fr = eufm10



\def\AAA{\bf}
\def\BBB{\bf}

\overfullrule=0in

\def\boxit#1{\hbox{\vrule
 \vtop{%
  \vbox{\hrule\kern 2pt %
     \hbox{\kern 2pt #1\kern 2pt}}%
   \kern 2pt \hrule }%
  \vrule}}
  \def\mathqed{  \vrule width5pt height5pt depth0pt}
\def\imp{\qquad \Rightarrow\qquad}

  \def\harr#1#2{\ \smash{\mathop{\hbox to .3in{\rightarrowfill}}\limits^{\scriptstyle#1}_{\scriptstyle#2}}\ }

\def\intave#1{\int_{#1}\!\!\!\!\! \! \!\!\!\!-}
\def\intavee{{\int \!\!\!\! \!\!-}}

 \def\GG{{{\bf G} \!\!\!\! {\rm l}}\ }

\def\ST{{ST-invariant}}

\def\den{\Theta}
\def\sh{subharmonic}
\def\iiff{\qquad\iff\qquad}
\def\bra#1#2{\langle #1, #2\rangle}
\def\bbf{{\bf F}}

\def\ss{\subset}

\def\half{\hbox{${1\over 2}$}}
\def\smfrac#1#2{\hbox{${#1\over #2}$}}
\def\oa#1{\overrightarrow #1}

\def\dist{{\rm dist}}

\def\log{{\rm log}}

\def\tr{{\rm tr}}
\def\max{{\rm max}}

\def\span{{\rm span\,}}

\def\det{{\rm det}}

\def\Sym{{\rm Sym}^2}

\def\arr{\longrightarrow}

\def\rn{\bbr^n}

\def\Int{{\rm Int}}

\def\Symn{{\Sym(\rn)}}

\def\Theorem#1{\medskip\noindent {\bf THEOREM \bf #1.}}
\def\Prop#1{\medskip\noindent {\bf Proposition #1.}}
\def\Cor#1{\medskip\noindent {\bf Corollary #1.}}
\def\Lemma#1{\medskip\noindent {\bf Lemma #1.}}
\def\Remark#1{\medskip\noindent {\bf Remark #1.}}
\def\Note#1{\medskip\noindent {\bf Note #1.}}
\def\Def#1{\medskip\noindent {\bf Definition #1.}}

\def\Ex#1{\medskip\noindent {\bf Example \bf    #1.}}

\def\pf{\medskip\noindent {\bf Proof.}\ }
\def\qed{\hfill  $\vrule width5pt height5pt depth0pt$}
\def\equdef{\buildrel {\rm def} \over  =}

   \def\cp{{\cal P}}

\def\ch{{\cal H}}

\def\cp{{\cal P}}

\def\gerG{{\fr{\hbox{g}}}}

\def\vf{\varphi}

\def\wt{\widetilde}

\def\and{\qquad {\rm and} \qquad}
\def\arr{\longrightarrow}
\def\ol{\overline}
\def\bbr{{\bf R}}\def\bbh{{\bf H}}
\def\bbc{{\bf C}}

\def\a{\alpha}
\def\b{\beta}
\def\d{\delta}
\def\e{\epsilon}

\def\l{\lambda}

\def\s{\sigma}

\def\D{\Delta}
\def\L{\Lambda}

\def\O{\Omega}

\def\lloc{L^1_{\rm loc}}
\def\dbar{\ol{\partial}}

\def\bo{\partial \Omega}

\def\Symn{\Sym(\rn)}
 
\def\USC{{\rm USC}}
\def\fa{{\rm\ \  for\ all\ }}

\def\cpt{\wt{\cp}}
\def\ft{\wt F}
\def\ob{\overline{\O}}

\def\HLCG{HL$_1$} 
\def\HLPCG{HL$_2$} 
\def\HLPCGG{HL$_3$} 
\def\HLDD{HL$_4$} 
\def\HLPUP{HL$_5$} 
\def\HLDDR{HL$_6$}
\def\HLHP{HL$_7$} 
\def\HLHPP{HL$_8$} 
\def\HLREST{HL$_9$} 
\def\HLPCON{HL$_{10}$} 
\def\HLSURVEY{HL$_{11}$} 
\def\HLBP{HL$_{12}$} 
\def\HLAE{HL$_{13}$} 
\def\HLTangII{HL$_{14}$}
\def\HLADP{HL$_{15}$} 
\def\HLCLASSICAL{HL$_{16}$}

\def\A{\Psi}
\def\Aa{\psi}
\def\intave#1{\int_{#1}\!\!\!\!\!\!\!-\ }

\def\AAA{1}
 \def\AA{2}
 \def\CC{3}
 \def\HH{4} 
 \def\LL{5}
 \def\BB{6}
 \def\KK{7}
 \def\BBB{8}
  \def\FF{9} 
  \def\PP{10}
  \def\FFF{11}
\def\JJ{12} 
\def\JJJ{13}
 \def\MM{14}
  \def\NN{15}

\def\ST{ST-invariant }
\def\T{\Theta}

\centerline
{
\titfont  TANGENTS TO SUBSOLUTIONS}
\medskip

\centerline{\titfont EXISTENCE AND UNIQUENESS,   I }

\bigskip

\centerline{\titfont F. Reese Harvey and H. Blaine Lawson, Jr.$^*$}
\vglue .9cm
\smallbreak\footnote{}{ $ {} \sp{ *}{\rm Partially}$  supported by
the N.S.F. }

 
\centerline{\bf ABSTRACT} \medskip

  {{\parindent= .44in
\narrower  \noindent
There is an interesting potential theory associated to each degenerate elliptic, fully  nonlinear equation $f(D^2u) = 0$. These include all  the potential theories  attached to  calibrated 
geometries.  This paper begins the study of tangents to the subsolutions in these theories, a topic
inspired by the results of Kiselman in the classical plurisubharmonic case.  
Fundamental to this study is a new invariant of the equation, 
called the {\sl Riesz characteristic}, which governs asymptotic structures.
The existence of tangents to subsolutions is established in general, as is the existence of an 
upper semi-continuous density function.   Two theorems establishing the strong uniqueness 
of tangents  (which means every tangent is a  Riesz kernel) are proved. They cover all O$(n)$-invariant  convex cone equations and their complex and quaternionic analogues, with the exception of the homogeneous Monge-Amp\`ere equations, where  uniqueness  fails.
They also cover a large class of geometrically defined subequations  which  includes those coming from calibrations. A discreteness result   for the sets where the density is $\geq c>0$ is also established in any case where strong uniqueness holds. A further result (which is sharp) asserts  the  H\"older continuity of subsolutions when the Riesz characteristic $p$ satisfies $1\leq p<2$.
Many explicit examples are examined.

The second part of this paper is devoted to the ``geometric cases''.
A Homogeneity Theorem and a Second Strong Uniqueness Theorem
are proved, and the tangents in the Monge-Amp\`ere cases are completely classified.

}}

\vfill\eject

\vskip.15in

\centerline{\bf TABLE OF CONTENTS} \bigskip

{{\parindent= .1in\narrower  \noindent

\qquad \AAA.  Introduction.   \smallskip

\qquad \AA.  The Radial Subequations Associated to a Subequation $F$.\smallskip

\qquad \CC.     ST-Invariant Cone Subequations --  The Riesz Characteristic.     \smallskip
 
 \qquad \HH.  Some Illustrative Examples.  \smallskip
 
 \qquad \LL.   $K_p$-Convexity and Monotonicity.
    \smallskip

\qquad \BB.    Monotonicity and Stability of Averages for $F$-Subharmonic Functions.    \smallskip

\qquad \KK.    Densities for $F$-Subharmonic Functions -- Upper Semi-Continuity.   \smallskip

\qquad \BBB.  Maximality of Subharmonics with Harmonic Averages. \smallskip

\qquad   \FF.      Tangents to Subharmonics.   \smallskip

\qquad   \PP.      Uniqueness of Averages of Both Tangents and of Flows.   \smallskip

\qquad   \FFF.     Existence of Tangents.   \smallskip
 
\qquad \JJ. Uniqueness of Tangents.\smallskip

\qquad  \JJJ.   The Strong Uniqueness Theorem I. \smallskip

\qquad \MM.  The Structure of the Sets  $E_c$ where the Density is $\geq c$.  \smallskip

\qquad \NN.  Subequations with Riesz characteristic $1\leq p<2$.  \smallskip

\smallskip

\hskip.5 in
Appendix A.   Subaffine Functions and a Dichotomy.

\smallskip

\hskip.5 in
Appendix B.     Uniform Ellipticity and $\cp(\d)$.

}}

\vskip .4in


\centerline{\headfont \AAA.\ Introduction.}
 \medskip 

The point of  this paper is to  introduce and study tangents for a  wide class of
 degenerate elliptic, fully nonlinear equations of the form $\bbf(D^2) = 0$ in $\rn$. 
It was inspired by   Kieselman's study [K$_1$] (cf.  [K$_2$])) of tangents to 
plurisubharmonic functions  in classical pluripotential theory.
The   aim   is to develop techniques 
 for studying the  behavior, in particular the singular behavior, of subsolutions 
 -- the upper semi-continuous
functions u which satisfy  $\bbf(D^2u) \geq 0$ in the viscosity sense. 
A number of quite general results are obtained.  These include existence,
 uniqueness and ``harmonicity''  of tangents for a wide range of equations.
 Densities for subsolutions are defined and shown to be upper semi-continuous,
 and a structure theorem is proved for the sets where the density is $\geq c>0$.
A key to the analysis is the notion of the Riesz characteristic of the equation.
This invariant is a real number $p\geq1$ which governs the asymptotic behavior
of singularities, and is easily computed in all of the examples, no matter how
degenerate (see Sections \CC \ and  \HH).

For this study we focus on the closed set  $F= \{A \in\Symn : \bbf(A)\geq0\}$
(cf. [Kr], [\HLDD]), and the operator 
$\bbf$ will play no role.
This set  is always assumed to have the following three properties:
\smallskip

(i) \ \ { (Positivity)}\ \ \qquad \ \ $F+\cp \ \ss\ F$ \ \ where $\cp\equiv \{A \geq0\}$.
\smallskip

(ii) \ \ { (ST-Invariance)} \ \ \ $F$ is invariant under a subgroup $G\ss {\rm O}(n)$ 

\hskip 1.6in
which acts transitively on the sphere
$S^{n-1} \ss\rn$.

\smallskip

(iii)  \ \ {(Cone Property)}\ \ \ $t F\ss F$ for all $t\geq0$.
\smallskip
\noindent
A closed  set  $F$ satisfying Positivity is called a subequation, and the viscosity $F$-subsolutions are called $F$-subharmonic functions. Each subequation $F$ has its own potential theory ([\HLDD], [\HLDDR]). 
For some of the results  here, in addition to these three conditions, $F$ is also assumed to
be convex. In this case distribution theory provides an alternate but equivalent 
foundation (Theorem \FF.5) for subsolutions, which is helpful.

The equations covered here include many classical examples coming from real,
complex and calibrated geometry, such as the Monge-Amp\`ere and Hessian equations.
The reader is encouraged to glance at Section \HH \  for some basic examples.

At the time of the first writing of this paper the authors were unaware of its connections 
to the important work of Armstrong, Sirakov and Smart [AS]. They   also studied conical
subequations $F\ss\Symn$ with the additional assumption that $F$ is uniformly elliptic.
This is a stringent assumption which eliminates many of the examples arising from geometry.
They also studied only solutions (as opposed to the much more general subsolutions considered here).
On the other hand they do not assume 
 invariance or convexity, which is extremely nice.
There are also  connections of our work to that of Labutin [La$_1$] who, like Armstrong, Sirakov and Smart,
studied uniformly elliptic equations. 
At the end of this introduction the overlap / lack of overlap is discussed in more detail.

We begin the paper by introducing  the algebraically defined  and easily computable  {\bf Riesz characteristic}
$p_F$ for $F$, which determines much  of the behavior of subsolutions examined here.  
The name comes from the fact that
 when $p\equiv p_F$ is finite, the classical $p^{\rm th}$ {\bf  Riesz kernel} $K_p(|x|)$,
 where 
 $$
K_p(t)\ =\ \cases
{
t^{2-p}          \qquad {\sl if} \ \ 1\leq p\ <2             \cr
\log \, t          \qquad {\sl if} \ \ p=2        \cr
-{1\over t^{p-2}}          \quad {\sl if} \ \  2<p<\infty.           \cr
}
\eqno{(\AAA.1)}
$$
is a solution of the non-linear equation $F$. In fact, every increasing radial solution 
is of the form $\den K_p(|x|) +C$ for constants  $\den\geq 0$   and $C\in\bbr$.
The signs in (\AAA.1) have been chosen so that $K_p(t)$ is always increasing.

 When  $p$  is finite,  there is  an associated tangential $p$-flow
on $F$-subharmonic functions $u$  at each point $x_0$, given for $x_0=0$ by
$$
u_r(x) \ = \ 
\cases
{
\ \ r^{p-2}  u(rx) \qquad\qquad {\rm if}\ \ p\neq 2, \  \ {\rm and}   \cr
  u(rx)  -  M(u,r)    \qquad {\rm if}\ \ p=2, 
}
\eqno{(\AAA.2)}
$$
where 
$$
M(u, r) \ \equiv \  \sup_{|x|\leq r} u.
\eqno{(\AAA.3)}
$$

The {\bf tangents} to   $u$ at $0\in\rn$ are defined to be  the set $T_0(u)$
cluster points of the flow (\AAA.2). When $F$ is convex, these cluster points
are taken in $\lloc(\rn)$. When $1\leq p_F<2$ (but $F$ not necessarily convex),
they can be taken in the local $\b$-H\"older norm for $\b < 2-p$. In either case,
$U\in T_0(u)$ if and only if there exists a sequence $r_j\downarrow0$ such
that $u_{r_j} \to U$ (in the  appropriate space). It is a basic result that tangents
are always entire $F$-subharmonic functions on $\rn$. In particular, the $\lloc$-limits
have unique upper semi-continuous representatives which  are viscosity $F$-subsolutions
(see Theorem \FF.5(b)).  A fundamental result, which is proved in Sections \FFF\ and \NN,  is the following.

\Theorem {\AAA.1. (Existence of Tangents)} {\sl If $F$ is convex or if $p_F<2$, then tangents always exist.}
\medskip

A natural question is whether tangents are actually solutions (as opposed to subsolutions).
The answer is no (if $p_F\geq 2$).  Classical pluripotential theory provides (self) tangent 
examples with large singular sets.  It also provides the remedy -- an appropriate concept enlarging
the space of (viscosity) solutions.

  An $F$-subharmonic function on $X^{\rm open}\ss\rn$
 is called $F$-{\bf  maximal}  if for each $F$-subharmonic function $v$ on $X$ and each
compact subset $K\ss X$,
$$
v\ \leq\ u \quad {\rm on}\ \ X-K   \imp  v\ \leq\ u \quad {\rm on}\ \ X.
$$
If $u$ is $F$-maximal on $X$, then on  any subdomain $Y\ss X$ where $u$ is continuous,
it is a viscosity solution (or ``$F$-harmonic''). In particular, it is always the Perron function for its boundary
values on any ball.  A second fundamental result is the following (see Theorem  \PP.2 and Corollary \PP.3).

\Theorem {\AAA.2. (Maximality of Tangents)} 
{\sl If $F$ is convex, then tangents are always maximal outside the origin in $\rn$.
If  $p_F<2$, then tangents are $F$-harmonic (maximal and continuous) outside the origin.}
\medskip

Existence and regularity (in the weakened form of maximality) for tangents
brings us to the natural question of uniqueness. Here there are several versions.

We say that {\bf uniqueness of tangents} holds for the subequation $F$
if for every  $F$-subharmonic function $u$  defined in a neighborhood of $0$, 
there is exactly one tangent to $u$ at 0.

We say that {\bf strong uniqueness of tangents} holds for  $F$
if for every such $u$,  the unique tangent is $\den(u) K_p(|x|)$, with $\den(u)\geq0$.

We say that {\bf homogeneity of tangents} holds for  $F$ if every tangent to an $F$-subharmonic 
is fixed by the tangential $p$-flow (\AAA.2).
\medskip

Since the flow takes a tangent to $u$ to another tangent to $u$, uniqueness of tangents implies
homogeneity of tangents.

Several important special cases where uniqueness holds are discussed in Section \JJ \
(Propositions \JJ.2,  \JJ.4 and  \JJ.5).

One of the main results of this paper is the Strong Uniqueness Theorem in  Section \JJJ.
Note that there is a natural action of the group O$(n)$ on $\Symn$.
The subequations $F\ss \Symn$ which are O$(n)$-invariant are exactly
those which are defined in terms of the eigenvalues of the matrices $A\in \Symn$.
Every such subequation has a complex and quaternionic counterpart defined 
on $\bbc^n$ and $\bbh^n$ by applying the same eigenvalue constraints to  the complex or quaternionic hermitian symmetric part  of $A$.

\vfill\eject

\Theorem{\AAA.3$\,$I. (Strong Uniqueness of Tangents I)} {\sl
Suppose $F$ is a convex  O$(n)$-invariant subequation, or the complex or quaternionic counterpart of such an equation.  Then, except for  the three basic
cases  $\cp,\cp^\bbc,\cp^\bbh$,  strong uniqueness of tangents holds for $F$.
}

\medskip

There do  exist {\bf non-convex}  O$(n)$-invariant subequations of every Riesz characteristic
for which strong uniqueness fails.  See Example \JJJ.14.

Theorem \AAA.3$\,$I establishes strong uniqueness for a wide range of equations.  These include the
$k^{\rm th}$ Hessian equations  ($k<n$)  and  $p$-convexity equations ($p$ real, $1\leq p\leq n$), 
the trace powers of the Hessian, equations coming from G\aa rding polynomials, 
and much more.  
Each of these has a complex and a quaternionic  counterpart to which Theorem \AAA.3$,$I applies.
However, there are many U$(n)$- and Sp$(n)\cdot$Sp$(1)$-invariant subequations, arising from calibrations
and Lagrangian geometry, which have no O$(n)$-invariant counterpart, so that Theorem \AAA.3$\,$I 
does not apply.     Results in these cases are provided by Theorems
 \AAA.3$\,$II and  \AAA.3$\,$III below, which require a completely different method of proof.

Suppose $F=F(\GG)$ is a subequation defined by a compact subset $\GG\ss G(p,\rn)$
of the Grassmannian of $p$-planes in $\rn$ (see Example \HH.4).

\Theorem{ \AAA.3$\,$II. (Strong Uniqueness II)} {\sl
  Fix $p\geq 2$ and $n\geq 3$.
Then strong uniqueness of tangents  to  $F(\GG)$-subharmonic functions
holds for:
\smallskip

(a) Every compact SU$(n)$-invariant subset 
$
\GG \ss G^\bbr(p, \bbc^n)
$
with the one exception $\GG = G^\bbc(1, \bbc^n)$,
\smallskip

(b) Every compact Sp$(n)\cdot {\rm Sp}(1)$-invariant subset 
$
\GG \ss G^\bbr(p, \bbh^n)
$
with three exceptions, namely 
 the sets of real $p$-planes which lie in 
 a quaternion line for $p=2,3,4$  (when $p=4$ this is $G^\bbh(1, \bbh^n)$),
\smallskip

(c) For $p\geq5$, every compact Sp$(n)$-invariant subset 
$
\GG \ss G^\bbr(p, \bbh^n).
$
 
}

\medskip

This result is based on a companion theorem which has further applications.
Given $\GG\ss G(p,\rn)$ as above, we
 say that $\GG$ has the  {\bf transitivity property} if 
 for any two  vectors $x,y\in\rn$ there exist $W_1,...,W_k\in\GG$
with $x\in W_1, y \in W_k$ and dim$(W_i\cap W_{i+1}) >0$ for all $i=1,...,k-1$. 
The subequations attached to Lagrangian, Special Lagrangian, Associative,
Coassociative, and Cayley geometries all have this property.

\Theorem{ \AAA.3$\,$III. (Strong Uniqueness III)} {\sl
If $\GG$ has the transitivity property, then strong uniqueness of tangents holds
for all $F(\GG)$-subharmonic functions.}

\medskip
Theorems  \AAA.3$\,$II and  \AAA.3$\,$III will be proved in Part II of this paper.
There homogeneity of tangents is  proved first,  and then 
strong uniqueness is established.  This method makes no use of uniform ellipticity,
and has its roots in pluripotential theory, not viscosity theory.

It is important to note that uniqueness of tangents does not always hold.  
In the basic case of convex functions ($F=\cp)$ we have uniqueness, but
strong uniqueness fails.  For classical  plurisubharmonic functions
(the complex counterpart: $F=\cp^\bbc$), the uniqueness  question was
 raised in [H]   and answered in the negative by  Kiselman [K$_1$],
 who actually characterized the sets which can arise as $T_0(u)$ for a plurisubharmonic
 function $u$ in $\bbc^n$.   In Part II of this paper
a similar result is obtained for the quaternionic counterpart $\cp^\bbh$.

The proof of Theorem  \AAA.3$\,$I involves several steps.
The first step is of a classical nature going back to standard potential theory for
the Laplacian and used by Labutin and Armstrong-Sirakov-Smart in viscosity theory.
In our formulation it involves various characterizations of radial $F$-harmonics.
For example, a result (Thm. \AA.4, \AA.7.), straightforward in 
the smooth case, but which fills a gap in the literature, characterizes the radial viscosity subsolutions
$u(x) = \psi(|x|)$ as the subsolutions of the one-variable  subequation
$$
R_F \ :\ \ \psi''(r) + {p_F -1  \over r} \psi'(r)\ \geq\ 0
\eqno{(\AAA.4)}
$$
This classical subequaton is reviewed in detail in Section \LL.
Several important facts are derived. For example, all subsolutions of (\AAA.4) are continuous,
which has the important consequence  that if a radial function is $F$-maximal,
then it is $F$-harmonic (a solution), and hence of the form $\den K_p(|x|)+c$.
Another consequence of (\AAA.4) is that quotients ${\psi(r) - \psi(t)  \over K(r)-K(t)}$
are jointly (or ''doubly'') monotone.  This can be applied to a general non-radial $F$-subsolution
$u$ by associating to $u$ several  radial functions which are also $F$-subharmonic
 (Lemmas \BB.1 and \BB.2).
The simplest is the maximum $M(u,|x|)$ defined by (\AAA.3),
which is a basic tool in [La$_*$] and  [AS].  We choose the following formulation (see Section \BB).

\Theorem {\AAA.4. (Double Monotonicity)} {\sl Let $u$ be $F$-subharmonic in a neighborhood of 
the origin in $\rn$. 
Then 
$$
{M(u,r)-M(u,s)  \over K(r)-K(s)} \quad {\rm is\ increasing\ in} \ \ r\ \ {\rm and}\ \ s.
\eqno{(\AAA.5)}
$$
for all $0<s<r$ where $M$ is defined.

 Furthermore, if $F$ is convex, the same statement holds with $M(u,r)$ replaced by 
 either 
 $$
 S(u,r) \ \equiv\ \intave S u(r\s) \,d\s 
 \qquad{\rm or}\qquad
 V(u,r) \ \equiv\ \intave B u(rx) \,dx 
 \eqno{(\AAA.6)}
$$
(the spherical or volume average)
where $B \equiv \{|x|\leq1\}$ is the unit ball,  $S\equiv \partial B$ is the unit sphere,
and $\intave S  = {1\over |S|} \int_S$ denotes the average or ``normalized'' integral.}

 \medskip
 This theorem has several immediate consequences for the functions $\Psi(u,r)$
 for $\Psi = M,S,V$.    In particular,  it leads to  the concept of densities  (see Corollary \LL.3).
  
  \Def{\AAA.5}  Suppose $u$ is $F$-subharmonic in a neighborhood of 
$0\in\rn$. Then the {\bf $M$-density}  of  $u$ at 0 is the decreasing limit
$$
\den^M(u,0) \ \equiv\ \lim_{s<r \downarrow 0} {M(u,r)-M(u,s)  \over K(r)-K(s)}.
$$
When $F$ is convex, there are also $\Psi$-densities
 $$
\den^\Psi(u,0) \ \equiv\ \lim_{s<r \downarrow 0} {\Psi(u,r)-\Psi(u,s)  \over K(r)-K(s)}.
$$
for $\Psi =S$ and $V$ as in (\AAA.6).

\medskip

Elementary results concerning these densities are established in  Lemma \LL.4.

When $F$ is convex, each $F$ subharmonic function is classically $\D$-subharmonic,
and so $\D u =\mu\geq0$ (a positive measure).  Thus we also have the standard
``mass density'' 
$$
\den^q(\mu,0) \ \equiv\ \lim_{r\downarrow0} {\mu\left( B_r(0)\right)  \over \a(q) r^q}
\qquad{\rm where}\ \ q=n-p.
$$
In this convex case all of the densities for   $M, S,V$ and $\mu$ are universally related, 
and when $p=2$ we have the further result that $\den^M=\den^S=\den^V$
(see Propositions \KK.1 and \KK.2).  

As noted, tangents need not be unique.  However, the averages of tangents are uniquely
determined by the density alone, even in the most degenerate cases.  This is step
two in the proof  of the Stong Uniqueness Theorem \AAA.3$\,$I.
It is also the key step in the proof  of existence (Theorem \AAA.1) and maximality
(Theorem \AAA.2).

In the classical case of pluripotential theory the Riesz characteristic is 2, 
and our next result, when $p=2$, is an extension of the work of Kiselman [K$_1$].

\Theorem{\AAA.6.  (Averages of Tangents)} {\sl
Suppose $F$ is convex and $u$ is an $F$-subharmonic function defined in a neighborhood of
the origin in $\rn$.  Let $p=p_F$ be  the Riesz characteristic of $F$.  If $p \neq2$,
then each tangent $U$ to $u$ at 0 has averages
$$
\eqalign
{
M(r) \ =\ \sup_S U(r\s) \ =\ &\den^M(u) K(r), 
\qquad    S(r)\ =\ \intave S U(r\s) \, d\s  \ =\ \den^S(u) K(r),\cr
 \ \ \  &{\rm and}\ \ \ V(r)\ =\ \intave B U(rx)\, dx  \ =\ \den^V(u) K(r)
 }
\eqno{(\AAA.7)}
$$
In particular,
$$
\den^\A(U)\ =\ \den^\A(u) \qquad {\sl for} \ \ \A\ =\ M, S,\  {\sl or} \ V
\eqno{(\AAA.8)}
$$

When $p=2$,  all the densities of $u$ and  any tangent $U$ to 
$u$ at 0,  agree, and will be simply denoted by $\den= \den(u)$.
Specifically, we have
$$
\den(u) = \den^M(U) = \den^S(U) =  \den^V(U) =  \den^M(u) =  \den^S(u) = \den^V(u).
\eqno{(\AAA.9)}
$$
Moreover, the averages of a tangent $U$ to $u$ are given by}
$$
M(r) \ =\ \den \,\log\, r, 
 \qquad  S(r)\ =  \den\, \log\, r +\intave S U,  
 \ \ \   {\rm and}\ \ \ V(r)\ =\  \den\, \log\, r +\intave B U 
\eqno{(\AAA.10)}
$$
\medskip

This result about spherical averages of tangents  has many applications, 
for example it is enough to prove
maximality of tangents (see Theorem \BBB.2).

\Theorem{\AAA.7. (Maximality Criterion)} 
{\sl 
Suppose $F$ is convex and $U$ is an $F$-subsolution on an annular region $A$ about 0.
If the spherical average $S(U, |x|)$ is an increasing $F$-solution on $A$, then $U$
is maximal on $A$.}

\medskip

Some of the remaining steps in the proof of Theorem \AAA.3$\,$I, which are given 
in detail in Section \JJJ, can be outlined as follows. By applying the maximality criterion
we conclude in Theorem \PP.2 that all tangents are $F$-maximal. 
Now if $F'$ is any subequation which contains $F$ and has the same Riesz characteristic,
then an $F$-tangent $U$ to $u$ is also an $F'$-tangent to $u$.  In the O$(n)$-invariant
(and the other cases of Theorem  \AAA.3$\,$I) it is somewhat surprising that there is a 
simple convex subequation of characteristic $p$ which contains all the others
(Proposition \JJJ.9).  This largest subequation is very nice -- in particular, it is uniformly elliptic.
This, together with Theorem \BBB.7, shows that tangents are harmonic for this largest subequation,
and that they are $C^1$.  One completes the proof of Theorem  \AAA.3$\,$I by showing that
for each tangent $U$ and rotation $g$, we must have $U=g^*U$ or otherwise one can 
produce a tangent which is not  $C^1$.

As with most notions of density in analysis, we have the following.

\Theorem {\AAA.8. (Upper Semi-Continuity of Density)} {\sl 
Suppose $u$ is $F$-subharmonic on an open set $X\ss\rn$.  Then
each of the densities 
$$
\den^M(u,x),  \ \ \den^S(u,x), \ \ \den^V(u,x)
$$
 considered above is an upper semi-continuous
function of $x$. 
Equivalently, for all $c\geq0$ and each $\den$ as above, the sets}
$$
E_c \ \equiv\ \{x : \den(u,x)\ \geq\ c\} \ \ {\sl are \ closed.}
$$

We also note that by standard geometric measure theory
$$
c\ch^{n-p}(E_c) \ \leq \ \mu(X).
$$

In many cases one can say much more about these high density sets $E_c$ for $c>0$.

For classical plurisubharmonic functions in $\bbc^n$  a deep theorem, due to L. H\"ormander, E. Bombieri  and  in its final form by Siu ([Ho$_1$], [B], [Siu]),
 states that  $E_c$ is a complex analytic subvariety.
One  straightforwardly deduces from this result  that for the 2-convexity subequation 
$\cp_2$ in $\bbr^{2n}$ the set $E_c$ is discrete, since $\cp^\bbc(J)\ss \cp_2$ for all orthogonal (parallel) complex structures $J$ on $\bbr^{2n}$.  
This very restrictive corollary has a quite general extension.

 \Theorem {\AAA.9. (Structure of High Density Sets)} {\sl Suppose strong uniqueness of tangents holds for $F$.
 Then for any $F$-subharmonic function $u$,  the set $E_c(u)$ is discrete.  
 }
 \medskip

Theorem  \AAA.9  is essentially sharp.  Suppose $\O$ is a domain with strictly
convex boundary. Given any finite subset  $E=\{x_j\}_{j=1}^N\ss\O$ , any set of numbers
 $\den_j>0$, $j=1,...,N$, and any $\vf\in C(\bo)$, there exists a unique continuous
  $u:\ob \to [-\infty, \infty)$ such that
  \smallskip
  
  (1) \ \ $u$ is $F$-harmonic on $\O-E$,

\smallskip
  
  (2) \ \ $u\bigr|_{\bo} = \vf$, and

\smallskip
  
  (3) \ \ $\den(u, x_j) = \den_j$ for $j=1,...,N$.
\smallskip
 
 \noindent
See Remark \MM.2 for more details.

\smallskip

The subequations with characteristic $1\leq p<2$ are very different in nature
 from those  where $p\geq2$.  They are discussed in detail in Section \NN.
 In particular, the following is proved. 

 \Theorem {\AAA.10. (H\"older Continuity $1\leq p<2$)}   {\sl
 Suppose  $F$ is a (not necessarily convex) 
  subequation with Riesz characteristic  $1\leq p<2$.   Then each $F$-subharmonic function is locally H\"older continuous with exponent $\a\equiv 2-p$.

  Furthermore, if $u$ is an $F$-subharmonic defined in a neighborhood of $0\in \rn$,
  then every sequence $\{u_{r_j}\}_{j=1}^\infty$ with $r_j\downarrow0$, has a subsequence
   which converges locally uniformly to an $F$-subharmonic function $U$ on $\rn$.
   In fact for each  $0<\b<2-p$ there exists a subsequence which converges locally 
   in $\b$-H\"older norm. Finally, when $F$ is convex, this limit  $U$ is $F$-harmonic on $\rn-\{0\}$.
  }
  
  \medskip
For the $k^{\rm th}$ Hessian equation the Riesz characteristic is  $p=n/k$.   For $k>n/2$, the H\"older continuity result  for this subequation is a fundamental theorem of Trudinger and Wang  [TW$_1$], and their proof
can be carried over  to more general convex equations.  However,  we do not require convexity 
in Theorem \AAA.10.

In Appendix A we examine the radial subequation for 
the "subaffine"  case $\cpt \equiv \{\l_{\rm max}\geq0\}$
and establish a basic dichotomy -- the Increasing/Decreasing Lemma.

In Appendix B we show that the subequation $\cp(\d) \equiv \{A+\d \tr(A)\geq0\}$ 
is uniformly elliptic in the  conventional sense.

While in Section \HH \ we give a number of examples to which our theory applies,
many more examples are given in the appendix to Part II. That appendix also 
constructs the maximal and minimal subequations of Riesz characteristic $p$
(showing, in particular, that these largest and smallest subequations exist).
There is a companion result describing the largest and smallest {\sl convex} subequations 
of characteristic $p$.  The largest is given in Proposition \JJJ.9. 
The smallest is given in Lemma A.1 of Part II.

It is  worth noting that
 the main results results in this paper (existence, strong uniqueness, maximality, etc.)
 apply to any subequation obtained by a linear change of variables, i.e., of the form $g^tFg$
for $g\in $ GL$_n(\bbr)$ (where $F$  is as assumed herein).
This means for cone subequations $F$ which are invariant under a conjugate subgroup  
$g^{-1}Gg$ where $G\ss {\rm O}(n)$ acts transitively on $S^{n-1}$.
Of course the notion of Riesz characteristic must be reformulated in this case,
and the Riesz kernel $K_p(|x|)$ must be replaced by its transform $K_p(|gx|)$.

\vskip .3in
\centerline
{
The Work of Armstrong, Sirakov and Smart
}
 \medskip
 In the very interesting paper [AS] the authors also study conical subequations
 $F\ss\Symn$ with the added assumption that $F$ is uniformly elliptic.  
 However, they do not  assume invariance or convexity.  An important part of their work
 (which is "automatic" in our case) 
proves the existence and uniqueness of  {\sl fundamental solutions}  --  $F$-harmonic functions $\Phi$ on 
 $\rn-\{0\}$, which are invariant under the flow $\Phi_r(x)= r^{p-2}\Phi(rx)$  for some $p\geq 1, p\neq 2$
 and bounded from above or below.   (When $p=2$ the log enters as it does here.)
 They show the existence and uniqueness of two families of  such solutions 
 (up to positive scalars and additive
 constants) among all entire punctured $F$-harmonics with a one-sided bound.
 In our degenerate cases two fundamental solutions are not always   available.
 In fact, they are if and only if both $F$ and the dual $\ft$ have finite Riesz characteristics.
 (See Proposition \CC.16 for a description of all such subequations.)

 One of the results in [AS] is closely related to the work here.  They prove existence and 
 strong uniqueness of tangents to { solutions} of uniformly elliptic equations. That is, under
 their assumptions that  $F$  is conical and uniformly elliptic,  they prove that:
 Any $F$-{\bf harmonic} function defined  on $B_\e - \{0\}$ and bounded above (or below), 
 has a unique tangent of the form  $\T \Phi$ for some $\T\geq 0$ (see [AS; (5.13) ff.]).

This paper addresses a much broader class of functions, namely subsolutions
to degenerate elliptic equations.  Naturally the equations must be in some ways
 restricted, but the results apply to a wide range of geometrically interesting cases.
 Here it is shown that tangents exist and are maximal,
 and that maximal plus continuous implies F-harmonic.
However, it is not true that maximal implies continuous in this general case.
It fails for example for $\cp_\bbc$,  as does  uniqueness
of tangents (not just strong uniqueness, see  Kiselman [K$_1$]).

Said differently, the step from maximal to $F$-harmonic does not always 
hold in the degenerate subharmonic case, and it is somewhat surprising that
strong uniqueness of tangents can actually be established for such a broad spectrum of 
interesting subequations with $p\geq2$.

We should add that  the techniques used in proving strong uniqueness in the 
non-O$(n)$-invariant cases are substantially different from
those  in the O$(n)$-case, and they appear in the sequel (Part II) to this paper.

 For the question of existence we need to assume convexity or that $1\leq p<2$.
This is quite reasonable since we are dealing with subsolutions and the equations are
only degenerate elliptic.  One needs a function space in which to extract convergent subsequences
just to get off the ground. These assumptions provide such spaces, namely $\lloc$ and H\"older.

The work in [AS] is related to earlier  results of Labutin [La$_*$] who studied the 
Pucci extremal equations.  He established among other things a removable singularity
result and an extension of a classical result of B\^ocher.  In this work the classical Riesz kernels
also play a prominent role.  There is a careful account of the
relationship to the work of Armstrong-Sirakov-Smart  given in [AS].

\vskip .3in
\centerline
{
Historical Reflections
}
 \medskip

 In 1982 the authors showed that for each calibration on a riemannian manifold 
 there is an associated family of minimal
 subvarieties -- forming a {\sl calibrated geometry} [\HLCG].  More recently [\HLPCG] it was
 discovered that the calibration also determines a potential theory of functions whose restrictions
 to each of the distinguished submanifolds are subharmonic. Although there is an analogue
in this setting of the $i\partial \dbar$  operator from complex geometry, 
that operator does not play a critical role in the development of the potential theory [\HLDD].
In fact, somewhat surprisingly,  a corresponding  potential theory can be established for any
collection of submanifolds determined by requiring their tangent spaces to be 
in an arbitrary given closed subset of the grassmannian.  Even more generally one has the potential
theory associated to an elliptic (possibly degenerate) nonlinear inequality $F(D^2u)\geq0$,
provided by viscosity subsolutions ([CIL]).

This raises the possibility of cross-fertilization between two well established and deep fields,
pluripotential theory (in several complex variables) and nonlinear elliptic theory.  This paper,
although not the first, can be viewed as an example of this phenomenon. The authors 
believe there are many more to come.

\vfill\eject


\centerline{\headfont \AA.\ The Radial Subequations Associated to a Subequation $F$.}
 \medskip 

In this section we first describe the ordinary differential inequality which governs $C^2$
radial (i.e., spherically symmetric) $F$-subharmonic functions.  
Our main result fills an apparent gap in the literature by extending
 this characterization to general upper semi-continuous radial $F$-subharmonics.
 Somewhat surprisingly this extension requires the attention of Lemma \AA.10 below.

Suppose $\psi(t)$ is of class $C^2$ on an interval contained in the positive real numbers.
We also consider $\psi$  as the function $\psi(|x|)$ of $x$ on the corresponding annular region in $\rn$.

\Lemma{\AA.1} {\sl
$$
D^2_x \psi \ =\ {\psi'(|x|)\over |x|} P_{[x]^\perp} + \psi''(|x|) P_{[x]}.
\eqno{(\AA.1)}
$$
where $P_{[x]} = {x\circ x \over |x|^2}$ 
denotes orthogonal projection onto the line $[x]$ 
through $x\neq0$ and $P_{[x]^\perp} =I- P_{[x]}$ denotes orthogonal projection onto
the hyperplane with normal $[x]$.
}
\pf
First note that $D(|x|)= {x\over |x|}$, and therefore  $D^2(|x|) = D({x\over |x|})
 = {1\over |x|} I - {x\over |x|^2}\circ {x\over |x|}
= {1\over |x|}(I- P_{[x]} ) = {1\over |x|} P_{[x]^\perp}$.
Hence, 
$$
D_x\psi \  =\  \psi'(|x|) {x\over |x|} 
\qquad {\rm and}
$$ 
$$D_x^2\psi \ \ =\  \ \psi'(|x|) 
D\left( {x\over |x|} \right) + \psi''(|x|)  {x\over |x|}\circ  {x\over |x|}
\ \ =\  \ 
{\psi'(|x|) \over |x|} P_{[x]^\perp} + \psi''(|x|)P_{[x]}. \qquad\qquad
\mathqed
$$

\Cor{\AA.2}  {\sl
The second derivative $D^2_x \psi$ has  eigenvalues ${\psi'(|x|) \over |x|}$
with multiplicity $n-1$ and   $ \psi''(|x|)$ 
with multiplicity 1.
}\medskip

Let  $F\ss \Symn$ be a pure second-order constant coefficient subequation.
Then by Lemma \AA.1 a radial $C^2$-function $u(x) = \psi(|x|)$ is $F$-subharmonic on an 
annular region in $\rn$ if and only if 
$$
D^2_x u\ =\ 
{\psi'(t) \over t} P_{e^\perp} +  \psi''(t) P_e \in F,\ 
\eqno{(\AA.2)}
$$
for $t=|x|$ in the corresponding interval in $(0,\infty)$.
We use $\l = \psi'(t)$ and $a=\psi''(t)$ as one-variable jet coordinates. Then the
basic one-variable subequation associated with $F$ is defined as follows.

\Def{\AA.3}  The {\bf radial subequation associated with} $F$ is the reduced  variable coefficient 
subequation $R_F$ on $(0,\infty)$ whose fibre at $t$ is
$$
(R_F)_t \ \equiv\ \left \{(\l,a)\in \bbr^2 : {\l \over t}P_{e^\perp} + aP_e  \in F,\ \forall\, |e|=1 \right\}.
$$
Thus for $C^2$-functions we have that
$$
u(x) \equiv \psi(|x|) \ \ {\rm is}\ \ F\ {\rm subharmonic}
\quad   \iff    \quad
\psi(t) \ \ {\rm is}\ R_F \ {\rm subharmonic} 
\eqno{(\AA.3)}
$$
  
We extend  this  to the viscosity setting where 
$F$-subharmonic functions  are just upper semi-continuous (see [C], [CIL], [HL$_{4,6}$] for definitions).
The proof given below of the implication $\Rightarrow$ is elementary, whereas the proof of
$\Leftarrow$ will require a lemma. Note that the equivalence: $u(x) = \psi(|x|)$ is upper semicontinuous $\iff$
$\psi(t)$  is upper semicontinuous, is obvious.

\Theorem{\AA.4. (Radial Subharmonics)} {\sl
The function $u(x) \equiv \psi(|x|)$ is $F$-subharmonic on an  annular region in $\rn$
if and only if 
 $\psi(t)$ is $R_F$-subharmonic on the corresponding open sub-interval of $(0,\infty)$.
 }

 \medskip

  \Remark{\AA.5}   In all but this section of the paper, the  subequations $F$  will  be
  assumed to be cones, 
   unless explicitly stated to the contrary.
    For such subequations the maximum principle
  holds, i.e., it holds for each $F$-subharmonic function $u(x)$ (see Theorem A.2).
  Consequently, if $u(x) = \psi(|x|)$ is a radial $F$-subharmonic on a  ball about 0,
  then $\psi(t)$ must be increasing in $t$.  This
  motivates focusing on an ``increasing''  version of Theorem \AA.4.
  
   \medskip

 We will use the  fact, which  is elementary to establish, 
 that for an upper semi-continuous function $\psi(t)$,
 $$
 \psi(t) \quad {\rm is \ increasing} \qquad\iff\qquad \psi\ \ {\rm is}\ \{\l\geq0\}\,{\rm subharmonic}.
 \eqno{(\AA.4)}
$$
(See [\HLCLASSICAL] for a proof.)

\Def{\AA.6}  The {\bf increasing radial subharmonic equation $R_F^{\uparrow}$ on $(0,\infty)$} is
defined by
$$
R_F^{\uparrow}  \ =\ R_F \cap  \{  \l \geq0\}. 
\eqno{(\AA.5)}
$$

In light of (\AA.3), it is obvious that for $C^2$-functions $\psi(t)$:
$$
\psi(t) \ \ {\rm is \ } R_F^{\uparrow} \, {\rm subharmonic}  \quad \iff\quad \psi(|x|) \ \ {\rm is \ } F\cap\{x\cdot p\geq0\}  \,{\rm subharmonic}
\eqno{(\AA.6)}
$$
where the variable coefficient first-order subequation $\{x\cdot p\geq0\}$ is the constraint
$x\cdot D_xu\geq0$ on $C^2$-functions.
The equivalence (\AA.6)  can be extended as  in Theorem \AA.4.

\Theorem{\AA.7. (Increasing Radial Subharmonics)}  {\sl
The function $u(x) \equiv \psi(|x|)$ is  an increasing, radial $F$-subharmonic  function if and only if 
$\psi(t)$  is $R_F^{\uparrow}$-subharmonic. }
\medskip


 
  \Remark{\AA.8}  We will sometimes blur the distinction between $\psi(t)$ and  $u(x) =\psi(|x|)$
 by calling $\psi(t)$ a radial (or increasing radial) $F$-subharmonic.
 
 \Remark{\AA.9}  The statement and proof of a theorem analogous to \AA.7 for
 decreasing radial subharmonics is left to the reader.

\medskip
\noindent
{\bf Proof of Theorem \AA.4.  ($\Rightarrow$):}  Suppose $u(x) \equiv \psi(|x|)$ is $F$-subharmonic.
If $\vf(t)$ is a test function for $\psi(t)$ at $t_0$, then $\vf(|x|)$ is a test function
for  $\psi(|x|)$ at any point on the $t_0$-sphere in $\rn$.
Therefore $D^2_{x_0}\vf \in F$.  Applying the formula for $D^2_{x_0}\vf$ in terms of 
$\vf'(t_0)$ and $\vf''(t_0)$, the equivalence (\AA.3),
 and the definition of $(R_F)_{t_0}$,  we have $J^2_{t_0} \vf \in R_F$.
This proves that $\psi(t)$ is $R_F$-subharmonic. 

\noindent
{\bf  ($\Leftarrow$):}   Suppose that $\psi(t)$ is $R_F$-subharmonic.  We must show 
that $u(x) \equiv \psi(|x|)$ is $F$-subharmonic. That is, given a test function $\vf(x)$ for
$u(x)$ at a point $x_0$, we must show that $D^2_{x_0}\vf \in F$.

Suppose that there exists a {\sl smooth} function $\overline\psi(t)$, defined near $t_0=|x_0|$,
such that $\overline\vf(x) \equiv \overline\psi(|x|)$ satisfies
$$
u(x)\ \leq\ \overline\vf(x) \ \leq \vf(x)
\eqno{(\AA.7)}
$$
near $x_0$.  Then $\overline\psi(t)$ is a test function for $\psi(t)$ at $t_0$.
Hence, the 2-jet of $\overline\psi$ at $t_0$ belongs to $R_F$. 
By Lemma \AA.1 and the discussion above, this implies that $D^2_{x_0}\overline\vf \in F$.
The inequality $\overline\vf(x) \leq\vf(x)$ (with equality at $x_0$) implies
that $D^2_{x_0}\vf = D^2_{x_0}\overline\vf +P$ for some $P\geq0$,
which proves that $D^2_{x_0}\vf\in F$ as desired.

To complete this argument by finding $\overline\psi(t)$ there is some flexibility given by Lemma 2.4 in [\HLDDR]
so that not all test functions $\vf(x)$ need  be considered.  
First we may choose new   coordinates
$z=(t,y)$ near $x_0$ so that $t\equiv |x| $. (Thus $t=$   constant defines the sphere of radius $t$ near $x_0$.)
Furthermore, we may assume that $\vf(z)$ is a polynomial of degree $\leq 2$ in 
$z=(t,y)$ and that it is a {\sl strict} local test function, i.e., $u(z)  < \vf(z)$ for $z\neq z_0$.
Now Lemma \AA.10 below  ensures the existence of $\overline\vf(x) = \overline\psi(|x|)$
satisfying (\AA.7). \qed

\medskip

Let $z=(t,y)$ denote standard coordinates on $\rn=\bbr^k\times \bbr^\ell$.
Fix a point $z_0=(t_0,y_0)$ and let $u(t)$ be an upper semi-continuous function
(of $t$ alone) and $\vf(z)$ a $C^2$-function, both defined in a neighborhood of $z_0$.

\Lemma {\AA.10} {\sl
Suppose $u(t) < \vf(z)$ for $z\neq z_0$ with equality at $z_0$. 
If $\vf(z)$ is a polynomial of degree $\leq2$, then there exists 
a polynomial $\overline\vf(t)$ of degree $\leq2$ with
$$
u(t)  \ \leq\ \overline\vf(t) \ \leq\ \vf(z) 
\qquad{\rm near}\ \ z_0.
\eqno{(\AA.8)}
$$
}
\pf
We may assume $z_0=0$ and $u(0)=\vf(0)=0$. Then
$$
\vf(z) \ =\  \bra pt + \bra qy + \bra {At} t + 2 \bra {Bt} y + \bra {Cy}y.
$$
We assume $u(t)  < \vf(t,y)$ for $|t| \leq \e$ and $|y| \leq \d$ with $(t,y)\neq (0,0)$.

Setting $t=0$,  we have $0=u(0) < \bra qy +\bra {Cy}y$
for $y\neq 0$ sufficiently small.  Therefore, $q=0$ and $C>0$ (positive definite). 
Now define
$$
\overline\vf(t)\ \equiv\ \bra pt + \bra {(A-B^tC^{-1} B)t}{t}.
\eqno{(\AA.9)}
$$
The inequalities  in (\AA.8) follow from the fact that for $t$ sufficiently small,
$$
\overline\vf(t)\ = \ \inf_{|y|\leq \d} \vf(z)\ =\ \bra  pt + \bra {At}t + \inf_{|y|\leq \d} \{ 2\bra {Bt} y+\bra{Cy}y\}.
\eqno{(\AA.10)}
$$
To prove (\AA.10) fix $t$ and consider the function $2\bra {Bt} y+\bra{Cy}y$.
Since $C>0$, it has a unique minimum point at the critical point $y= - C^{-1}Bt$. 
The minimum value is $-\bra{B^tC^{-1}Bt} t$. If  $t$ is sufficiently small, the critical point
$y$ satisfies $|y|<\d$, which proves (\AA.7).\qed
 \medskip

\noindent
{\bf Proof of Theorem \AA.7.}
 The arguments given for  Theorem \AA.4  along with the following missing steps 
 provides the proof.  If  $\vf(t)$ is a 
test function for $\psi(t)$ at a point $t_0$, then $\vf(|x|)$ is a test function for $\psi(|x|)$  at $x_0$
whenever $|x_0|=t_0$.  Now
$$
D_{x_0}  \vf \ =\ \vf'(|x_0|) {x_0\over |x_0|}\quad {\rm and\ hence}\quad
x_0\cdot D_{x_0}\vf \ =\ |x_0|\vf'(|x_0|).
\eqno{(\AA.11)}
$$
Hence, if $\psi(|x|)$ is $\{p\cdot x\geq0\}$-subharmonic, then $\psi(t)$ is $\{\l\geq0\}$-subharmonic, and thus increasing.  Conversely,
if $\psi(t)$ is increasing and $\vf(x)$ is a test function for $\psi(|x|)$ at $x_0$,
then $\overline \vf (t) \equiv \vf({tx_0\over |x_0|})$ is a test function for $\psi(t)$ at $t_0=|x_0|$.
Hence, ${\overline \vf }'(t_0)\geq 0$.  However, ${\overline \vf }'(t_0)= (D_{x_0}\vf)\cdot x_0$.\qed

\vskip.3in


\centerline{\headfont \CC.   ST-Invariant Cone Subequations -- The Riesz Characteristic}
\medskip 

This section is devoted to investigating the cone subequations which satisfy a weak form of invariance
which will be referred to as {\bf spherical transitivity  (ST)}.  Two characteristic numbers 
$(p,q)$ will be associated with each such  subequation $F$.  They uniquely determine the radial subequation
for $F$ and, as we  shall show in this and the following sections, can be easily computed in any example.
Moreover, we give a complete description of all possible examples (of ST-invariant subequations
with characteristics $(p,q)$) in the second subsection here.  Most readers will prefer to come back  to this subsection.  Although it adds important perspective to the scope of ST-cone subequations, it is not used in
the subsequent results of the paper.

Recall from the introduction that  a subequation $F\ss\Symn$ is said to be {\bf  \ST}
if there exists a subgroup $G\ss {\rm O}(n)$ 
which acts transitively on the sphere $S^{n-1}\ss\rn$ and leaves 
$F$  invariant (under the induced action of $G$ on $\Symn$).
\medskip

For an \ST cone subequation $F$, 
$$
{\rm the \  slices\ \ }  F\cap \span \{P_{e^\perp},  P_e\} \ \ {\rm  for}\ \  e\in S^{n-1} \ {\rm are\ all\  isomorphic.}
\eqno{(\CC.1)}
$$
Note that $ \span \{P_{e^\perp},  P_e\} = \span \{I,  P_e\}$
and that the induced action on $\Symn$ sends $P_e$ to $P_{g(e)}$.
In particular,
$$
\l P_{e^\perp} +\mu  P_e \ \in\ F \ \  \ {\rm  for \ one }\ \  e\in S^{n-1}  
\quad\iff\quad
\l P_{e^\perp} +\mu  P_e \ \in\ F \ \  \ {\rm  for \ all }\ \  e\in S^{n-1}.
\eqno{(\CC.2)}
$$
 This weakening of \ST   will be referred to as {\bf weak invariance}.

 \bigskip
 
 \centerline{\bf
 The Riesz Characteristics
 }
 
 \medskip

 We begin by focusing on the first of the two characteristics $(p,q)$.
 Although there is an abundance of  interesting \ST cone subequations in dimensions $\geq3$, there are not many
 increasing radial subequations.  In fact they are described by a single ``characteristic'' number 
 $p$ between 1 and $\infty$, which determines a one-variable subequation as follows.

 \Def{\CC.1}  For each $p$ with  $1\leq p <\infty$, 
the {\sl  increasing radial subequation} $R_p^{\uparrow}$ is defined by
 $$
 R_p^{\uparrow}   \ :\ \ a+ {(p-1) \over t} \l \ \geq\ 0 \and \l\ \geq\ 0,
  \eqno{(\CC.3)}
$$ 
while  for $p=\infty$, the subequation $R_\infty^{\uparrow}$ is first-order and defined
 by $ R_\infty^{\uparrow}  = \{\l \geq0\}$.

\def\op{\overline{p}}

 \Def{\CC.2. (The Increasing Riesz Characteristic)} Suppose $F$ is an \ST cone subequation.
 The {\bf increasing characteristic  $p_F$ of $F$} is defined to be 
 $$
 p_F \ \equiv \ \sup\{\op : P_{e^\perp} - (\op-1)P_e \in F\}.
  \eqno{(\CC.4a)}
$$ 
Equivalently, for finite Riesz characteristic, $p_F$ is the unique number $p$ such that
 $$
P_{e^\perp} - (p-1)P_e \  \in \ \partial F.
  \eqno{(\CC.4b)}
$$

\Prop{\CC.3. (Increasing)} {\sl
 Suppose that $F$ is an \ST cone subequation.
 Then  the increasing radial subequation $R_F^{\uparrow}$
 equals $R_p^{\uparrow}$ where $p=p_F$ is the increasing Riesz characteristic of $F$
 }

\pf
Using Definitions \AA.3, \AA.6, \CC.1 and \CC.2, we must show that for $\l\geq0$
$$
{\l\over t} P_{e^\perp} + a P_e \ \in\ F  \iiff  a+{p-1\over t}\l \ \geq\ 0.
$$
Set $-(\op-1) \equiv {at/ \l}$, so that 
$
{\l\over t} P_{e^\perp} + a P_e \in F \iff   P_{e^\perp}  -(\op-1)P_e \in F  
$.
Then  $\op\leq p \iff -{at\over \l}  \leq p-1 
\iff a+{p-1\over t}\l\geq0$.\qed

\medskip

Note that by Definition \CC.2, the positivity condition for $F$, and  the fact that $0\in F$,
we have that $p_F\geq 1$.  Thus  $1\leq p_F\leq \infty$.

The only equation with $p_F=1$ is $\cp$.  At the other extreme we have $p_F=\infty$.  Here there
is a test which is very simple to apply in all the \ST \ examples, namely: 
$p_F =\infty$ iff $-P_e \in F$.
Hence, determining when $p_F<\infty$ is also simple, namely: $p_F <\infty$ iff $-P_e \notin F$.

\Lemma{\CC.4} {\sl  For \ST cone subequations $F$
\smallskip
(a) \ \ $p_F =1 \quad\iff\quad P_{e^\perp} \in \partial F \quad\iff\quad F=\cp$
\smallskip
(b) \ \ $p_F =\infty \quad\iff\quad -P_{e} \in  F \quad\iff\quad  -P_{e} \in  \partial F$.
\smallskip
(c) \ \ $p_F <\infty  \quad\iff\quad -P_{e} \notin  F \quad\iff\quad P_{e} \in  \Int \ft \quad\iff\quad  \cp \ss \Int \ft$.
}
\medskip

Actually, as noted above, it is easy to compute the exact value of $p_F$ in all
the examples.

\medskip
\noindent
{\bf Proof of (a).}  Note first that $p_F>1 \iff  P_{e^\perp} -\e P_e \in F$ for all small  $\e>0$.
Now if $F$ contains an element $A$ with at least one eigenvalue strictly negative, 
then by positivity and the cone property there is an element $A' = P_{e^\perp} -\e P_e \in F$.
Hence $F\neq \cp \Rightarrow p_F>1$. 

\medskip
\noindent
{\bf Proof of (b).}  Note first that $-P_e \in F \ \  \Rightarrow\ \   \a P_{e^\perp} - P_e \in F \ \ \forall\, \a\geq0
 \ \  \Rightarrow\ \  P_{e^\perp} - (p-1)P_e \in F \ \ \forall\,  p\geq1
 \ \  \Rightarrow\ \  p_F =\infty.$
 On the other hand $-P_e \notin F \ \  \Rightarrow\ \   \e P_{e^\perp} - P_e \notin F \  \  \forall\, \e\geq0
 \ {\rm small}
 \ \  \Rightarrow\ \  P_{e^\perp} - (p-1)P_e \notin F \ \  \forall\,  p \ {\rm large}
 \ \  \Rightarrow\ \  p_F <\infty$.   To complete the proof of (b) note that $-P_e \in \Int F$ cannot occur
 unless $F=\Symn$ since $-P_e \in \Int F \Rightarrow 0\in \Int F$. 
 
\medskip
\noindent
{\bf Proof of (c).} Since  $\sim(-F) = \Int \ft$, the first part of (c) follows from the first part of (b).
For any subequation $G$ (such as $\ft$), $A\in \Int G \ \ \Rightarrow\ \ A+\cp \ss \Int G$.
Finally, $P_e +\cp = \cp$, proving that $P_e \in \Int \ft \ \ \Rightarrow\ \ \cp\ss \Int \ft$.
 \qed
\medskip

The primary application of the Riesz characteristics  (and the reason for choosing the name)
is the fact     that the solutions 
of the associated increasing radial equation $R_p^{\uparrow}$  are given by the Riesz kernels.

\Prop{\CC.5} {\sl An \ST\ subequation  $F$ has finite Riesz characteristic $p=p_F$ if and only if the increasing radial harmonics for $F$ are:
$$
\den K_p(|x|) + C
\eqno{(\CC.5)}
$$
where $\den \geq0$, $C\in\bbr$,  and $K_p(t)$ is the $p^{\rm th}$ Riesz function defined on 
$1\leq t <\infty$ by }
$$
K_p(t)\ =\ \cases
{
t^{2-p}          \qquad {\sl if} \ \ 1\leq p\ <2             \cr
\log \, t          \qquad {\sl if} \ \ p=2        \cr
-{1\over t^{p-2}}          \quad {\sl if} \ \  2<p<\infty.           \cr
}
\eqno{(\CC.6)}
$$
\pf
From (\CC.4b) it is easy to see that $u(x) \equiv \psi(|x|)$ is $F$-subharmonic if and only if 
$\psi(t)$ is $R_p^{\uparrow}$-subharmonic.
The ordinary differential equation given by equality in (\CC.3) is easily solved,
and $\den K_p(t)+C$ are the viscosity solutions.
One can  check directly using Lemma \AA.1 that
$$
D^2 \overline K_p(|x|)\ =\  {1\over |x|^p}\left(  P_{[x]^\perp} - (p-1)P_{[x]}   \right)
\and
D\overline K_p \ =\ {x\over |x|^p}
\eqno{(\CC.7)}
$$
where $K_p$ has been renormalized to 
$$
\overline K_p\  \equiv\  \smfrac 1 { |p-2|} K_p \quad {\rm if\ \ } p\neq 2
\ \ {\rm and\ }\ 
\overline K_2 \ =\ K_2\  \qquad\mathqed
\eqno{(\CC.8)}
$$

 The sign of $K_p(t)$ has been chosen so that $K_p(|x|)$ is a {\sl increasing} or {\sl downward-pointing}
 $F$-harmonic on $\rn-\{0\}$.
 The actual  normalization in (\CC.6) is simpler when the  focus is on the function $u$, while the
 normalization in (\CC.8) is simpler when the  focus is on the first and second derivatives
 of $u$.



%

The second of the two numbers $(p,q)$  can also be defined in several equivalent ways.

\Def{\CC.6. (The Decreasing Riesz Characteristic)}
For each \ST cone subequation $F$, this characteristic, denoted $q_F$, is defined by
$$
q_F \ =\ \sup\left\{ \bar q : -P_{e^\perp} + (\bar q-1) P_e \notin F       \right\},  
\eqno{(\CC.9a)}
$$
or equivalently  $q_F$ is the unique number $q$ such that
$$
-P_{e^\perp} + (q-1) P_e \ \in\ \partial F,
\eqno{(\CC.9b)}
$$
or finally, $q_F$ can be defined to be the increasing characteristic of the dual subequation, i.e.
$$
q_F\ =\ p_{\wt F}
\eqno{(\CC.9c)}
$$

Since $\partial \ft=-\partial F$, the equivalence of  (\CC.9c) follows easily from (\CC.4b).
Thus the decreasing characteristic of $F$  might also be called the {\bf dual characteristic} of $F$.

For each $1\leq q<\infty$ set
$$
R^{\downarrow}_q \ :\ \ a + {q-1\over t}\l \ \geq\ 0 \and \l\ \leq\ 0
\eqno{(\CC.10)}
$$
while for $q=\infty$ the subequation $R^{\downarrow}_q$ is first-order
and defined by $R^{\downarrow}_\infty = \{\l\leq0\}$.

Then the decreasing versions of Propositions \CC.3, Lemma \CC.4(c)  
and Proposition \CC.5  state the following.

\Prop{\CC.7. (Decreasing)}
$$
R^{\downarrow}_F \ =\  R^{\downarrow}_q \qquad{\rm with}\ \ q\ \equiv\ q_F.
\eqno{(\CC.11a)}
$$
$$
F \ \ {\rm has\ finite\ decreasing\ characteristic\ \ } q_F
\quad\iff\quad
P_e\in \Int F,
\eqno{(\CC.11b)}
$$
which in turn holds if and only if the decreasing radial $F$-harmonics are
$$
-\den K_q(|x|)+C \qquad {\sl where}\ \  \den\geq0 \ \ {\sl and}\ \  C\in\bbr, 
\ \ {\sl and}\ \  q=q_F.
\eqno{(\CC.11c)}
$$

\medskip
\noindent
{\bf Remark.}
In summary we have that:  
\smallskip\noindent
(1) For some $p$ finite, $K_p(|x|)$ is an increasing (or downward-pointing)
$F$-harmonic on $\rn-\{0\} \iff  -P_e\notin  F\iff F$ has finite increasing characteristic.
\smallskip\noindent
(2)  For some $q$ finite, $-K_q(|x|)$ is an decreasing (or upward-pointing)
$F$-harmonic on $\rn-\{0\} \iff  P_e\in \Int F\iff F$ has finite decreasing characteristic.
\smallskip\noindent
(3) Both $K_p(|x|)$ and $-K_q(|x|)$ are $F$-harmonic on $\rn-\{0\} \iff  F$
has both characteristics $(p,q)$ finite $\iff \ -P_e\notin F$ and $P_e\in \Int F$.

These criteria hold for a significant number of  degenerate (non uniformly elliptic) subequations.
See the next section and Appendix A in Part II.)
However, in case (3) if either $F$ or $\ft$ is convex, then both are uniformly elliptic. 
Conversely, uniform ellipticity always implies that $(p,q)$ are both finite even in the non-convex case.

\medskip
Finally, combining both characteristics we have

\Prop{\CC.8}
{\sl
If $F$ has characteristics $(p,q)$, then the radial subequation for $F$ is }
$$
F\ =\ R^\uparrow_p \cup R_q^\downarrow.
\eqno{(\CC.12)}
$$

\Remark {\CC.9. (Boundary Convexity and the Riesz Characteristic)}
The finiteness of the two characteristics of $F$, which is so easy to ascertain, 
is equivalent to automatic boundary convexity for all domains.

\Prop{\CC.10}
{\sl
The boundary $\bo$ of every smoothly bounded domain $\O\ss\ss\rn$ is

\smallskip
\noindent
(a) \ \ strictly $F$-convex 
$\quad\iff\quad
p_{\ft} \ =\ q_F \ <\ \infty
\quad\iff\quad
P_e\in \Int F$,
\smallskip
\noindent
(b) \ \ strictly $\ft$-convex 
$\quad\iff\quad
p_F \ =\ q_{\ft} \ <\ \infty
\quad\iff\quad
-P_e\notin  F$,
\smallskip
\noindent
(c) \ \ both strictly $F$- and $\ft$-convex 
$\quad\iff\quad
(p_F, q_F)$ is finite
$\quad\iff\quad
P_e\in \Int F$ and $-P_e \notin F$.
}

\pf  We first prove (b).
By Lemma 5.3(ii$'$) in [\HLDD], $\bo$ is strictly $\ft$-convex at $x\in\bo$ for all domains $\O$ if and only if
$$
\forall \, B \in \Sym(W), \ \ B+tP_e \in \Int\ft \fa t\geq\ {\rm some\ } t_0.
\eqno{(\CC.13)}
$$
where $|e|=1$ and $W= e^{\perp}$.  Now (\CC.13) $\Rightarrow \ \ P_e\in \Int \ft$
$\Rightarrow {1\over t} B +P_e \in \Int \ft$ for all $t\geq $ some $t_0$
 $\Rightarrow$  (\CC.13). 
 Thus (\CC.13) is equivalent  to $p_F <\infty$ by Lemma \CC.4(c).
 The proof of (a) follows by duality, and (a) and (b) together imply (c).
  \qed

Results  in [\HLDD] immediately imply the following.

\Theorem{\CC.11. (Universal Solvability of the Dirichlet Problem)}
{\sl
Suppose that $F$ is an \ST cone subequation for which both  Riesz
characterstics $p_F$ and $q_F$ are finite (or equivalently for which 
the simple condition $P_e\in \Int F$ and $-P_e \notin F$ holds).  Then for every domain $\O\ss\ss\rn$ with
smooth boundary $\bo$, and for every $\vf\in C(\bo)$, there exists a unique
$h\in C(\ob)$ such that
\medskip

(1) \ \ $h$ is $F$-harmonic on $\O$, and
\medskip

(2)\ \ $h\bigr|_{\bo} \ =\ \vf$.
}

\Remark{\CC.12} In fact Theorem \CC.11 holds for any constant coefficient
second-order subequation $F$  if and only if  its asymptotic cone subequation
$\oa F$  satisfies $P_e\in \Int F$ and $-P_e \notin F$ for all $|e|=1$.

\bigskip
\centerline{\bf
A Description of all ST-Invariant Cone Subequations
}
\medskip

Although it is always easy to compute the characteristics $(p,q)$ of a given  $F$,
it is still enlightening to give a description (or construction) of all the possible ST-invariant cone subequations
with characteristics $(p,q)$.

The following specific examples are instrumental in this description.
 For $A\in\Symn$ let 
$\l_1(A) \leq\cdots\leq \l_n(A)$ denote the ordered eigenvalues of $A$, and 
set $\l_{\rm min}(A) \equiv \l_1(A)$ and $\l_{\rm max}(A)\equiv \l_n(A)$.  We then
define
$$
\cp^{\rm min/max}_p \ \equiv \  \left\{A : \l_{\rm min}(A)  + (p-1) \l_{\rm max}(A) \ \geq\ 0 \right \}
\eqno{(\CC.14)}
$$
$$
\cp^{\rm min/2}_p \ \equiv \  \left\{A : \l_{\rm min}(A)  + (p-1) \l_{2}(A) \ \geq\ 0 \right \}
\eqno{(\CC.15)}
$$
It is clear  that both of these are O$(n)$-invariant cone subequations.
Both $A\equiv P_{e^\perp}-(p-1) P_e$ and $B\equiv  - P_{e^\perp} + {1\over p-1} P_e$
have the property that $\l_{\rm min} +(p-1)\l_{\rm max}  =0$, which shows that 
$A,B \in  \partial \cp_p^{\rm  min/max}$ and hence $\cp_p^{\rm  min/max}$ has characteristics
$(p,q)$ where $q$ satisfies $(p-1)(q-1)=1$.
Similarly,  $\cp_p^{\rm  min/2}$ has characteristics $(p,\infty)$ if $n\geq 3$.

Our general discussion is a characterization in terms of these two examples and their duals.

\Prop{\CC.13}
{\sl
Suppose that $F$ is an ST-invariant (not necessarily convex) cone subequation.
Then $F$ has a finite (increasing) Riesz characteristic $p$ if and only if 
$$
\cp_p^{\rm min/2} \ \ss\ F\ \ss\ \cp_p^{\rm min/max}.
\eqno{(\CC.16)}
$$
Equivalently, $K_p(|x|)$ is an increasing (or downward-pointing) radial $F$-harmonic.
In particular, both the  "smallest" and the "largest" subequations, $\cp_p^{\rm min/2}$
and  $\cp_p^{\rm min/max}$, have Riesz characteristic $p$.
}
\pf
Let $A(p) \equiv P_{e^\perp}-(p-1) P_e$.
If $F$ satisfies (\CC.16), then $A(p)\in \cp_p^{\rm min/2} \Rightarrow A(p)\in F$,
and  $A(p) \notin \Int \cp_p^{\rm min/max} \Rightarrow A(p)\notin \Int F$,
which proves that $A(p)\in \partial F$, and hence $F$ has characteristic $p$.

 Each $A\in \Symn$ can be written as a sum $A = \l_1 P_{e_1} + \cdots + \l_n P_{e_n}$
using the ordered eigenvalues of $A$. Set $B_0 \equiv \l_1 P_{e_1} + \l_2 P_{e_1^\perp}$,
and $B_1 \equiv \l_1 P_{e_1} + \l_n P_{e_1^\perp}$, 
and note that $B_0\leq A \leq B_1$.

If $A\in \cp^{\rm min/2}_p$, then $\l_1 + (p-1)\l_2\geq0$. Thus, $B_0 \in \cp^{\rm min/2}_p$.
Since $\cp^{\rm min/2}_p$ and $F$ have the same increasing radial profile $E^\uparrow$
given by (\CC.1)   (and $\l_2\geq0$), we conclude that $B_0\in F$.  
However, $B_0\leq A$ proving that $A\in F$.

For the other inclusion, pick $A \in F$.  Since $F\ss \cpt$, we have $\l_{\rm max} \geq0$.
Now $A\leq B_1$ implies $B_1\in F$. Again $F$ and $\cp^{\rm min/max}_p$ 
have the same increasing radial profile $E^\uparrow$ given by (\CC.1).  Therefore,
$B_1 \in \cp^{\rm min/max}_p$.  This implies by definition that $A \in \cp^{\rm min/max}_p$.\qed
\medskip

This imposes a constraint on the decreasing characteristic $q$ of $F$.

\Cor{\CC.14}
{\sl
The characteristics of $F$ satisfy 
$$
(p-1)(q-1)\ \geq\ 1.
\eqno{(\CC.17)}
$$
}
\pf
It follows from  Definition (\CC.9a) that  if one shrinks a subequation, then its decreasing characteristic 
goes up.   Thus if $F$ has characteristic $p$, we have $\cp_p^{\rm min/max} \supset F$
and so the decreasing characteristic $q$ of $F$  satisfies $q-1 \geq q_{ \cp_p^{\rm min/max}} -1 = 1/(p-1)$.
\qed\medskip

\noindent
{\bf Remark.}   The only ST-invariant cone subequation with given characteristics $(p,q)$ satisfying
(\CC.17) is $ \cp^{\rm min/max}_p$.   This can be proved by using Proposition \CC.15 below,
but details are omitted here.

\medskip

It is just as easy to describe all examples with dual characteristic $q$. 
 First note that the duals of the two subequations in (\CC.16) are given by
 $$
 \cpt_p^{\rm min/2} : \l_{\rm max}(A) + (p-1) \l_{n-1}(A)\ \geq\ 0, 
 \eqno{(\CC.18)}
$$
 $$
 \cpt_p^{\rm min/max} : \l_{\rm max}(A) + (p-1) \l_{\rm min}(A)\ \geq\ 0.
 \eqno{(\CC.19)}
$$
Note that the increasing characteristics of these two subequations are both $\infty$,
and the decreasing characteristics are $p$ by (\CC.9c).

Applying Proposition \CC.13 to $\ft$ now yields the following result.

\Prop{\CC.15}
{\sl
Suppose that $F$ is an ST-invariant (not necessarily convex) cone subequation.
Then $F$ has a finite (decreasing) Riesz characteristic $q$ if and only if }
$$
\cpt_q^{\rm min/max} \ \ss\ F\ \ss\ \cpt_q^{\rm min/2}.
\eqno{(\CC.20)}
$$
\pf
$$
\cp_q^{\rm min/2} \ \ss\ \ft \ \ss\ \cp_q^{\rm min/max}
\quad\iff\quad
\cpt_q^{\rm min/max} \ \ss\ F\ \ss\ \cpt_q^{\rm min/2}. \qquad\mathqed
$$
\medskip

Finally, it is possible to describe all the ST-invariant cone subequations with both characteristics finite.

\Prop{\CC.16}
{\sl 
Suppose that $F$ is an ST-invariant  cone subequation.
Then $F$ has both Riesz characteristics $(p,q)$ finite if and only if
$$
\cp_p^{\rm min/2} \cup \cpt_q^{\rm min/max} 
\ \ss\ F\ \ss\ 
\cp_p^{\rm min/max} \cap \cpt_q^{\rm min/2}.
\eqno{(\CC.21)}
$$
Such subequations exist if and only if 
$$
(p-1)(q-1) \ \geq\ 1,
\eqno{(\CC.22)}
$$
and so in particular if this constraint holds for $(p,q)$, then both}
$$
\cp_p^{\rm min/2} \cup \cpt_q^{\rm min/max} 
\ \ {\sl and}\ \ 
\cp_p^{\rm min/max} \cap \cpt_q^{\rm min/2} \ \ 
{\sl have\ characteristics}\ \ (p,q).
\eqno{(\CC.23)}
$$

\pf 
Note that (\CC.21) holds if and only if both (\CC.16) and (\CC.20) hold.  Thus 
by Propositions \CC.13 and \CC.15, $F$ has finite Riesz characteristics $(p,q)$ if and only if 
(\CC.21) holds.

Corollary \CC.14 states that if $F$ has characteristics $(p,q)$, then (\CC.22) must hold.
Now suppose that (\CC.22) holds.  Then 
$$
\cpt_q^{\rm min/max} \ \ss\ \cp_p^{\rm min/max}
\and
\cp_p^{\rm min/2} \ \ss\ \cpt_q^{\rm min/2}
\eqno{(\CC.24)}
$$
because $\l_{\rm max} + (q-1) \l_{\rm min} \geq0 \ \Rightarrow
\l_{\rm min} + (p-1) \l_{\rm max} \geq0$ if $p-1 \geq 1/(q-1)$; and
$\l_{\rm min} + (p-1) \l_{\rm 2} \geq0 \ \Rightarrow\ 
\l_{n-1} + (p-1) \l_{\rm max} \geq0\ \Rightarrow\ 
\l_{\rm max} + (q-1) \l_{n-1} \geq0$ if $q-1 \geq 1/(p-1)$.
Finally, (\CC.24) implies that 
$ \cp_p^{\rm min/2} \cup  \cpt_q^{\rm min/max} \ss 
 \cp_p^{\rm min/max} \cap  \cpt_q^{\rm min/2}$
 so that both of these subequations have characteristics $(p,q)$.\qed


\vfill\eject

\centerline{\headfont \HH.   Some Illustrative Examples.}
\medskip 

For the basic subequations the Riesz characteristic is quite easy to compute.
We shall illustrate this with a selection of examples of differing types. A  detailed
discussion of subequations of characteristic $p$, and further results, are given in Appendix A of Part II.

For  $A\in \Symn$ we let
$$
\l_1(A)\ \leq\ \l_2(A)\ \leq\ \cdots \ \leq\ \l_n(A)
\eqno{(\HH.1)}
$$
denote the {\bf ordered eigenvalues} of $A$.

\Ex{\HH.1. (The p-Convexity Equation)}  For each real number $p$ with $1\leq p\leq n$,
the smallest (see Lemma A.2 in Part II) convex  cone subequation with characteristic $p$ 
is also one of the most basic:
$$
\cp_p\ \equiv\ \{ A : \l_1(A)+\cdots + \l_{[p]} (A) + (p-[p]) \l_{[p]+1}(A) \ \geq\ 0\}.
\eqno{(\HH.2)}
$$
For $p$ an integer
the $\cp_p$-subharmonic functions are characterized by the fact that their restrictions
to minimal submanifolds of dimension $p$ are intrinsically  subharmonic .
For this and a discussion of
the geometry associated with this equation, see [\HLPCON].
(Results for integer $p$ go back to H. Wu [Wu], [Sh].)
Note, by the way, that $\cp_1=\cp$ is the homogeneous Monge-Amp\`ere subequation
and $\cp_n=\D$ is the standard Laplacian.

There are complex and quaternionic analogues $\cp_p^\bbc$ and $\cp_p^\bbh$
defined by (\HH.2) but using the eigenvalues of the complex (respectively quaternionic)
hermitian symmetric part of $A=D^2u$.
When $p=1$ this yields the homogeneous complex and quaternionic Monge-Amp\`ere subequations. 
The $\cp_p^\bbc$-subharmonic functions are characterized by the fact that their  restrictions to complex
$p$-dimensional submanifolds are $\D$-subharmonic.
 The Riesz characteristics of  $\cp_p^\bbc$ and $\cp_p^\bbh$ are $2p$ and $4p$ respectively. See Lemma \HH.8 below.

\Ex{\HH.2. (The Elementary Symmetric  or Hessian Equations)} For each integer $k$, $1\leq k\leq n$, let
$\s_k(A)$ denote the $k^{\rm th}$ elementary symmetric function of the eigenvalues of
$A\in\Symn$.  The convex cone subequation
$$
\Sigma_k\ =\ \{ A : \s_1(A)\geq 0,  \s_2(A)\geq 0,  ... , \s_k(A)\geq 0\}
\eqno{(\HH.3)}
$$
has Riesz characteristic
$$
p_{\Sigma_k} \ \equiv \ {n\over k}
\eqno{(\HH.4)}
$$
These subequations, often called {\sl hessian equations}, 
 have been the focus of much study (e.g., [TW$_*$], [La$_*$]).
There are again complex and quaternionic analogues $\Sigma_k^\bbc$ and 
$\Sigma_k^\bbh$ with Riesz characteristics $2n/k$ and $4n/k$ respectively.

\Ex{\HH.3. (The $\d$-Uniformly Elliptic  Equations)}  The $\d$-uniformly elliptic regularization
of the basic subequation $\cp\equiv \{A\geq0\}$ (cf. Example \HH.10) is
$$
\cp\left(   \d     \right)\ \equiv \
\left    \{A : A+ \smfrac  {\d}{n}  \tr(A)  I \ \geq\ 0    \right  \}.
\eqno{(\HH.5)}
$$
These are convex cone subequations with Riesz characteristic $p=n(1+\d)/(n+\d)$.
Given $p$ with $1\leq p\leq n$ and setting
$$
\d \ =\ {n(p-1) \over n-p}
\eqno{(\HH.6)}
$$
Lemma A.2 states that $\cp(\d)$ is the largest O$(n)$-invariant convex cone subequation  with Riesz characteristic
$p$. There are again complex and quaternionic analogues described in Example  \HH.7  below.

\Ex{\HH.4.  (Geometrically Defined Subequations)}  
These important examples account for our choice of normalization in defining the Riesz characteristic.
Fix a compact subset $\GG\ss G(p,\rn)$
in the Grassmannian of $p$-planes in $\rn$, and define
$$
F(\GG) \ \equiv \ \{A : \tr_W(A)\ \geq\ 0 \ \ {\rm for\ all\ \ } W\in\GG\}\}
\eqno{(\HH.7)}
$$
where $\tr_W(A)$ denotes the trace of $A\bigr|_W$.  
Assuming the ST-invariance of $\GG$,
 the Riesz characteristic is easily seen to be
$$
p_{F(\GG)}\ =\ p.
\eqno{(\HH.8)}
$$
Many interesting subequations arise this way.  When $\GG = G(p,\rn), G^\bbc(p, \bbc^n)$
and $G^\bbh(p, \bbh^n)$ we retrieve the integer cases in   \HH.1 above.
There are many other interesting examples.
One such  is LAG $\ss G^\bbr(n, \bbc^n)$, the set of Lagrangian
$n$-planes in $\bbc^n$. Closely related are the sets of isotropic $p$-planes,
and $p$-planes satisfying certain CR (Cauchy-Riemann) conditions.
Also of interest is SLAG $\ss$ LAG,
 the special Lagrangian planes (cf. [\HLCG]).
This latter is an example of a subequation associated to a calibration [\HLPCG]. 
 Other particularly interesting
examples come from the associative and coassociative calibrations in $\bbr^7$ and the Cayley
calibration in $\bbr^8$. All the specific subequations in this paragraph have the property that they
are \ST, i.e.,  invariant under a subgroup $G\ss {\rm O}(n)$ which acts transitively on the sphere $S^{n-1}\ss\rn$. 

These geometrically defined subequations will be the sole focus of Part II of this paper.

\medskip
\Ex{\HH.5.  (Branches of G\aa rding  Operators)} In many of the cases above,
one can associate a homogeneous polynomial operator  $\Phi(D^2 u)$. 
When the polynomial   $\Phi$ is G\aa rding
hyperbolic with respect to the identity $I$ (which is typically the case),
the equation has many branches [G], [\HLHP], [\HLHPP].

The simplest case is $\cp=\cp_1$ where the operator is $\Phi(A) = \det_\bbr(A)$. Here the branches are
given by $\{\l_k(A)\geq0\}$ (see (\HH.1)).  Unfortunately, in this case the branches for $k>1$ have
infinite characteristic. 

For the general G\aa rding polynomial $\Phi(A)$ of degree $m$, there are ordered eigenvalues,
$$
\L_1(A) \ \leq\ \L_2(A)\ \leq\ \cdots\ \leq \ \L_m(A), \qquad {\rm and} \ \ \Phi(A) \ =\ \L_1(A)\cdots \L_m(A).
\eqno{(\HH.1)'}
$$
Just as with $\det_\bbr(A)$,  the $k^{\rm th}$  branch is defined  by $\{\L_k(A)\geq 0\}$
for $k=1,...,m$.  The Riesz characteristics $p_1\leq \cdots\leq p_m$ of these respective  branches are determined by 
the eigenvalues of $P_e$ (assuming ST-invariance).  They are exactly the numbers 
$1/\L_j(P_e)$, $j=1,...,m$ arranged in increasing order (see Proposition A.10 in [\HLTangII]).
Therefore the number of branches with finite Riesz characteristic equals the number of non-zero eigenvalues of 
$P_e$.  Only the first and smallest branch is convex, and it is uniformly elliptic $\iff$ all branches are
uniformly elliptic
$\iff \Phi(P_e)>0$.

 G\aa rding operators are plentiful.   For instance, in each
of our first three examples there is an associated G\aa rding operator, and hence each comes
equipped with branches.  To illustrate, for the case where $p$ is an integer in Example \HH.1, we have
$$
\Phi(A) 
\ =\ \prod_{i_1<\cdots<i_p} (\l_{i_1}(A)+\cdots+\l_{i_p}(A))  \ =\   \det\left ( D_A : \L^p \rn \arr \L^p \rn\right).
\eqno{(\HH.9)}
$$
Said differently, $\L_I(A) = \l_{i_1}(A) + \cdots + \l_{i_p}(A)$ are the eigenvalues.
Here $D_A$ is the extension of $A$ as a derivation.  The $k^{\rm th}$ branch is  given by requiring
that the $k^{\rm th}$ ordered $p$-fold sum of the $\l_i$'s  be $\geq0$.  
One easily computes that {\sl the
first ${ {n-1} \choose {p-1}}$ branches have Riesz characteristic $p$} and the remaining branches
have infinite characteristic.

In Example \HH.2 the G\aa rding operator is $\Phi(A) = \s_k(A)$. 
Although the eigenvalues $\L_j(A)$ of $\Phi$ do  not
have an explicit formula in terms of the standard eigenvalues of $A$, 
the eigenvalues of $A=P_e$ are all zero except for one which equals $k/n$.  
Hence, $\Sigma_k$ has characteristic $n/k$ and all other branches have characteristic $\infty$.

In Example \HH.3 the eigenvalues are
$$
\L_k (A) \ =\ \l_k(A) +{\d\over n(1+\d)} \tr (A),\quad k=1,...,n.
$$
Hence, each of  the $k^{\rm th}$ branches  $\{\L_k(A)\geq 0\}$, for 
$k\geq2$, has the same Riesz characteristic $p = n(1+{1\over \d})$,  which is finite but larger than $n$,
while as noted above, the first branch $\cp(\d)$ has characteristic $n(1+\d)/(n+\d)$.

\medskip
\Ex{\HH.6.  (Trace Powers of the Hessian)} Consider the non-convex cone subequation
$$
F \ \equiv\ \left\{A :  \tr\left( A^q\right) \ \geq\ 0\right\}
$$
where $q>1$ is an odd integer. More generally one could define $A^q$ for any $q>0$ by
using the function $t^q$ for $t\geq0$ and $-|t|^q$ for $t<0$.  In all cases one computes
that the Riesz characteristic is
$$
p_F\ =\   1+ (n-1)^{{1\over q  }}
$$
More generally, for $k\in [1,n]$ and $q>0$ real numbers, there is the subequation
$$
F \ \equiv\ \{A : \l_1^q(A) + \cdots +\l_{[k]}^q(A) + (k-[k])\l_{[k]+1}^q \ \geq\ 0\}
$$
with $t^q$ defined as above.  Here the  Riesz characteristic is
$$
p_F\ =\   1+ (k-1)^{{1\over q  }}
$$

\medskip
\Ex{\HH.7.  (Complex and Quaternionic Analogues)}  Suppose $F\ss\Symn$ is an O$(n)$-invariant
subequation. Then $F$ can be defined by the constraint set $E\ss\rn$ imposed by $F$ on the eigenvalues $\l(A) = (\l_1(A),...,\l_n(A))$. Thus $A\in F \iff \l(A) \in E$.
The equation $F$ has  {\sl complex
and quaternionic analogues} $F^\bbc$ and $F^\bbh$, defined on $\bbc^n = (\bbr^{2n}, J)$ and 
$\bbh^n = (\bbr^{4n}, I,  J, K)$ respectively,   as follows.  For $A\in \Sym(\bbr^{2n})$ consider
the hermitian symmetric part
$$
A_\bbc \ \equiv\ \half (A-JAJ) 
$$
whose eigenspaces are complex lines with ordered eigenvalues 
$\l_1(A_\bbc)\leq\cdots\leq \l_n(A_\bbc)$.  One now defines $F^\bbc$ by applying
the eigenvalue constraints $E$ of $F$ to these eigenvalues of $A_\bbc$. The story in the quaternionic case is 
parallel and uses the quaternionic hermitian symmetric part
$A_\bbh \ \equiv\ {1\over 4}  (A- IAI -JAJ -KAK)$ and eigenvalues $\l_k(A_\bbh)$.

\Lemma{\HH.8} {\sl If $F$ is an O$(n)$-invariant  cone subequation with Riesz characteristic $p$, then
the Riesz characteristics of $F^\bbc$ and $F^\bbh$ are}
$$
p_{F^\bbc} \ =\ 2p \and  p_{F^\bbh} \ =\ 4p.
$$
\pf
We consider the complex case.  If  $A = P_{e^\perp} - (p-1) P_e  \in \Sym(\bbr^{2n})$,
then one computes  that 
$$
A_\bbc \ =\ P_{\bbc e^\perp}  -\left({p\over 2} -1\right) P_{\bbc e}
\and
A_\bbh \ =\ P_{\bbh  e^\perp}  -\left({p\over 4} -1\right) P_{\bbh  e}
\eqno{(\HH.10)}
$$
which displays the eigenvalues of $A_\bbc$ and $A_\bbh$.
\qed

\medskip
\Ex{\HH.9.  (The Subequation Determined by a  G\aa rding Operator and a Universal Eigenvalue Constraint)}  The procedures above can be greatly generalized. Note, to begin, that 
given an O$(m)$-invariant subequation $F$, the eigenvalue set $E\equiv \l(F)$ is closed, invariant under permutation of coordinates and $\bbr^m_+$-monotone.  Conversely, any such eigenvalue set $E$
determines an O$(m)$-invariant subequation $F=\l^{-1}(E)$.
Each such $E$ is a {\sl universal eigenvalue subequation} in the sense that, for each degree-$m$
G\aa rding operator $\Phi$ on $\Symn$, the set  $F \equiv \l_{\Phi}^{-1}(E)$ is a subequation on $\rn$,
where $\l_\Phi: \Symn \to \bbr^m$ is the eigenvalue map associated to $\Phi$.  See Proposition A.8 
in Appendix A of [\HLTangII] for the details and further discussion.


\medskip
\Ex{\HH.10.  (The $\d$-Uniformly Elliptic Regularization of a Subequation)} Given a cone
subequation $F\ss\Symn$ and $\d>0$, define
$$
F(\d) \ \equiv\ \left\{ A : A+\smfrac \d n \tr(A) I \ \in\ F\right\}.
\eqno{(\HH.11)}
$$
This equation satisfies the uniformly elliptic condition:
$$
F(\d) + \cp(\d) \ \ss\ F(\d).
\eqno{(\HH.12)}
$$
One computes that 
$$
F\ \ {\rm has\ Riesz\ characteristic\ \ } p
\qquad\iff  \qquad
F(\d)\ \ {\rm has\ Riesz\ characteristic\ \ } {  pn(1+\d)  \over   n+\d p  }
$$

\vskip .3in
 

\vfill\eject

\centerline{\headfont \LL.    K$_p$-Convexity and Monotonicity.  } \medskip

In this section we give a fairly complete discussion of the classical one-variable results that underlie this paper.
They concern the properties of  subsolutions to the one-variable subequation $R_p$ introduced below.

In the  following section we will prove that associated to each $F$-subharmonic function $u$
there are three functions of $r$ (denoted $M(r), S(r)$ and $V(r)$), which are subsolutions 
of $R_p$, and which capture  much of the asymptotic behavior of $u$.
By Lemma \LL.1(3) below this will imply the key double monotonicity result,  Theorem \BB.4, which 
is needed for defining the notion of density and for proving our main theorems.

Fix a real number $p$ with $1\leq p <\infty$, and for $r>0$
consider the {\bf one-variable Riesz kernel} 
$$
K_p(r)\ \equiv\ {1\over (2-p)} r^{2-p} \ \ {\rm  if} \ \ p\neq 2\ \ {\rm and}\ \ K_2(r)\ =\ \log\, r.
\eqno{(\LL.1)}
$$
With this normalization 
$$
K_p'(r) \ =\ {1\over r^{p-1}} \qquad{\rm for\ all}\ \ 1\leq p  < \infty.
$$
Note that $K_p(r)$ is a strictly increasing solution to the subequation
$$
 R_p \ : \ \psi''(r)+{p-1\over r} \psi'(r) \ \geq \ 0\ \ \ \  {\rm on}\ (0,\infty).
\eqno{(\LL.2a)}
$$
Alternatively,
$$
 R_p \ : \  {d\over dr} \biggl( r^{p-1} \psi'(r)  \biggr)
 \ =\   {d\over dr} \left(  { \psi'(r) \over K_p'(r)} \right)
 \ \geq \ 0\ \ \ \  {\rm on}\ (0,\infty).
\eqno{(\LL.2b)}
$$
All solutions of $R_p$ are of the form
$$
h(r) \ \equiv\ C K_p(r) + k \ \ {\rm with} \ \ C, k\in \bbr.  \qquad {\bf  (Riesz\ Harmonics)}.
\eqno{(\LL.3)}
$$
Note that  $h(r)$ is increasing if and only if $C\geq 0$.

The change of variables
$$
s\ =\ K_p(r) \ \ {\rm  along\ with\ its\ inverse}\ \ r\ =\ K_p^{-1}(s)
\eqno{(\LL.4)}
$$
play an important role.  The inverse $r(s) =  K_p^{-1}(s)$ is defined on the range of $K_p$ which is
the interval $(0,\infty)$ when $1\leq p<2$,  all of $\bbr$ when $p=2$, 
and $(-\infty,0)$ for $2<p<\infty$.

\Lemma {\LL.1. (The Equivalences)}
{\sl 
 The following conditions on an upper semi-continuous
  function $\psi(r)$, defined on a subinterval  of $(0,\infty)$, are equivalent.

\medskip

(1) \ {\bf ($R_p$-Subharmonic) } \ \ $\psi(r)$ satisfies the subequation $R_p$ defined by (\LL.2).

\medskip

(2) \ ({\bf $K_p$-Convexity)} \ \ $\psi(r)$ is $K_p$-convex, meaning that under the change of variables 

\qquad (\LL.4), the function $f(s) \equiv \psi(r(s))$   is a convex function of $s$.

\medskip

(3) \ {\bf ($K_p$-Monotonicity) } \ \ 
$$
{\psi(r) - \psi(t) \over K_p(r) - K_p(t)}\ \ {\rm is\ non-decreasing\ in\ } r \ {\rm and} \ t\  (r\neq t).
$$

\medskip

(4) \ {\bf ($R_p$-Comparison) } \ \ If $\psi(r) \ \leq CK_p(r) + k$ for $r=s$ and $r=t$,
then the

\qquad  inequality holds for all $r$ between $s$ and $t$.
}

\pf
 Now  $f(s) \equiv \Aa(r(s))$ implies $\Aa(r) = f(K_p(r))$.
 First, assume $\psi$ is smooth.   Then
  $$
  \Aa'(r) = f'(s) {K_p'(r)}  \quad {\rm  and\  hence\ } \quad
  {d\over dr} \left(  { \psi'(r) \over K_p'(r)} \right) \ =\ f''(s)K_p'(r).
  \eqno{(\LL.5)}
$$ 
For general $\psi$, the fact that: $\vf(r)$ is a test function for $\psi$ at $r_0$ if and only if 
$\vf(r(s))$ is a test function for $\psi(r(s))$ at $s_0\equiv s(r_0)$, reduces the proof to the smooth case.
Since viscosity convexity $f''(s)\geq0$  is equivalent to classical convexity (see for example,  [\HLDD, Prop 2.6]), this proves that (1) $\iff$ (2).

Now (3) is just monotonicity of the slopes of secant lines to the function $f(s)\equiv \psi(r(s))$,
and hence it is equivalent to the convexity of $f(s)$.
Assertion (4) is just the statement that $f(s)$ is convex if and only if $f$  satisfies comparison 
with affine functions $Cs+k$.\qed
\medskip

\Cor{\LL.2} {\sl Let $\psi(r)$  satisfy the equivalent conditions in Lemma \LL.1.  Then
\medskip
(a) \ \  The function $\psi(r)$ is locally Lipschitz continuous.

\medskip
(b) \ \ The left and right hand derivatives $\psi_{\pm}'(r)$ exist.
}
\pf
The corresponding statements for the function $f(s)\equiv \psi(r(s))$ with 
$r(s) = K^{-1}(s)$ are standard classical facts about  the convex function $f$.\qed

\vskip .3in

\centerline
{
\bf Densities.
}
\medskip

The remainder of the appendix is devoted to describing properties of a function $\psi(r)$,
 defined on an interval $(0, r_0)$ (with $r_0=\infty$ possible), under  the 

\medskip
\noindent
{\bf Hypothesis:}  \ \ $\psi$ satisfies the equivalent conditions (1) -- (4) of 
Lemma \LL.1 

\medskip

The    Properties (1) and/or (3)   enable us to introduce the following.

\Cor{\LL.3 (Existence of   Densities)}  {\sl  The decreasing limits
$$
\den^\psi \ =\ \lim_{\eqalign{&r,t\to0\cr &t>r>0}}  {\Aa(t) -\Aa(r) \over K(t) -K(r)} 
\ =\ \lim_{r\to 0}  { \psi_\pm'(r)       \over    K'(r)     }
\eqno{(\LL.6)}
$$
exist and define the {\bf density} $\den^\psi$.  Moreover, if $\psi$ is increasing, then  $0\leq \den^\psi <\infty$.  
}
\pf
To see that the two decreasing limits in (\LL.6) agree divide the numerator and denominator of
${\psi(r+\d) - \psi(r)  \over K(r+\d) - K(r)}$ by $\delta$ and let $\d\to0$.
\qed

\medskip
In this one-variable context, rather than in its later applications, it might be better to call this 
``the derivative of $\psi(r(s))$ at $r=0$''.

Note that the monotonicity quotient in (3) remains unchanged if $\psi$ is replaced by
a translate $\psi-c$ with $c\in\bbr$.  In particular, the densities of $u$ and $u-c$ are the same.
This point is critical in establishing the following.

\Lemma{\LL.4}
{\sl  If $\psi$ is increasing, then there  exists $c\in \bbr$ and $r_0$ such that 
$$
{\psi(r) - c \over K(r)} \ \ \ {\rm decreases\ to\ \ \ } \den^\psi
\ \ {\rm as}\ \ 0< r< r_0\ \ {\rm decreases\ down\ to\ \ } 0.
$$
Moreover, 
$$
 {\psi(r) -\psi(0)  \over K(r)} \quad {\rm decreases\ to \ \ } \den^\psi\ \ {\rm if}\ \ 1\leq p<2,\ \ {\rm and}
 \eqno{(\LL.7a)}
$$
$$
\lim_{r\to0} {\psi(r) \over K(r)} \ =\ \den^\psi\quad {\rm if\ \ } 2\leq p<\infty.
\eqno{(\LL.7b)}
$$
(Note: if we set $\psi(0)=0$ when $1\leq p<2$, then we have  in all cases ($1\leq p<\infty$) that $\lim_{r\to0} \psi(r)/K(r) = \den^\psi$.)}

\pf
For any value of $p$, $1\leq p<\infty$, there is exactly one point in$[0,\infty]$
where $K$ vanishes.   However there are three cases:
$K(0)=0$ if $1\leq p<2$, $K(1) = 0$  if $p=2$, and $K(\infty)=0$ if $2<p$.
First  let us suppose that the function $\psi$ is defined and finite on an interval
containing the point where $K$ vanishes.
Then one can take $r$  in (\LL.6)  to be that point,  and the Proposition follows immediately from
the double monotonicity in (3). Specifically, in the three cases we obtain that:
$$
(1\leq p<2)\qquad
{\psi(t) - \psi(0) \over K(r)}  \ \ \ {\rm decreases\ to\ \ \ } \T^\psi \ \ {\rm as\ \ } 
0<t< r_0\ \ {\rm decreases\ to\ \ }0,
\eqno{(\LL.7a)'}
$$
$$
(p=2)\qquad
{\psi(t) - \psi(1) \over K(t)}  \ \ \ {\rm decreases\ to\ \ \ } \T^\psi \ \ {\rm as\ \ } 
0<t< r_0\ \ {\rm decreases\ to\ \ }0,
\eqno{(\LL.7b)'}
$$
$$
(2<p)\qquad
{\psi(t) - \psi(\infty) \over K(t)}  \ \ \ {\rm decreases\ to\ \ \ } \T^\psi \ \ {\rm as\ \ } 
0<t< r_0\ \ {\rm decreases\ to\ \ }0,
\eqno{(\LL.7c)'}
$$
This leaves us with an extension problem in the last two cases.   Namely we must prove that there  exists an $r_0>0$ such that 
the restriction of the given $\psi$ to $(0, r_0)$
$$
(p=2)\qquad 
{\rm has\  an\  extension\ \ }   \overline \psi  \ \ {\rm to}\ \ (0,1] \ \ {\rm satisfying\  Lemma \ \LL.1(3),\  and} \qquad
\eqno{(\LL.8b)}
$$
$$
\eqalign
{
(2<p)\qquad
&{\rm has\  an\  extension\ \ }   \overline \psi  \ \ {\rm to}\ \ (0,\infty) \ \  {\rm satisfying\  Lemma \ \LL.1(3)} 
\qquad\qquad
\cr
&{\rm   with}\ \ \psi(\infty) <\infty.
}
\eqno{(\LL.8c)}
$$

Suppose that $\psi$ has domain containing $(0,r_0]$ 
(and if $p=2$ that $r_0<1$ since if $r_0\geq1$ in this case we are finished.)
Make the change of variables $s_0 \equiv K^{-1}(r_0)$  as in Lemma \LL.1(2).
Since $2\leq p<\infty$, we have $s_0<0$. The convex increasing function $f(s)$ on 
$(-\infty, s_0)$ can be extended to a convex increasing function $\overline f$ on $(-\infty, 0]$
 by defining $\overline f$ to be the affine function
 $$
 a(s)\ \equiv\ f_{-}'(s_0)(s-s_0) + f(s_0)
 \eqno{(\LL.9)}
 $$
on $s_0\leq s\leq 0$.  Since the graph of $a(s)$ is a supporting line for
the epigraph of $f$ over $(-\infty, s_0)$, this extension $\overline f$ is convex
and increasing on $(-\infty, 0]$. 

Observe now that by translating our original $\psi$ by a suitable additive constant,
we can insure that $\overline f < 0$  on $(-\infty, 0]$. 
Now set  $\overline \psi(t) \equiv \overline f(K(t))$,
where $0<t\leq 1$ if $p=2$ and $0<t\leq  \infty$ if $2<p$.
Finally, by (\LL.7b) where $\psi(1)$ is finite,  and (\LL.7c) where $\psi(\infty)$ is finite,
the fact that $K(0)=-\infty$ implies that $\lim_{r\to0} \psi(r)/K(r) =\den^\psi$.
\qed

\Remark{\LL.5}  The subequation $R_p : \psi''+ {p-1\over r} \psi'(r)\geq0$
is linear and could have been interpreted in the distributional sense as well as the viscosity sense.


\vfill\eject

\centerline{\headfont  \BB.  Monotonicity and Stability of Averages for F-Subharmonic Functions.}
\medskip

In this section  we discuss three of the basic ways of taking an average of an $F$-subharmonic 
function, and show that each average produces  a radial $F$-subharmonic.  Since the radial $F$-subharmonics are just one-variable $R_p$-subharmonics (Proposition  \CC.3), they are well understood and  
 enjoy all the properties of Lemma \LL.1.  In particular,  they satisfy the double monotonicity described in
 Theorem \BB.4 below, which provides the vehicle for defining the densities explored in the next section.
 Finally, the stability of these averages under the tangential flow is established in Lemma \BB.5.

We assume as always
  that the subequation $F$ is an \ST cone with invariance group $G\ss {\rm O}(n)$.
We further assume that the Riesz characteristic $p$ of $F$ is finite.  This is because when $p=\infty$,
the increasing radial subequation $R^{\uparrow}_F$ is simply $g'(t)\geq0$
(Proposition \CC.3 and Definition \CC.1).
Thus, when $p=\infty$, all increasing functions $g(t)$ determine increasing 
radial subharmonics $g(|x|)$, and no sensible notion of density is possible.

To begin we set some
notation.  Let $B_r(x_0) = \{x : |x-x_0|\leq r\}$ denote the ball of radius $r$ about $x_0$,
and set $S_r(x_0) \equiv \partial  B_r(x_0)$.  
Let $A(a,b; x_0) \equiv \{x  : a<|x-x_0|<b\}$ denote an annular region centered at $x_0$.
Here and elsewhere, when $x_0=0$, reference to it will be dropped from the notation.
Thus, $B_r = B_r(0)$ and $S_r =\partial B_r$.  Similarly we set $B=B_1$
and $S=\partial B$.
 
The first average only requires that  $F$ be an \ST cone (not necessarily convex).
We denote the {\sl (spherical) maximum} for an $F$-subharmonic function $u$ defined on a region containing 
$S_r(x_0)$ by 
$$
M(u, x_0; r) \ \equiv \  \sup_S u(x_0+rx),
\eqno{(\BB.1a)}
$$
Note that if $u$ is $F$-subharmonic on $B_R(x_0)$,
then by the maximum principle
$$
\eqalign
{
&M(u, x_0; r)   =  \sup_B u(x_0+rx) =  \sup_{B_r(x_0)} u  \cr
&{\rm and\ hence\ is\ increasing\ for\ \ } 0\leq r\leq R.}
\eqno{(\BB.1b)}
$$

 By the ST-invariance of $F$
$$
M(u, x_0; |x|) \ \equiv \ \sup_{ g\in G} u(x_0+gx).
\eqno{(\BB.2)}
$$

We now simplify by setting $x_0=0$ and using the abbreviated notation $M(r) \equiv M(u;r) = M(u,0;r)$ when the meaning is obvious.

\Lemma{\BB.1} {\sl
If $u$ is $F$-subharmonic on an annular region $A(a,b)$, then $M(|x|)$ is a  radial
$F$-subharmonic function on $A(a,b)$.  
If $u$ is $F$-subharmonic on $B_R$, then  $M(r)$ is also increasing in $r$.
}
\pf   Let $u_g(x)  \equiv  u(gx)$ with $g\in G$.
 Then
$
M(|x|)\ =\ \sup_{g\in G}  u_g(x).
$
Since $F$ is  $G$-invariant, each $u_g$ is $F$-subharmonic. 
Therefore, by the standard ``families locally bounded above'' property for $F$,
it suffices to prove that
\medskip

\centerline { $M(t)$ is upper semi-continuous.}
\medskip
\noindent
This is done as follows.
For each $\d>0$,
$
N_\d \ \equiv\ \{x : u(x) < M(t)+\d\}
$
is an open set containing $\partial B_t$.  Hence the annulus $A(t-\e, t+\e)$
is contained in $N_\d$ for $\e>0$ small.  Thus $M(r) < M(t) +\d$ if  $t-\e\leq  r\leq t+\e$,
proving  that $M(t)$ is upper semi-continuous, and hence $M(|x|)$ is $F$-subharmonic.\qed

\medskip

For the other averages we make the further standing assumption that  $F$ is convex.
In this case we note the following.
$$
F\ {\rm is\ an\ \ST \ convex\ cone} \qquad\Rightarrow\qquad F \ \ss\ \D\ \equiv\ \{\tr A=0\}
\eqno{(\BB.3)}
$$
\pf
 If $F\cap \{\tr\,A =c<0\}$ is non-empty, then invariance plus convexity implies  that
$ - {c\over n}I \in F$.  Now by the  cone property, $-\l I\in F$ for all $\l>0$.
This along with positivity implies that $F=\Symn$. Since $\tr(P_{e^\perp} - (p-1)P_e) = n-p$,
 the condition $F\ss \{\tr\,A \geq0\}$ implies $p_F \leq n$.
 Therefore,
$$
F\ {\rm is\ an\ \ST \ convex\ cone} \qquad\Rightarrow\qquad 1\ \leq\ p_F\ \leq\ n.
\eqno{(\BB.4)}
$$

We now define the {\sl  spherical } and {\sl volume averages} of $u$ at $x_0$  by
$$
S(u, x_0; r) \ \equiv \  {1\over |S|} \int_{\s\in S}  u(x_0+r\s)\, d\s \ \equiv\ \intave S  u(x_0+r\s)\, d\s,
\eqno{(\BB.5a)}
$$
$$
V(u, x_0; r) \ \equiv \  {1\over |B|} \int_{\ x\in B}  u(x_0+rx)\, dx \ \equiv\  \intave B   u(x_0+rx)\, dx.
\eqno{(\BB.5b)}
$$

Note that for any upper semi-continuous function $u$, each of these functions is jointly
 upper semi-continuous in $(x_0, r)$ since $u(x_0+rx)$ is the infimum of $\vf(x_0+rx)$
 taken over continuous functions $\vf\geq u$.

\Lemma{\BB.2} {\sl
Suppose that $u$ is $F$-subharmonic on the  annulus $A(a,b)$.
Then   $S(u; |x|)$ is a radial $F$-subharmonic on $A(a,b)$.
If $u$ is $F$-subharmonic on the ball $B_R$, then  both $S(u; |x|)$ and  
$V(u; |x|)$ are   increasing radial $F$-subharmonic
functions on $B_R$  (with    limiting values $S(u,0) = V(u,0) =u(0)$ at $x=0$).
}
\pf
As noted above $S(u;r)$ and $V(u;r)$ are upper semi-continuous in $r$, and hence
so are the functions  $S(u;|x|)$ and $V(u;|x|)$ of $x$ defined on $B_R$.
The statement about their limiting values at $x=0$ is a standard fact  about 
$\D$-subharmonic functions.
It remains to  show that  $S(u;|x|)$ and $V(u;|x|)$ are $F$-subharmonic on $B_R$.
  Note that
$$
S(|x|) \ = \  \int_G u(gx)\, dg
\eqno{(\BB.6)}
$$
for a suitably  normalized invariant measure $dg$ on $G$, and that
$$
V(|x|) \ = \   n\int_0^1 S(\rho |x|) \rho^{n-1}\,d\rho
\qquad{\rm since}\qquad |B|\ =\ {1\over n} |S|.
\eqno{(\BB.7)}
$$

To prove (\BB.7), set $|x| = r$ and compute $V(r) = {1\over |B|} \int_B u(ry)\,dy$ using polar coordinates.
Now since $F$ is a convex cone subequation, averages such as (\BB.6) and 
(\BB.7) preserve $F$-subharmonicity.  This is explained further by Theorem \FF.5 and   Remark \FF.6.
\qed

\Remark {\BB.3}
By Theorem \AA.4 and  Theorem \AA.7 the lemmas above could have been restated by concluding that the 
functions $M(r)$, $S(r)$ and $V(r)$ are $R^{\uparrow}_F$-subharmonic on $(0,R)$,
or $R_p$-subharmonic on $(a,b)$ in the annular cases.

The properties of upper
semi-continuous functions $\psi(r)$ satisfying $R_p$ have been presented in detail
in Section \LL.  We make full use of those results by applying them to the three functions 
$M(u, x_0, r), S(u, x_0, r)$ and $V(u, x_0, r)$, where $u$ is an $F$-subharmonic function.
This includes the $K_p$-convexity, the $K_p$-monotonicity and the $R_p$-comparison properties 
of Lemma \LL.1.
\medskip

In particular, the $K_p$-monotonicity, part (3) of Lemma \LL.1, gives the following
basic result.

\Theorem {\BB.4. (Double Monotonicity)} {\sl Let $u$ be $F$-subharmonic in an annular region about
the origin in $\rn$. Then 
$$
{M(u,r)-M(u,s)  \over K(r)-K(s)} \quad {\rm is\ increasing\ in} \ \ r\ \ {\rm and}\ \ s.
\eqno{(\AAA.4)}
$$
for all $0<s<r$ where $M$ is defined.

 Furthermore, if $F$ is convex, the same statement holds with $M(u,r)$ replaced by 
$S(u,r)$; or by $V(u,r)$  provided that $u$ is $F$-subharmonic on a ball about the origin.}
\medskip

It is an important fact that each of these averages is stable under limits in $L^1$.
This basic  classical fact can be found in [Ho$_2$, Sec III.3.2]. We state it here in slightly
different form needed later for tangents.

\Lemma{\BB.5. (Stability of Averages)} {\sl Suppose $u_j$ is a sequence of 
$\D$-subharmonic functions on $B_R$ converging in $L^1(B_R)$ to a $\D$-subharmonic function $U$.  
Then for $0<r<R$,
\medskip

(1) \ \ \ $M(U,r) = \lim_{j\to\infty} M(u_j,r)$,

\medskip

(2) \ \ \ $S(U,r) = \lim_{j\to\infty} S(u_j,r)$,

\medskip

(3) \ \ \ $V(U,r) = \lim_{j\to\infty} V(u_j,r)$,
}

\pf Taking $K\equiv B_r$ in (3.2.7) of  Theorem 3.2.1 in  [Ho$_2$] gives us that
$$
\limsup_{j\to\infty} M(u_j, r) \ \leq\ M(u,r).
$$
Suppose there exists  $C < M(u,r)$ such that  $M(u_j,r) \leq C$ for all 
$j$ sufficiently large. Then in the $L^1$-limit  we would have $u-C\leq0$ a.e. on $B_r$. However, for
$\D$-subharmonic functions, this implies that $u-C\leq0$ everywhere on $B_r$, contrary to the
definition of $M(u,r)$.  We conclude that $\limsup_{j\to\infty} M(u_j, r) =M(u,r)$. The fact that this is
also true for all subsequences proves (1).  

As discussed in the paragraph prior to Proposition 3.2.14 in [Ho$_2$], the Theorem
3.2.13 can be applied to spherical measure $\s_r$ on $\partial B_r$.
Thus $u_j \s_r$ converges to $U\s_r$ in the weak topology of measures, yielding (2).
Finally, (3) is implied directly  by the hypothesis of $L^1(B_r)$-convergence.  \qed


 \vskip .3in

\centerline{\headfont \KK.   Densities for $F$-subharmonic Functions -- Upper Semi-Continuity.}
\medskip

From the results of the last section and  Corollary \LL.3 we have three densities, 
$$
\den^M(u,x)  \quad \den^S(u,x) \and \den^V(u,x)
$$
associated to an $F$-subharmonic function $u$ defined in a neighborhood of the origin.
For  the second two densities, we must assume that $F$ is convex.  Under this convexity assumption
there exists a fourth , even more classical density.
\bigskip

\centerline{\bf The Mass Density}
\medskip

Note that by (\BB.3) $u$ is classically $\D$-subharmonic.
Thus $\D u$ is a measure $\mu\geq0$, which means $\D u$ has a ``mass density''.
 Given a measure $\mu\geq0$ defined
in a neighborhood of a point $x_0\in \rn$, and $0< k \leq n$,  the limit
$$
\den^k(\mu, x_0) \ \equiv\ \lim_{r \downarrow 0} {\mu \left (B_r(x_0)\right)\over \a(k) r^k} ,
\eqno{(\KK.1)}
$$
if it exists, is called the {\sl $k$-dimensional mass density of $\mu$ at $x_0$.}  (See, for example,
[F, 2.10.19] for discussion and   definition of the constants $\a(k)$.)   When $k$ is an integer,
$\a(k) = |B^k|$, the volume of the unit ball in $\bbr^k$.  Suppose $\den^k(\mu, x)$ exists everywhere
or replace $\lim$ by $\limsup$ in (\KK.1). Fix an open set  $X$, a constant $c>0$,
and define $E_c \equiv \{x \in X:  \den^k(\mu, x)\geq c\}$. Then the Hausdorff $k$-measure satisfies
(cf. [Si, page 11])
$$
c \ch^k(E_c) \ \leq\ \mu(X).
$$

\bigskip
\centerline{\bf   Comparing Densities}
\medskip

The next proposition states that: {\sl All densities but $\den^M$ ``agree''}, where ``agree'' 
means ``are equal up to universal factors''.

\Prop{\KK.1}  {\sl
Suppose that $u$ is $F$-subharmonic near $x_0$ where $F$ is convex with characteristic $p$, and set $\mu=\D u$.
  Then when  $p\neq 2$, 
$$
\den^S(u, x_0) \ =\ {n-p+2\over n} \den^V(u,x_0) \ =\ {\a(n-p)\over n|p-2|\a(n)} \den^{n-p}(\mu,x_0),
\eqno{(\KK.2)}
$$
and when $p=2$ we have that}
$$
\den^S(u, x_0) \ =\ \den^V(u, x_0) \ =\ {\a(n-2)\over n\a(n)} \den^{n-2}(\mu,x_0).
\eqno{(\KK.3)} 
$$

The discussion of all densities is completed by showing that the maximum density and the
spherical density are in general ``comparable'', and in fact equal when $p=2$.

\Prop{\KK.2} {\sl
Suppose that $u$ is $F$-subharmonic near $x_0$ where $F$ is convex 
and of characteristic $p$.
Then there exists a constant $C=C(p,n) > 1$ such that}
$$
\den^M(u,x_0) \ \leq\ \den^S(u,x_0) \ \leq\ C\den^M(u,x_0) \qquad {\rm if}\ \ 2<p<\infty, \ \ {\rm and} 
\eqno{(\KK.4)}
$$ 
$$
\den^S(u,x_0) \ \leq\ \den^M(u,x_0) \ \leq\ C\den^S(u,x_0) \ \qquad {\rm if}\ \ 1<p<2, \ \ {\rm while} 
\eqno{(\KK.5)}
$$ 
$$
\quad
\den^M(u,x_0) \ = \ \den^S(u,x_0) \ \qquad\quad {\rm if}\ \ p=2.
\eqno{(\KK.6)}
$$ 

\Remark {\KK.3} 
Kiselman  proved  the   equality in (\KK.6) in the plurisubharmonic case where
 $F=\cp^\bbc$ on $\bbc^n$ (see page 161, line 2 ff. in [K$_1$])
by using Harnak's Inequality for $\D$-subharmonic functions.
The same proof works for any convex $F$ of characteristic $p=2$.
Note that for $p=1$ the left inequality in (\KK.5) holds but the right inequality fails, even for linear
functions.

\medskip
\noindent
{\bf Proof of Proposition \KK.1.}   We give the proof of the first equality for all $p$ using $(\BB.7)$. 
Taking  $x_0=0$ and dropping $u$ and $x_0$ from the notation, it says  that 
$$
V(r) \ =\ n \int_0^1 S(rt) t^{n-1}\, dt 
\eqno{(\KK.7)} 
$$
Hence,  we have
$$
{V(r) \over K(r) } \ =\ n\int_0^1 {S(rt) \over K(rt)} {K(rt)\over K(r)} t^{n-1} \,dt.
$$
When $p\neq 2$, $K(rt)/K(r) = 1/ t^{p-2}$, so that letting $r\downarrow 0$ and integrating yields
the first equality in (\KK.2).  When $p=2$, 
$$
{K(rt)\over K(r)} \ =\ 1 + {\log\, t  \over \log\, r},
$$
so letting $r\downarrow 0$ and integrating yields the first equality
$\den^V(u)\ =\ \den^S(u)$ in (\KK.3).

For the proof of the second equalities we show that the mass density 
$\den^{n-p}(\mu)$ ($\mu=\D u$) is the desired multiple of the spherical density $\den^S(u)$. 
Recall the classical fact that 
$$
\mu(B_r)\ =\ (n-2) |S|{S_{-}'(r)   \over K_n'(r)    }.
\eqno{(\KK.8)} 
$$
(See (3.2.13)$'$ in [Ho$_2$,  Thm. 3.2.16] for a proof.)
Since 
$$
{n-2 \over K_n'(r)} \ =\ r^{n-1}\ =\ { |p-2| r^{n-p}  \over K_p'(r) } \quad {\rm when}\ \ p\neq2,
$$
we have
$$
r^{p-n}\mu(B_r)\ =\  |p-2| |S|{S_{-}'(r)   \over K_n'(r)    }  \quad {\rm when}\ \ p\neq2.
\eqno{(\KK.8)'} 
$$
If $p=2$, this holds with $|p-2|$ replaced by 1.
Finally, letting $r\downarrow 0$ and using (\LL.6) completes the proof.\qed

 \medskip
 \noindent
 {\bf Proof of Proposition \KK.2}
 For simplicity let $x_0=0$.  
Note that for all $p$ and $r$  we have $S(u, r) \leq M(u,r)$. 
On the other hand,  $K(r)<0$ when $p \geq 2$ and $K(r)>0$ when $p<2$.
Dividing by $K(r)$ and letting $r\downarrow0$ then gives the inequalities on the left
as well as the inequality $\den^M(u) \leq \den^S(u)$ when $p=2$
(since $u$ and $u+c$ have the same density, we can assume that $u(0)=0$ when 
$p<2$.)

The remainder of the proof is a consequence of Harnak's inequality.
The standard form of this inequality is for a function $v\leq 0$ which is 
$\D$-subharmonic on $B_\rho$.  It says, with $\vf$ defined by
$$
\vf(\l)\ \equiv\   {  1-\l  \over (1+\l)^{n-1}     } \qquad{\rm for}\ \ 0<\l<1,
$$
that
$$
M(v, \l r) \ \leq \  \vf(\l)S(v, r) \fa 0< r \leq \rho.
\eqno{(\KK.9)} 
$$
(See, for example, Prop. 4.2.2 in [D].)
For an arbitrary $\D$-subharmonic function $v$, the function $v-M(v, r)$ is $\leq0$ on $B_r$.
Hence,  (\KK.9) gives the following more general form of Harnak's inequality
$$
M(v, \l r) - M(v,r) \ \leq\ \vf(\l) \bigl( S(v,r) - M(v, r)\bigr) \fa 0< r\leq \rho.
\eqno{(\KK.10)} 
$$
for functions not necessarily $\leq0$.

Suppose first  that $p>2$.  
We may assume $u(0)=-\infty$ since otherwise the assertion  is  trivial.
Then $u$ is negative near 0, and we can apply the standard form   (\KK.9)
of  Harnak's  inequality to obtain
$$
{M(u,\l r)  \over K(\l r) } \ \geq\ \l^{p-2} \vf(\l) {S(u,r)  \over K(r)}.
$$
Letting $r\downarrow0$ gives $\den^M(u,0) \geq c \den^S(u,0)$
where $c =  \l^{p-2} \vf(\l) >0$.  This gives (\KK.4) with    $C=1/c$. (Note that $c\equiv \sup_\l \l^{p-2} \vf(\l) $ provides 
the best constant $C$.)

Suppose now that $1<p<2$. Replace $u$ by $u(x)-u(0)$ so that $u(0)=0$.
Since densities are unchanged by adding a constant, 
we have  $\den^M(u,0) = \lim_{r\downarrow0}M(u,r)/K(r)$ and
  $\den^S(u,0) = \lim_{r\downarrow0}S(u,r)/K(r)$ by Corollary \LL.4. 
  Since $u$ may not be $\leq0$, we must use the general
 form (\KK.10) of Harnak.  Dividing by $K(r)$ gives
$$
{ (1+\l)^{n-1}   \over  1-\l   } \left({M(u,\l r)\over K(r)}  -  {M(u,r)\over K(r)}  \right)
\ \leq\ 
{S(u,r) \over K(r)} - {M(u,r) \over K(r)}.
\eqno{(\KK.11)} 
$$
Using the fact that $K(\l r) = \l^{2-p}K(r)$ and letting $r\downarrow0$ gives
$$
\psi(\l) \den^M(u,0) \ \leq\ \den^S(u,0) \qquad {\rm with}\ \ 
\psi(\l) \ =\  1+{  (1+\l)^{n-1} \over 1-\l} (\l^{2-p}-1).
$$
Now direct calculation shows that $\lim_{\l\downarrow0} \psi'(\l) = \infty$, and so
$c\equiv \sup_{0<\l<1} \psi(\l) >0$. This gives the desired result with $C=1/c$.

It remains to prove that $\den^S(u) \leq \den^M(u)$ when $p=2$.
Set $\l={1/ e}$ in (\KK.11) and note the fact that $K(r) = \log\, r = \log\, {r\over e} +1 = K(\l r) +1
= K(\l r) (1+o(r))$.  Then taking the limit as $r\to 0$ in (\KK.11) yields
$0=\den^M(u)-\den^M(u) \leq \den^S(u)-\den^M(u)$ by Lemma \LL.4.  This completes
the proof of Proposition \KK.2.\qed

\vskip.3in

\centerline{\bf The Upper Semi-Continuity of Density.}\medskip

\Theorem{\KK.4}  {\sl Each  of the densities
$\den^M(u,x)$,  $\den^S(u,x)$,  and  $\den^V(u,x)$
 considered above is an upper semi-continuous function of $x$.}
\pf  
Because of Proposition \KK.1 there are only two cases to consider.  We must show that 
$$
{\limsup_{\eqalign{&x\to x_0\cr &x\neq x_0}}} \ \den(u, x) \ \leq\ \den(u, x_0).
\eqno{(\KK.12)}
$$
Set $x_0=0$.  Assume $0< |x| < r < t$.  Then
$$
\den^\Aa(u,x)  \ \leq\  {\Aa(u,x,t) -\Aa(u,x,r)  \over K(t) -K(r)}.
\eqno{(\KK.13)}
$$

\medskip
\noindent
{\bf Case 1. $\Aa = M$.}  
By using the facts that  $B_t(x) \ss B_{t+|x|}(0)$ 
and $B_{r-|x|}(0) \ss B_r(x)$, we see that
 the last quantity above is
$$
\leq \  {\sup_{B_{t+|x|}(0)} u -  \sup_{B_{r-|x|}(0)} u \over K(t) -K(r)}.
$$
The function $M(u,0,r) \equiv \sup_{B_r(0)} u$   is continuous 
(see Corollary \LL.2(a)) and   increasing.
Therefore,  
$$
{\limsup_{\eqalign{&x\to 0\cr &x\neq 0}}} \ \   \den^M(u,x) \ \leq\ 
{\sup_{B_{t}(0)} u -  \sup_{B_{r}(0) }  u  \over K(t) -K(r)} \qquad 0<r<t.
$$
Finally, the limit of the RHS as $r,t\to0$ equals $\den^M(u,0)$. This proves the first case.

\medskip
\noindent
{\bf Case 2. $\Aa = V$.}  
It  suffices to note that 
$\lim_{x\to0} V(u,x,t) = V(u,0,t)$, which   follows since
$V(u,x,t) = \intave B u(x+ty)\,dy$ and $u$ converges in $L^1(B)$ to $u(ty)$ as $x\to0$.
\qed

\medskip
\noindent
{\bf Note.} 
By using Theorem 3.2.13 in [Ho$_2$], one can show  that $u(x+t\s) \,d\s$
 converges weakly in measure to $u(t\s)\,d\s$ as $x\to0$. This gives a direct proof 
 that $S(u,x,t)$ is continuous in $x$ at $x=0$ without using Proposition \KK.1).
 \qed

\Cor{\KK.5} {\sl
For all $c >0$, the set}
$$
E_c \ \equiv\ \{x : \den(u,x)\ \geq\ c\} \ \ {\sl is\ closed.}
$$

\Note{\KK.6} When $p=1$ the set   where $\den(u)=0$ is just the set
of differentiability points of $u$ (see (\LL.5) in Part II).


\vfill\eject

\centerline{\headfont \BBB.   Maximality of Subharmonics with Harmonic Averages.}
\medskip 

In this section we extend  the standard notion of maximality in  pluripotential theory
  to each $F$-potential theory. This notion extends the notion of being $F$-harmonic,
  but is still very close to it.
In fact, a maximal function is harmonic if and only if it is continuous.
Our main result, Theorem \BBB.2, is key for the study of tangents. 
It provides a new criterion for an $F$-subharmonic function to be $F$-maximal.
An excellent reference for pluripotential theory is [Kl].

\Def{\BBB.1}  An $F$-subharmonic function $u$ on an open set $X\ss\rn$ is said to 
be {\bf $F$-maximal on $X$} if for each $F$-subharmonic function $v$ on $X$ and each
compact subset $K\ss X$,
$$
v\ \leq\ u \quad {\rm on}\ \ X-K   \imp  v\ \leq\ u \quad {\rm on}\ \ X
\eqno{(\BBB.1)}
$$
Note that by replacing $v$ with $\max\{u,v\}$, condition (\BBB.1) is equivalent to
$$
v\ \geq\ u \ \  {\rm on}\  X \quad {\rm and}\quad
v\ =\ u \  \ {\rm on}\  \ X-K   \imp  v\ =\ u \quad {\rm on}\  X
\eqno{(\BBB.1)'}
$$

\smallskip
Most of the previous results  come together in the proof of the next result.

\Theorem{\BBB.2. (The Maximality Criterion)} {\sl
Suppose that $F$ is an \ST convex cone subequation, and $U$ is an  
$F$-subharmonic function on the annulus $A = \{x : a<|x|<b\}$. If the spherical average
$$
 S(U,t) \ \equiv\  \intave{S} U(t \s) \,d\s
\ \ {\sl determines \ an\ increasing \  } F\ {\sl harmonic\ } S(U, |x|)\ {\sl on\ } A(a,b),
\eqno{(\BBB.2)}
$$
then the function
}
$$
U 
\ \ {\sl is \ } F \ {\sl maximal\ on\ } A.
\eqno{(\BBB.3)}
$$

\pf
The hypothesis on $U$ can be restated as the condition 
$$
 S(U,t) \ \ {\sl is \ } R^{\uparrow}_F \ {\sl harmonic\ on\ } (a,b),
\eqno{(\BBB.2)'}
$$
by  Theorem \AA.7.  By Proposition \CC.3, $R^\uparrow_F =  R^\uparrow_p$, so by Proposition \CC.5
this proves that (\BBB.2)$'$ is equivalent to
$$
 S(U,t) \ =\ \den K(t) + c \quad{\rm on}\ \ (a,b)
\eqno{(\BBB.2)''}
$$
for  constants $\den\geq0$ and $c\in\bbr$. Now by the homogeneity  of $S$ and $K$, this is equivalent to
$$
 {S(U,t) -S(U,r) \over K(t)-K(r)} \ =\ \den\ \geq\ 0
\fa r<t\ \ {\rm in}\ \ (a,b)
\eqno{(\BBB.2)'''}
$$
for some constant $\den\geq0$. 

As in (\BBB.1)$'$ assume that $v$ is $F$-subharmonic on $A$ 
with $v \geq U$ and that outside a  compact subset $K\ss A$
we have  $v=U$.  By the fundamental double  monotonicity Theorem  \BB.4
we have that for $a<r<t<b$,
$$
{S(v,t) - S(v,r) \over K(t)-K(r)} \quad{\rm is\ increasing\ in\ \ } r\ \ {\rm and}\ \ t.
\eqno{(\BBB.4)}
$$
Since $v=U$ outside $K$, this quotient equals $\den$ if both $r$ and $t$ are 
sufficiently close to $a$ or sufficiently close to $b$.  Hence, this quotient equals
$\den$ for all $r<t$ in $(a,b)$.  That is,  $S(v,t)$ satisfies (\BBB.2)$'''$.  It follows that 
$S(v,t)$, in addition to   $S(U,t)$, satisfies (\BBB.2)$''$.  Therefore,
$$
S(v,t) \ = \   S(U,t) + c   \quad \forall\, t\in (a,b)
\eqno{(\BBB.5)}
$$
Taking $t$ close to $a$ shows that $c=0$.  Now the fact that $S(v,t)=S(U,t)$
for all $t\in (a,b)$ combined with the inequality $U\leq v$ implies that 
$U=v$ on $A$, thus proving that $U$ is $F$-maximal on $A$.\qed

The following additional facts about $F$-maximal functions are standard in  
pluripotential theory, where $F=\cp^\bbc$.  The proofs easily adapt to the more
general subequation $F$, but since these results are not part of the viscosity 
literature, we inlcude them for the convenience of the reader.
 Throughout the remainder of this section
$F$ is an arbitrary subequation, i.e., a closed set $F\ss \Symn$ which satisfies $F+\cp\ss F$.

\Prop{\BBB.3}
$$
{\sl If}\ \ u\ \ {\sl is } \ F\ {\sl harmonic\ on }\ X, \ {\sl then}\ \ u\ \ 
 {\sl is } \ F\ {\sl maximal\ on }\ X
\eqno{(\BBB.6)}
$$
This is immediate since comparison holds for $F$ (cf. [\HLDD, Thm. 6.5]).
The only thing standing in the way of the converse is the continuity of $u$.

\Ex{\BBB.4}  The subequation $F=\cp^\bbc$ of pluripotential theory has many
functions, such as $\log|z_1|$ on $\bbc^2$,  which are maximal but not $F$-harmonic.
In fact any function $u(z_1)$, which is $\D$-subharmonic on a domain  $X_0\ss \bbc$, 
when considered as a function $\overline u (z) \equiv u(z_1)$ on $X=X_0\times \bbc^{n-1}$ with $n\geq2$, is $\cp^\bbc$-maximal. (If $v(z) \leq u(z_1)$ on $X-K$, then by applying the maximum principle to $v$ on slices $z_1 = $ constant, we get $v(z) \leq \overline u(z)$ on $X$.)
Now  $\overline u (z) \equiv u(z_1)$ is $\cp^\bbc$-harmonic  if and only if $u$ is continuous, 
however,  $u$ is not necessarily continuous even if it is bounded.

\Prop{\BBB.5}
\smallskip
\centerline
 {\sl
If $u$ is $F$-maximal and continuous on $X$, then $u$ is $F$-harmonic on $X$.
}

\pf
 This is the standard ``bump-function'' argument which occurs
for example as far back as [BT] or in [I].  It goes as follows.  Suppose $u$ is not $F$-harmonic
but is $F$-maximal, and therefore $F$-subharmonic.  Then $v\equiv -u$ is not $\ft$-subharmonic.
Therefore, by Lemma 2.4 in [\HLDDR], there exist $x\in X$, $\e>0$ and 
a quadratic polynomial $Q(y)$ such that $v(y) < Q(y) -\e|y-x|^2$ on $\overline{B_r(x)} - \{x\}$
with equality at $y=x$, but $D^2_xQ \notin \ft$, i.e., $-D^2_xQ \in \Int F$.  Thus, 
$w\equiv -Q+\d$ is strictly $F$-subharmonic at $x$, and hence in a neighborhood $B_r(x)$.
Pick $\d>0$ sufficiently small that $v(y)< Q(y)-\d = -w(y)$ on $\partial B_r(x)$.  Then 
$w(y) < u(y)$ on $\partial B_r(x)$, but $w(x) = u(x) +\d$.
This proves that $u$ is not maximal.\qed
\medskip

$F$-harmonic functions may not be closed under decreasing limits.
 For instance in Example \BBB.4 each $u(z_1)$ which is $\D$-subharmonic
 is the decreasing limit of functions $u_j(z_1)$  which are smooth and 
 $\D$-subharmonic.  The extensions ${\overline u}_j \to \overline u$ to $\bbc^n$ give an
 example for the case $F=\cp^\bbc$.
 
 This defect is corrected by enlarging the space of $F$-harmonic functions to the
 space of $F$-maximal functions. (This is the smallest such enlargement by Theorem \BBB.7 below.)

\Prop{\BBB.6}
\smallskip
\centerline
 {\sl
If $u$ is the decreasing limit of a sequence of $F$-maximal functions, then $u$ is $F$-maximal.
}
\pf Suppose  $\{u_j\}$ are $F$-maximal and  $u_j \downarrow u$ on an open set $X$.  Fix a compact
set $K\ss X$.  Then $v\leq u$ on $X-K$ $\Rightarrow$ $v\leq u_j$ on $X-K$ $\Rightarrow$
$v\leq u_j$ on $X$ $\Rightarrow$ $v\leq u$ on $X$.\qed
\medskip

This fact has a strong converse.

\Theorem{\BBB.7}
\smallskip
\centerline
 {\sl
If $u$ is locally $F$-maximal, then $u$ is locally the  decreasing limit}

\centerline {\sl
 $u=\lim_{j\to\infty} u_j$ of $F$-harmonic  functions  $u_j$.
}
\medskip

The proof of this fact requires a lemma.

\Lemma {\BBB.8} {\sl
Suppose $u$ is $F$-subharmonic on $X$, $\O^{\rm open}\ss\ss X$, and $v\in\USC(\ob)$ is 
$F$-subharmonic on $\O$.  If $v\leq u$ on $\bo$, then
$$
\overline v \ \equiv\ 
\cases
{
\max \{u,v\} \ \ {\sl on}\ \ \ob   \cr
\ \ \ \ \ u\qquad\ \   {\sl on}\ \ X-\ob 
}
$$
is $F$-subharmonic on $X$.
}
\pf  Sup-convolution provides a decreasing sequence $u^\e \downarrow u$ 
of continuous  $F$-subharmonic functions which are defined on  subdomains
which contain $\ob$ and increase to $X$.  Set
$$
 v^\e_\d \ \equiv\ 
\cases
{
\max \{u^\e +\d,v\} \ \ {\sl on}\ \ \ob   \cr
\ \ \  \ \ u^\e+\d\qquad\ \   {\sl on}\ \ X-\ob.
}
$$
Since $\{v< u^\e+\d\}$ is a relatively open subset of $\ob$ containing $\bo$, the function
$v^\e_\d$ is $F$-subharmonic on domains containing $\ob$ which increase to $X$
as $\e\downarrow 0$. Finally, $v^\e_\d\downarrow \overline v$ as $\e,\d \downarrow 0$,
proving that $\overline v$ is $F$-subharmonic on $X$.\qed

\medskip

Using this Lemma \BBB.8 the definitions (\BBB.1) and (\BBB.1)$'$ of $F$-maximality on $X$
can be further refined as follows:

For each domain $\O\ss\ss X$ and $v\in\USC(\ob)$ which is $F$-subharmonic on $\O$,
$$
v\ \leq\ u \quad {\rm on}\ \ \bo \imp v\ \leq\ u \quad {\rm on}\ \ \ob
\eqno{(\BBB.1)''}
$$

Using this definition of $F$-maximality together with the fact that on balls $B\ss\rn$ the 
Dirichlet problem is uniquely solvable by the Perron function, it is easy to prove Theorem \BBB.7.
\medskip
\noindent
{\bf Proof of Theorem \BBB.7.}
Suppose $u$ is maximal on $X$ and $\overline B\ss X$ is a closed ball.  Choose
$\vf_j \in C(\partial B)$ such that $\vf_j\downarrow u\bigr|_{\bo}$.
Let $u_j \in C(\ob)$
denote the solution to the Dirichlet Problem on $\overline B$
 with $u_j\bigr|_{\bo} = \vf_j$ and $u_j$ $F$-harmonic on $B$.  
Since $u_j$ is the Perron function for boundary values $\vf_j$, we have $u\leq u_j$
for all $j$ and $u_j\downarrow v$ is decreasing.  Thus $u\leq v$.   
Also $v\bigr|_{\partial B} = \lim u_j \bigr|_{\partial B}  = \lim \vf_j = u \bigr|_{\partial B}$, 
and  $v$ is $F$-subharmonic on $B$.
Thus,  by (\BBB.1)$''$ above, $v\leq u$ on $\overline B$.  Hence, $u=v=\lim u_j$.\qed


\vskip .3in

\centerline{\headfont \FF.   Tangents to Subharmonics.}
\medskip 

Now we come to the main topic of the paper --   introducing the notion of tangents to $F$-subharmonics. 
In this section the  \ST cone subequation   $F$ on $\rn$ is assumed to be convex.
We shall work at a fixed point, which for simplicity is assumed to be the origin.
 That is,
given an $F$-subharmonic function $u$ defined in a neighborhood of 0, we 
define the notion of tangent functions  to $u$ at 0.   A required clarification is
given  by Proposition \FF.4.   The basic properties of a tangent $U$ to $u$ at 0 are then  established
in Theorems \PP.4 and \FFF.2.

\Def{\FF.1} Suppose that $u$ is $F$-subharmonic on the ball about the origin of radius $\rho$.
The {\bf tangential $p$-flow} (or {\bf $p$-homothety}) determined by the Riesz characteristic $p=p_F$ of $F$ 
is defined as follows.
\medskip

(a) \ \ $u_r(x)\ =\ r^{p-2}u(rx) \qquad$ if \ \ $p\ >\ 2$, 
\medskip

(b) \ \ $u_r(x)\ =\  {1\over r^{2-p}} \left[  u(rx) - u(0) \right]   \qquad$ if \ \ $1\ \leq \ p\ < \ 2$, and
\medskip

(c) \ \ $u_r(x)\ =\ u(rx) -  M(u,r)  \qquad$ if \ \ $p\  =\ 2$

\Remark{\FF.2}  
Suppose $1\leq p<2$.  Since $u(0)=M(u,0)$ is finite,
some readers may prefer to assume once and for all in part (b) that $u(0) =0$ so 
that the $p$-flows for all $p\neq2$ are the same, namely that
$$
u_r(x) \ =\ r^{p-2}u(rx)\qquad {\rm if} \ \ p\neq 2.
\eqno{(\FF.1)}
$$
Others may wish to make this assumption in the proofs.
\medskip

Note that
$$
\eqalign
{
{\rm The\ functions}\ u_r\  {\rm are} \  F\, {\rm subharmonic\  on \ } B_{\rho/r},  \cr
 \ {\rm and\ as\  } r\ \to\ 0, \ {\rm these\  balls\ expand\ to\ }\rn.\quad
 }
$$

An upper semi-continuous  function $U(x)$ on $\rn$ taking values in $[-\infty, \infty)$  is invariant under this flow if and only if there exists an u.s.c. function $g$ on the
unit sphere $S$ such that 
$$
U(x) \ =\ |x|^{p-2} g\left(  {x\over |x|}\right) \qquad{\rm in\ the\ cases\  where \ \ } p\neq2,
$$
while in the case where $p=2$, we leave it to the reader to prove that
$$
U(x)\ =\ \den \log |x| + g\left(  {x\over |x|}\right) \qquad {\rm with}\ \ \sup_{S^{n-1}} g \ =\ 0\ \ \ {\rm and}\ \ \ 
\den\geq 0\ \ {\rm a \ constant.}
$$
Functions of this form will be said to have {\bf Riesz homogeneity} $p$.

Under our assumptions on $F$ each $F$-subharmonic function $u$ is $\lloc$ since it is 
$\D$-subharmonic by (\BB.3).

\Def{\FF.3. (Tangents)}  
Suppose that $u$ is an $F$-subharmonic function defined in a neighborhood of the 
origin.  For each sequence $r_j \searrow 0$ such that
$$
\overline U \ \equiv\ \lim_{j\to\infty} u_{r_j} \ \ {\rm converges\  in\ } \lloc(\rn),
\eqno{(\FF.2)}
$$
the point-wise defined function 
$$
 U(x) \ \equiv\ \lim_{r\to 0} \, {\rm ess} \! \! \!\! \sup_{B_r(x) \quad } \overline U
\eqno{(\FF.3)}
$$
is called a {\bf tangent to $U$ at 0}.  We let $T_0(u)$ denote the set of all such tangents $U$.
(We will refer to $\overline U$, satisfying (\FF.2), as an  {\bf $\lloc$-tangent} 
when the distinction between the function $U$ and the equivalence class of functions
$\overline U$ is important.)

Our first result clarifies this  Definition.

\Prop{\FF.4} {\sl
Each tangent $U$ to $u$ at 0 is an entire $F$-subharmonic function on $\rn$.
Moreover, $U$ belongs to the $\lloc$-class $\overline U\in \lloc(\rn)$
and is the unique $F$-subharmonic function in this  $\lloc$-class.
}
\medskip

To prove Propostion \FF.4 we use the following result established in [\HLPUP, Cor. 5.4]
(see [\HLBP] for generalizations.) We say that a subequation $F$ {\bf can be defined using fewer of the variables in $\rn$}
if there exist an $(n-1)$-dimensional subspace $W\ss\rn$ and a subequation $F'\ss\Sym(W)$
which determines $F$ by: $A\in F \iff A\bigr|_W \in F'$.

An important point here is that the same representative $u$ of the $\lloc$-class $\overline u$
 (given by (\FF.4)) is the correct representative, no matter  which subequation $F$ is being considered.

\Theorem{\FF.5. (Distributional versus Viscosity Subharmonics)} {\sl
Suppose $F$ is a convex cone subequation which cannot be defined using fewer of
the variables in $\rn$. 
\medskip
\item{(a)} 
If $u$ is $F$-subharmonic in the viscosity sense, then $u$ is $\lloc$ and $F$-subharmonic in the distributional sense.

\medskip
\item{(b)}
If $\overline u$ is an  $F$-subharmonic distribution,  then $\overline u \in \lloc$ and 
the limit
$$
u(x) \ =\ \lim_{r\to0} \, {\rm ess} \! \! \!\! \sup_{B_r(x) \quad} \overline u
    \qquad{\rm exists\ at\  each \ point}
\eqno{(\FF.4)}
$$

and defines an upper semi-continuous function $u$ in the $\lloc$-class $\overline u$
 which is  $F$-subhar-
 
monic  in the viscosity sense.
Moreover, $u$ is the unique such representative of $\overline u$.
}

\medskip
\noindent
{\bf Remark \FF.6.} We refer the reader to Sections 3,4, and 5 of [\HLPUP] for a fuller discussion of this result
and the definition of an $F$-subharmonic distribution (Definition 4.1 and Proposition 4.3).
However, the terminology used in  [\HLPUP]  is somewhat different.  Here we use the terminology
emplyed in [\HLBP].   In  [\HLPUP]  a convex cone subequation $F$ is called a ``positive cone''
and denoted $\cp^+$.  The polar cone is denoted by $\cp_+$.
A convex cone subequation which cannot be defined using fewer of the variables in $\rn$ is called
an {\sl elliptic cone''.}
\medskip

From the distributional point of view it is straightforward to see that averages, or more generally
convolution, of an $F$-subharmonic  function $u$ with any non-negative measure is again $F$-subharmonic.

\medskip
\noindent
{\bf Proof of Proposition \FF.4.}  
 We use these facts about the \ST convex cone subequation $F$:

\centerline
{(1) \ \ $F\ \ss\ \Delta$ \qquad (2)\ \ $1\ \leq\ p_F\ \leq\ n$
}
\smallskip
\centerline{
\qquad (3) \ \ $F$ cannot be defined using
fewer of the variables in $\rn$.
}

\medskip
\noindent
{\bf Proof}.  Properties (1) and (2) have already been noted in (\BB.3) and (\BB.4).
For Property (3) note that the ST-invariance of $F$ rules out the possibility that $F$
could be defined using fewer of the variables in $\rn$.  Because of (3) one can apply Theorem \FF.5.   

  Suppose 
$\overline U = \lim_{j\to\infty}  u_{r_j}$ in $\lloc(\rn)$ is an $\lloc$ tangent to $u$ at 0.  
Since $F$ is a cone, each $u_r$ is viscosity $F$-subharmonic, and hence in $\lloc$ and 
distributionally $F$-subharmonic by Part (a) of Theorem \FF.5.
Hence, in the limit, $\overline U$ is distributionally $F$-subharmonic.
Now apply Part (b) of Theorem \FF.5 to $\overline U$ to complete the proof.\qed
\medskip

In light of Proposition \FF.4 we frequently drop the distinction between $U$ and 
$\overline U$.

\vskip .3in


\centerline{\headfont  \PP.  Uniqueness of  Averages of Tangents and of Flows.} 
\medskip

Most of the   properties of tangents  can be deduced from the following result, which 
proves that averages of tangents are always  unique by showing that they are radial harmonics.

\Theorem{\PP.1.  (Averages of Tangents)} {\sl
Suppose that $u$ is an $F$-subharmonic function defined in a neighborhood of
the origin in $\rn$.  Let $p=p_F$ be  the Riesz characteristic of $F$. 

 If $p \neq2$,
then each tangent $U$ to $u$ at 0 has averages
$$
\eqalign
{
M(r) \ \equiv \ \sup_S U(r\s) \ = \ &\den^M(u) K(r), 
\qquad    S(r)\ \equiv \ \intave S U(r\s) \, d\s  \ =\ \den^S(u) K(r),\cr
 \ \ \  &{\rm and}\ \ \ V(r)\ \equiv\ \intave B U(rx)\, dx  \ =\ \den^V(u) K(r)
 }
\eqno{(\PP.1)}
$$
In particular,
$$
\den^\A(U)\ =\ \den^\A(u) \qquad {\sl for} \ \ \A\ =\ M, S,\  {\sl or} \ V
\eqno{(\PP.2)}
$$

When $p=2$,  all the densities of $u$ and  any tangent $U$ to 
$u$ at 0,  agree, and will be simply denoted by $\den= \den(u)$.
Specifically, we have
$$
\den(u) = \den^M(U) = \den^S(U) =  \den^V(U) =  \den^M(u) =  \den^S(u) = \den^V(u).
\eqno{(\PP.3)}
$$
Moreover, the averages of a tangent $U$ to $u$ are given by
$$
M(r) \ =\ \den \,\log\, r, 
 \qquad  S(r)\ =  \den\, \log\, r +\intave S U,  
 \ \ \   {\rm and}\ \ \ V(r)\ =\  \den\, \log\, r +\intave B U,
\eqno{(\PP.4)}
$$
with
$$
-C\den \ \leq \  \intave S U\ \leq \ 0
\quad{\rm and}\quad
-(C+1) \den \ \leq\ \intave B U, \quad 
{\rm where} \ \ C \ =\ {1\over \vf \left({1\over e}\right) }\ >\ 1.
\eqno{(\PP.5)}
$$
and where $\vf(\l) = (1-\l)/(1+\l)^{n-1}$.}
\medskip

When $p\neq 2$, these formulas show  that any two tangents have the same  
maxima  $M(r)$ and the same spherical averages $S(r)$  and volume averages $V(r)$,
all being the appropriate density times $K(r)$. 
When $p=2$, $M(r)$, $S(r)$ and $V(r)$ all agree with $\den \log\,r$ modulo an additive 
constant, but the constant depends on the tangent $U$, not just on $u$.

  In all cases, for each tangent $U$, the function $S(U,|x|)$ is $F$-harmonic 
  on $\rn-\{0\}$ since $\den K(|x|) +C$ is $F$-harmonic there (Proposition \CC.5).

Combining Theorem \BBB.2 and Theorem \PP.1
is one of the main ingredients of the paper and has the following immediate consequence. 

\Theorem{\PP.2} {\sl   
Every tangent to an $F$-subharmonic function is $F$-maximal.}

\medskip
Applying Proposition  \BBB.5 yields the following.

\Cor{\PP.3}
{\sl Every continuous tangent to an $F$-subharmonic function is $F$-harmonic.}
\medskip

 Theorem \PP.1, the uniqueness of averages of tangents,  follows from the stability of averages (Lemma \BB.5)
and the uniqueness  of the averages of a flow. Its proof is given at the end of this section.

We may assume that $u(0)=0$ if $1\leq p<2$ (Remark \FF.2), and that $u(0)= -\infty$ 
if $2\leq p <\infty$.

\Theorem{\PP.4. (Averages of Flows)} {\sl
For $p\neq 2$ and $\A = M, S$ or $V$,
$$
\lim_{s\downarrow0} \A(u_s, r) \ =\ \den^\A(u) K(r).
\eqno{(\PP.6)}
$$

For $p=2$, if $\A=M$ we also have
$$
\lim_{s\downarrow0} M(u_s, r) \ =\ \den^M(u)  K(r)\ =\ \den^M(u) \log\, r.
\eqno{(\PP.6a)}
$$
In this case the limit is decreasing and uniform in $r\leq R$.
For $\A = S$ or $V$ we have 
$$
\liminf_{s\downarrow0} S(u_s, r) \ \geq\ \den^M(u) (\log\, r - C), \ \ {\rm and}
\eqno{(\PP.6b)}
$$
$$
\liminf_{s\downarrow0} V(u_s, r) \ \geq\ \den^M(u) (\log\, r - C-1)
\eqno{(\PP.6c)}
$$
with $C$ as in (\PP.5).}

\medskip

Direct calculations from the definitions of the flow and the averages establish the next
result.

\Lemma {\PP.5}  {\sl For $\A = M, S$,  or $V$:}
$$
\A(u_s, r) \ =\ s^{p-2} \A(u,sr)  \ =\   {\A(u, sr) \over K(sr)} K(r)  \qquad{\rm when }\ p\neq 2, \ \ \ {\rm and}
\eqno{(\PP.7)}
$$
$$
\A(u_s, r) \ =\  \A(u,sr) - M(u,s)  \ =\   {\A(u, sr) - M(u,s)\over K(sr)-K(s)} K(r) \qquad{\rm when }\ p=2
\eqno{(\PP.8)}
$$

 \pf
 For example, when $\A$ is the volume average $V$ and $p\neq 2$, we have
 $$
 V(u_s, r)\ =\ {1\over |B|} \int_B u_s(rx)\, dx 
 \ =\     {s^{p-2}\over |B|} \int_B u(rsx)\, dx 
 \ =\  s^{p-2} V(u, rs).\qquad
$$
The remaining calculations are left to the reader.\qed

\medskip

\noindent
{\bf Proof of Theorem \PP.4.} 
By Lemma \LL.4 the identity (\PP.7)
 implies (\PP.6) for $p\neq 2$.  In  the case
where $p=2$ the limit  (\PP.6a) for the maximum
follows from (\PP.8) by the double monotonicity Theorem \BB.5.
 The limit  (\PP.6c) for $V$ follows from the limit
 (\PP.6b) for $S$ since 
 $
 V(u_s, r)  = n\int_0^1 S(u_s, t) \, t^{n-1} \, dt
 $
 by (\BB.7), and 
 $
 n\int_0^1(\log\, rt - C) \, t^{n-1} \, dt = \log \, r - C - 1.
 $
 
 It remains to prove (\PP.6b).  
Harnak's inequality  in the form (\KK.10) with $v=u_s$ and $\l=1/e$ states that
$$
C \left( M\left( u_s,  \smfrac r e \right) - M(u_s, r)\right) + M(u_s, r) \ \leq \ S(u_s, r).
$$
We know the limit of the terms involving $M$  as $s\downarrow 0$. This gives
$$
C \den^M(u) \left( \log\, \smfrac re - \log\, r\right) + \den^M(u) \log\,r \ \leq\ 
\liminf_{s\downarrow0} S(u_s, r)
$$
as desired.
\qed
\medskip

\noindent
{\bf Proof of Theorem \PP.1.}  
The density statements for $u$ are contained in Propositions \KK.1 and \KK.2.
The density statements for $U$ follow from the formulas in Theorem \PP.1 and the
density statements for $u$.  The  formulas in Theorem \PP.1 follow immediately 
from the formulas in Theorem \PP.4 for the averages of  flows and the stability
of averages (Lemma \BB.5), with the exception of (\PP.4)  for 
$S$ and $V$, and the estimates in (\PP.5).

The estimates (\PP.6b) and (\PP.6c) and the Stability Lemma \BB.5 show that 
for any tangent $U$ to $u$ at 0,
$$
\den^M(u) \left(  \log \, r -C   \right)  \ \leq\ S(U,r)
\and 
\den^M(u) \left(  \log \, r -C -1  \right)  \ \leq\ V(U,r)
$$
for all $0<r<\infty$.  Also we have that 
$
V(U,r) \leq S(U,r) \leq M(U,r) = \den^M(u)\log\,r$.

Since $V(U,e^t)$ and $S(U,e^t)$ are entire convex functions of $t$, the linear inequalities
$$
\den(t-C) \ \leq\ S(U,e^t) \ \leq\ \den t 
\and
\den(t-C-1) \ \leq\ V(U,e^t) \ \leq\ \den t 
$$
imply that $S(U, e^t) = \den (t+k)$ and $V(U,e^t) = \den(t+k'-1)$ where $k$ and $k'$ satisfy
$-C \leq k, k'\leq 0$.\qed


\vfill\eject

\centerline{\headfont   \FFF. Existence of Tangents.}
\medskip

We now address the basic existence question.  Again $F$ is assumed here to be convex.
However, in the case where $1\leq p<2$ much stronger results are true even if 
$F$ is just a cone and not necessarily convex. These stronger results are established in Section \NN.

\Theorem{\FFF.1. (Existence of Tangents)} {\sl
Suppose that $u$ is $F$-subharmonic on a ball $B_\rho$.
For each $R > 0$ there exists $\d>0$ such that the family $\{u_r\}_{0<r\leq \d}$
 is unformly  bounded above and bounded  in  norm in $L^1(B_R)$.
In particular, the set $\{u_r\}_{0<r\leq \d}$ is precompact in $L^1(B_R)$.
}

\pf
An upper bound for $u$ can be chosen to be any number  greater than $\den^M(u) K(R)$
by (\PP.6) if $p\neq 2$ and by (\PP.6a) if $p=2$.  Consequently the boundedness in
$L^1(B_R)$ is equivalent to a lower bound for $V(u_s, R)$ which is uniform in $s$.
This lower bound can be chosen to be any number less  than $\den^V(u) K(R)$ if $p\neq 2$,
or $\den^M(u) (\log\,R -C-1)$ if $p=2$, by (\PP.6) and (\PP.6c) respectively in Theorem \PP.4.
\qed

\bigskip
The basic properties of the tangent set $T_0(u)$ are contained in the following theorem.
Again see Section \NN\ for the stronger versions of parts (2) and (4) 
using the H\"older topology instead of the $\lloc$-topology.
when $1\leq p<2$.

\Theorem{\FFF.2} {\sl
Suppose that $u$ is an $F$-subharmonic function defined in a neighborhood of
the origin in $\rn$.
Then the tangent set $T_0(u)$ to $u$ at $0$ satisfies:

\medskip

\item{(1)}  \ \  $T_0(u)$ is non-empty.

\medskip

\item{(2)}  \ \  $T_0(u)$ is a compact subset of $\lloc(\rn)$.

\medskip

\item{(3)}  \ \  $T_0(u)$ is invariant under the homothety $U \to U_r$.

\medskip

\item{(4)}  \ \  $T_0(u)$ is a connected subset of $\lloc(\rn)$.

\medskip
 }

\pf
 Parts (1) and (2)  are immediate from Theorem \FFF.1.
The arguments for parts (3) and (4) 
are given in [S, Proposition 1.1.1]. 
We include them  here for completeness.  To prove (3) note that
 $U(x) = \lim_{r_j\downarrow 0} u_{r_j }(x)$  implies
 $U_r(x) = \lim_{s_j\downarrow 0} u_{s_j }(x)$ with $s_j= r r_j$.
 To prove (4) suppose $u_{r_j} \to U_0$ and $u_{t_j} \to U_1$ with
 $U_0$ and $U_1$ elements of disjoint open sets $N_0$ and $N_1$
 which cover  $T_0(u)$.  We can assume $r_j<t_j$ for all $j$ and choose
 $s_j$ between $r_j$ and $t_j$ with $u_{s_j} \notin N_0\cup N_1$.
 (Note that $s\mapsto u_s$ is a continuous map into $\lloc$.)
 By Theorem \FFF.1  the sequence  $u_{s_j}$ has a convergent 
 subsequence, and its limit is in neither $N_0$ nor $N_1$, a contradiction.\qed



\vfill\eject

\centerline{\headfont   \JJ. Uniqueness of Tangents.}
\medskip 

In this section we discuss some basic situations where 
tangents are unique. Our main uniqueness results are are stated
and proved in subsequent sections.  As in Sections \FF-\FFF we assume that $F$ 
is  convex  with  finite  Riesz characteristic  $p$.

\Def{\JJ.1}  Suppose $u$ is an $F$-subharmonic function defined in a neighborhood of the origin.
\medskip

\item{(a)} \ If $T_0(u) =\{U\}$ is a singleton, then we say that {\bf uniqueness of tangents holds for $u$.}
If  uniqueness of tangents holds for all such $u$, we say the that {\bf uniqueness of tangents holds for $F$.}
\medskip

\item{(b)} \ If $T_0(u) =\{\den  K(|x|)\}$ with $\den\geq0$ a constant, 
then we say that {\bf strong uniqueness of tangents holds for $u$.}
If  strong uniqueness of tangents holds for  all such $u$, then
 we say that {\bf strong uniqueness of tangents holds for $F$.}

\medskip

\item{(c)} \ If every tangent $U$ to $u$ satisfies $U_r=U\   \forall\,r$, then
 we say that {\bf homogeneity of tangents holds for $u$.}
 If  homogeneity of tangents holds for all such $u$, then we say that 
{\bf homogeneity of tangents holds for $F$.}

\medskip

Now (b) $\Rightarrow$ (a) $\Rightarrow$ (c).  The first implication is obvious. For the second, note 
 that (a) can be rephrased since
$$
T_0(u) \ =\ \{U\}\qquad\iff\qquad 
\lim_{r\to0} u_r\ \ {\rm exists\ in\ } \lloc(\rn)\ {\rm and\  equals\ } U.
\eqno{(\JJ.1)}
$$
Thus by (a), $u_{r_j}$ and  $u_{rr_j}$ have the same limit $U$, but $u_{rr_j}$
has limit $U_r$, which proves (c).

In general, $S(u,r)\leq M(u,r)$. Therefore, 
$$
{\rm For} \ \ 2\leq p\leq n, \ \den^M(u) \leq \den^S(u),\ \ \quad{\rm and\ for\ }\ 1\leq p <2, \  \den^S(u) \leq \den^M(u)
\eqno{(\JJ.2)}
$$
by (\LL.7) since $K>0$ in the first case and $K<0$ in the second case.   However, if strong uniqueness holds for $u$,
then all densities ``agree''  because of Proposition \KK.1 
and the following.
$$
{\rm If\ for\ some\ } \ \den \geq0, \ T_0(u) = \{\den K\}, \ 
\ {\rm then\ \ } \den^M(u)=\den^S(u) = \den.
\eqno{(\JJ.3)}
$$
This follows  from (\PP.2)  and the fact that $\den^M(K) = \den^S(K) =1$.

There are two classical cases where strong uniqueness holds,
that will prove useful later. For the sake of completeness we include 
proofs. 

\Prop{\JJ.2. (Radial Subharmonics)} {\sl
Suppose that $u(x) = f(|x|)$ is a radial $F$-\sh\ function defined on a neighborhood of 0.
Then
$$
\lim_{r\to0} u_r \ =\  \den(u) K_p(|x|)
$$
in $\lloc(\rn)$ and uniformly on compact subsets in  $\rn-\{0\}$.
Thus, $T_0(u) =\{\den K_p\}$.
}

\pf
Since $u$ is radial, we have that $u_r(x) = M(u_r, |x|)$, but by Theorem  \PP.4 we know that
$\lim_{r\downarrow0} M(u_r, |x|) = \den K_p(|x|)$
uniformly in $0< |x|\leq R$.
\qed

 \Remark {\JJ.3}  The conclusion of convergence in $C(\rn-\{0\})$ only requires $F$ to be an \ST
cone subequation with finite characteristic. It does not require convexity.

\Prop{\JJ.4. (Newtonian Case)} {\sl
Suppose $u$ is a $\D$-subharmonic function defined on a neighborhood of  0.  Then}
$$
\eqalign
{
\lim_{r\to 0} u_r (x) \ &=\ -{\den(u) \over |x|^{n-2}}\qquad \ \ \
{\rm in}\ \ \lloc(\rn)\qquad {\rm when}\ \ n\geq3, \ \ {\rm and}  \cr
\lim_{r\to 0} u_r (x) \ &=\  \den(u) \log\, |x|\qquad
{\rm in}\ \ \lloc(\rn)\qquad {\rm when}\ \ n = 2
}
$$

\pf
Each such $u$ is of the form $u=v+h$ where $v=K * v$ is a Newtonian potential
and $h$ is harmonic near the origin. (Take the measure $\nu$ to be a cut-off of the measure
$\mu=\D u$ in a small ball about the origin.) This reduces the proof to the case $v\equiv K * \nu$.
(In the $n=2$ case  $u_r$ and $v_r+h_r$ differ by $M(v,r) + M(h,r) -M(u,r)$, but this error has limit
zero.)

Now one checks that: for $n\geq3$, $(K * \nu)_r = K * (({1\over r})_*\nu)$ and   for $n=2$,
$(K * \nu)(rx) = K * (({1\over r})_*\nu)(x) + \nu(1) \log \, r$, so that 
$M(K * \nu, r) = M(K *  ({1\over r})_*\nu, 1) + \nu(1)\log\,r$.
Now $\lim_{r\to0} ({1\over r})_*\nu$ always exists weakly in the space of measures and equals
$\den [0]$, where $\den = \lim_{r\to0} \nu(B_r)$ is the zero-dimensional density of $\nu$ at 0.
Since $K\in\lloc(\rn)$, the limit of $(K * \nu)_r $ exists in $\lloc(\rn)$ and equals $K * (\den [0]) = \den K$.
(Note that for $n=2$, $M(K * ({1\over r})_*\nu, 1)$ has limit $M(\den \log\,|x|, 1)=0$.)\qed
\medskip

In the $n=2$ case there is a different proof following  Kiselman [K$_1$].
Note that by (\PP.4) we have $M(U,r) =\den \log\,r$ for any tangent $U$ to $u$ at 0.
In particular, $U(x) -\den \log\,|x|$ is $\leq0$ on $\bbr^2$ and $\D$-subharmonic on $\bbr^2-\{0\}$.
Hence, it can be extended to $\bbr^2$ as a subharmonic function, and then by Liouville's Theorem
it must be constant.  Since $M(u_r,1) = 0$ for all $r$ small, $M(U,r)=0$, proving that the constant 
is zero.\qed

\medskip

Proposition \JJ.4 can be partly generalized.

\Prop{\JJ.4$'$. (Riesz Potentials, $p> 2$))} {\sl
Suppose $u = K_p * \nu$ where $\nu\geq 0$ is a compactly supported measure.
Then
$$
\lim_{r\to 0} u_r \ =\ -{\den(\nu) \over |x|^{p-2}}\qquad{\rm in}\ \ \lloc(\rn)
$$
where, up to a universal constant, $\den(\nu) \ =\ \lim_{r\to 0} \nu(B_r)$}
\pf
Ignoring constants, we have (cf. [L])
$$
\D u \ =\ (\D K_p) * \nu \ =\ K_{p+2} * \nu\ \equiv \ \mu.
$$
Note that
$$
K_n * \mu \ =\ K_n * K_{p+2} * \nu \ =\ K_p * \nu \ =\ u.
$$
We  compute  that
$$
u_r(x) \ =\ r^{p-2} u(rx)\ =\ r^{p-2} (K_p * \nu) (rx)  
\ \ {\rm is\ equal\ to\ \ } 
K_p * \left\{ \left({\smfrac 1 r}\right)_* \nu\right\},
$$
and observe that
$
\lim_{r\downarrow 0} \left({\smfrac 1 r}\right)_* \nu \ =\ \Theta(\nu) [0].
$
\qed\medskip

We complete this section with a final case where strong uniqueness  holds.

\Prop{\JJ.5. (Zero Density)}  {\sl
Suppose that $u$ is $F$-subharmonic in a neighborhood of the origin and $F$ is convex
with $p>1$.
If any of the densities of $u$ is zero at 0, then all the densities of $u$ vanish at 0,
and in this case
$$
\lim_{r\to0} u_r \ =\ 0\quad {\rm in} \ \lloc(\rn).
\eqno{(\JJ.4)}
$$
If $F$ is not convex but $1\leq p<2$, then  $\den^M(u,0)=0$ implies that}
$$
\lim_{r\to0} u_r \ =\ 0\quad {\sl locally\  in} \ \a\ {\sl Holder\  norm,} \ \a=2-p.
\eqno{(\JJ.5)}
$$

\pf
The equality of zero densities is a direct consequence of Propositions \KK.1 and \KK.2,
while (\JJ.4) follows from Theorem \PP.4.

The proof of the final assertion of Proposition \JJ.5 is postponed
as it follows immediately  from (\NN.9).\qed


\vskip.3in

\centerline{\headfont \JJJ.   The Strong Uniqueness Theorem I.}
\medskip 

In this section we give two proofs of one of our two main results  concerning strong uniqueness.
Recall that every O$(n)$-invariant subequation $F$ has  
complex and quaternionic analogues $F^\bbc$ and $F^\bbh$,  which are invariant
under U$(n)$ and Sp$(n)$ respectively  (see Example \HH.7).

\Theorem {\JJJ.1} {\sl
Suppose that $F$ is  O$(n)$-invariant and convex  with finite Riesz characteristic $p$.  
Then, except for the case $F=\cp$, 
strong uniqueness of tangents holds for $F$. Furthermore, 
except for the cases $\cp^\bbc$ and $\cp^\bbh$,
strong uniqueness  of 
tangents also holds  for the complex and quaternionic analogues $F^\bbc$ 
and $F^\bbh$ of $F$.
}

\Remark {\JJJ.2}  For the subequations $\cp, \cp^\bbc$ and $\cp^\bbh$,  strong uniqueness fails dramatically.  Nonetheless, tangents are classified in these cases.  This is discussed
in  Part II of this paper.
\pf
 Let $u$ be $F$-subharmonic
in a neighborhood of the origin and choose $U\in T_0(u)$.
Then
$$
U(x) \ =\ \lim_{j\to\infty} u_{r_j} (x)
$$
for a sequence $r_j\downarrow 0$, where  the flow $u_{r_j} (x)$, given in Definition \FF.1, depends on $p$.

Theorem \PP.2 states that 
$$
U \in T_0(u) 
\qquad\Rightarrow\qquad
U\ \ {\rm is \ F\ maximal \ on\ \ } \rn-\{0\}, \ {\rm and}
\eqno{(\JJJ.1)}
$$
$$
U \in T_0(u) \ \ {\rm and}\ \ 
  U\in C(\rn-\{0\})
  \qquad\Rightarrow\qquad
  U\ \ {\rm is\ \ F\ harmonic\ on}\ \ \rn-\{0\}.
  \eqno{(\JJJ.2)}
$$
We first  prove the theorem under the additional assumption that
$F$ is uniformly elliptic. 
(Note, however, from Section \HH \  that there many examples of subequations $F$ which are not uniformly elliptic, but for which the theorem still applies.)

\Prop{\JJJ.3}   { \sl
If, in addition to the hypotheses of Theorem \JJJ.1, $F$ is uniformly elliptic, then
strong uniqueness of tangents holds for $F$.
}

\pf Two regularity results are needed for $F$.  They   can be stated as follows.

\medskip\noindent
{\bf Fact \JJJ.4.}  A sequence $\{u_j\}$ of $F$-harmonics on $X^{\rm open}\ss\rn$,
which is bounded in $L^\infty(K)$ for each compact $K\ss X$, is precompact in $C(X)$.

\medskip\noindent
{\bf Fact \JJJ.5.}  Each  $F$-harmonic function is $C^1$. 
\medskip
The reader is referred to [CC] and [T] for these results.
 Also for Fact \JJJ.4 one can use the Krylov-Safanov H\"older Estimate 4
in [E] which holds with $f=0$ because of the First Linearization on  page 107.

Recall that $F$ is assumed to be invariant under a subgroup $G\subseteq {\rm O}(n)$
which acts transitively on $S^n$.

\Lemma{\JJJ.6} {\sl
\smallskip
(a)  Suppose $U\in T_0(u)$. Then $g^*U\in T_0(g^*u)$ for each $g\in G$, and the densities

\qquad
$\den^S(g^*U) = \den^S(U) = \den^S(u) = \den^S(g^*u)$ are all equal.
\smallskip
(b)  If  $U\in T_0(u)$ and $V\in T_0(v)$, then $\max\{U,V\} \in  T_0(\max\{u,v\})$.
\smallskip
(c)  If  $U\in T_0(u)$ and $g\in G$, then $\max\{U,g^* U\}\in T_0(\max\{u,g^* u\})$.
}
\medskip
The straightforward proofs are omitted.
\medskip

The proof of Proposition \JJJ.3 will progress in three stages.  First we establish strong uniqueness
for continuous tangents, then for tangents which are locally bounded, and finally for 
general tangents.

The proof that $U=\den K_p$ for $U\in C(\rn-\{0\})$ is as follows.  
Note that for $g\in G$,   $\max\{U, g^*U\} \in C(\rn-\{0\})$, and 
therefore by  Lemma \JJJ.6 and (\JJJ.2), 
$$
\max\{U, g^*U\} \ {\rm is }\ F\,{\rm harmonic\ on \ }\rn-\{0\} \ {\rm for\ each\ } g\in G.
\eqno{(\JJJ.3)}
$$
By the $C^1$-regularity result Fact  \JJJ.5 we have that
$$
\max\{U, g^*U\} \ {\rm is }\ C^1 \,{\rm  on \ }\rn-\{0\} \ {\rm for\ each\ } g\in G.
\eqno{(\JJJ.4)}
$$

Although the maximum of two $F$-subharmonics is always subharmonic,
it is unusual for the maximum of two distinct $F$-harmonics to be $F$-harmonic.
In fact we have the following.

 \Lemma{\JJJ.7} {\sl  Let $f$ be a function on the unit sphere in $S\ss\rn$
 with the property that $\max\{f, g^*f\} \in C^1(S)$ for all $g\in G$. 
 Then $f=$ constant.}
 \pf
 We begin with the case $G={\rm O}(n)$. If we can prove constancy on every
 great circle in $S^{n-1}$, we are done. 
So we are  immediately reduced to the case $n=2$.
 Lifting to the covering $\bbr\to S^1$ we are then reduced to the following
elementary fact:
 \smallskip
 \centerline{\sl
 Let $f:\bbr\to\bbr$ be a  $2\pi$-periodic function with the property
that  for all $a\in \bbr$, }

\centerline{\sl
the function 
$
F_a(x) \ \equiv \ \max\{f(x), f(x+a)\}
$
is differentiable.  Then $f\equiv$ constant.}
\smallskip   \def\gerH{{\fr{\hbox{h}}}}\def\gerK{{\fr{\hbox{k}}}}\def\gerS{{\fr{\hbox{s}}}}\def\gerO{{\fr{\hbox{o}}}}
\noindent
We see this as follows.   If $f$ is not constant, there exists a point $x$ with $f'(x)>0$. Since it is 
periodic, there must also exist a point $y$ with $f'(y)<0$.  Set $a=y-x$.  
Then the left hand derivative of $F_a$ is $<0$ (if it exists), and the 
right hand one is $>0$. This completes the argument for $G={\rm O}(n)$.

Consider now the general case of a closed subgroup $G\ss {\rm O}(n)$. Fix   
$x\in S^{n-1}$ and decompose the Lie algebra as $\gerG = \gerK \oplus \gerH$
(orthogonal with respect to the Killing form of $\gerS\gerO (n)$),
where $\gerK =   \gerG \cap \gerS\gerO(n-1)$ is the Lie algebra of the subgroup
$K\equiv \{g\in G : g(x)=x\}$. Now the differential of the $G$-action at $x$ gives
an isomorphism $\gerG \cong T_x(S^{n-1})$ and every 1-parameter subgroup $\vf_t\ss G$
generated by an element of $\gerG$, the orbit is a great circle. The argument made above 
for O$(n)$ now applies, and Lemma \JJJ.7 is proved.\qed

\medskip

Taken together, these two lemmas  prove that the punctured harmonic $U(x)$ is radial
 (constant on spheres about the origin).
Therefore, by Proposition \CC.5, $U=\Theta K+C$, and by (\PP.1), $C=0$.  This completes the proof of Proposition \JJJ.3 if $U\in C(\rn-\{0\})$.
\medskip

For the next step we establish the following strengthening of Proposition \BBB.5
which reduces the case $U\in L^\infty_{\rm loc}(\rn-\{0\})$ to the case $U\in C(\rn-\{0\})$.


\Prop{\JJJ.8} {\sl
Suppose $F$ is uniformly elliptic.  Then each locally bounded $F$-maximal function is
$F$-harmonic.
}
\pf  Suppose $u$ is an $F$-maximal $L_{\rm loc}^\infty$-function on a domain $X\ss \rn$.
By Theorem \BBB.7 for any compact set $K\ss X$, $u$ is  the decreasing limit of a sequence $\{u_j\}_j$
of $F$-harmonic functions on  a neighborhood of $K$. By Fact \JJJ.4, the limit $u$
is continuous, and hence $F$-harmonic by Proposition \BBB.5.\qed

\medskip
This completes the second stage of the proof of Proposition \JJJ.3 where $U\in L^\infty_{\rm loc}(\rn-\{0\})$.
It remains to prove the last stage where $U$ is a general tangent.

By Lemma \JJJ.6(b), for each $N>0$ we have 
$U^N \equiv \max\{U, NK_p\} \in T_0( \max\{u, NK_p\} )$.
Since $U^N \in L^\infty_{\rm loc}(\rn-\{0\})$, $U^N$ is a multiple  of $K_p$.
We now observe that $U^N$ decreases down to  $U$ as $N\to\infty$. Hence, if each $U^N$ is a multiple of the Riesz kernel, then so is $U$.  This completes the proof of Proposition \JJJ.3.\qed
\medskip

The last result needed for the proof of Theorem \JJJ.1 in the O$(n)$-invariant case is the following
Proposition, which reduces the case of our general $F$ of characteristic $p$, to a specific 
maximal such equation, which is uniformly elliptic.

\Prop {\JJJ.9}  {\sl
The subequation
$$
\cp^{\rm largest}_p \ \equdef\ \left\{A :    A+{p-1\over n-p}(\tr A) I\geq0      \right\}
$$
contains all the O$(n)$-invariant {\rm convex}  cone subequations $F$ 
of Riesz characteristic $p$, and has  Riesz characteristic $p$ itself.  Since 
$$
\cp^{\rm largest}_p \ =\ \cp(\d) \qquad{\rm with} \ \ \d\ =\ {(p-1)n\over n-p}
$$
(see Example \HH.3), the subequation $\cp^{\rm largest}_p$ is uniformly elliptic
when $p>1$.
}

\pf
Suppose $A=\l_1 P_{e_1} +\cdots+ \l_n P_{e_n}$ is in diagonal
form with $\l_1\leq \cdots\leq \l_n$. Then by definition (\HH.5) we know that
$$
A\ \notin\ \cp(\d)
\qquad\iff\qquad
\bra A {P_{e_1} + \smfrac \d n I} \ =\ 
\l_1 + \smfrac \d n \left (\l_1 +\cdots + \l_n  \right) \ <\ 0.
$$
If $\mu' = \pi(\l')$ is a permutation of $\l' = (\l_2,...,\l_n)$, then
$A_\pi \equiv \l_1 P_{e_1} + \mu_2 P_{e_2} +\cdots+ \mu_n P_{e_n}$
also belongs to the open half-space $H$ defined by $\bra{A}{ P_{e_1} +{\d\over n} I} < 0$,
and $H$ is disjoint from $\cp(\d)$.  Averaging $A$ over these permutations yields
$B\equiv  \l_1 P_{e_1} + {\Sigma\over n-1}p_{e_1^\perp}$ where $\Sigma\equiv \l_2+\cdots+\l_n$.
Since $B\in H$ we have $B\notin \cp(\d)$.  Hence setting $e\equiv e_1$ and using the fact
that $\cp(\d)$ is a cone, we can rescale to obtain $B'\equiv P_{e^\perp} - (p'-1) P_e \notin \cp(\d)$.
Since the characteristic of $\cp(\d)$ is equal to $p$,
this proves that $p'>p$.

Now if $A\in F$, then since $F$ is O$(n)$-invariant and convex, the average $B\in F$.
Finally since $F$ is a cone, $B' \in F$. Since $p'>p$, this proves that $F$ has Riesz characteristic
$>p$, contrary to assumption.
\qed
\medskip

Proposition \JJJ.9 says that if $U$ is a tangent to an $F$-subharmonic function,
where $F$ satisfies the hypotheses, then $U$ is $\cp^{\rm largest}_p$-tangent.
Since the subequation $\cp^{\rm largest}_p$ is uniformly elliptic, Proposition \JJJ.3 applies, which completes
the proof of Theorem \JJJ.1  in the orthogonally invariant case.

\Note{\JJJ.10}  Some (in fact, many) readers may be uncomfortable with the
assertion that $\cp(\d)$-harmonics have the regularity of viscosity solutions to equations
which are convex and uniformly elliptic in the conventional sense. A  discussion
of this point is given in Appendix B.
\medskip

 Consider now the complex analogue $F^\bbc$ of $F$ on $\bbc^n$.
 Then we have $F^\bbc \ss \cp^\bbc(\d)$, the complex analogue of the
 subequation defined in Proposition \JJJ.9. Now for  any $A\in \Sym_\bbr(\bbc^n)$
 one has that $\tr(A) = 2\tr_\bbc(A_\bbc)$ and $\l_1(A) \leq \l_1^\bbc(A_\bbc)$.
 Hence, $\cp({\d\over 2})\ss   \cp^\bbc(\d)$ as subsets of $\Sym(\bbr^{2n})= \Sym_\bbr(\bbc^n)$.
   It follows that $\cp^\bbc(\d)$ is uniformly elliptic (for $p>1$).
 The arguments given above now  go through to establish the theorem in this case.

 The case of the quaternionic analogue $F^\bbh$ is proved in exactly the same
 way.  This completes the proof of Theorem \JJJ.1.\qed


\bigskip

 For the interested reader we present a second argument for  Theorem \JJJ.1
 where the passage from maximal to harmonic is based on regularization via the group $G$ -- a technique
 which is discussed,  for example,  in [HS]. 
 
 \medskip
 \noindent
 {\bf A Slightly Different Proof of Theorem \JJJ.1}   
 Let $u$ be $F$-subharmonic
in a neighborhood of the origin and choose $U\in T_0(u)$.
 For clarity  of exposition we work in the case $p>2$.
Then
$$
U(x) \ =\ \lim_{j\to\infty} r_j^{p-2}u(r_j x)
$$
for a sequence $r_j\downarrow 0$.  Let $\chi = \chi_\e:G\to [0,\infty)$ be a family of smooth functions converging to the 
$\d$-function at the identity in $G$, and for any function $f$ which is $\lloc$ in $\rn-\{0\}$
and in $L^1(S^{n-1}(r))$ for all $r$, define
$$
f^\e(x) \ \equiv \ \int_G f(gx) \chi(g)\, dg
$$
where $dg$ is Haar measure with unit volume on $G$. The following lemma is proved below.

\Lemma {\JJJ.11}
$$
U^\e(x) \ =\ \lim_{j\to\infty}  r_j^{p-2}u^\e(r_j x)
$$

Now by the Fubini Theorem, $U^\e$ satisfies
$$
\eqalign
{
S(U^\e, r) \ &=\ \int_{|x|=1} U^\e(rx) \, dx \ =\ \int_{|x|=1} \left\{ \int_G U(grx)  \chi(g)\, dg\right\}\, dx\cr
&=\  \int_G \left\{\int_{|x|=1} U(rgx) \, dx \right\} \,\chi(g) \, dg 
\ =\  \int_G \, S(U,r) \,\chi(g) \, dg \cr
&=\ S(U,r) \ =\ \den^S K(r).
}
$$
From this we conclude that  $U^\e$ is maximal by Theorem  \BBB.2.  The next lemma is also proved below.

\Lemma {\JJJ.12}  {\sl $U^\e$ is continuous and converges to $U$ in $\lloc(\rn-\{0\})$ as $\e\to0$.}

\medskip

Note that the continuity of $U^\e$ implies that
it is $F$-harmonic (Proposition \BBB.5).

We now fix $g_0\in G$ and define
$$
V^\e(x) \ \equiv \ U^\e(g_0x)  \ =\ \lim_{j\to\infty}  r_j^{p-2}u^\e(r_j g_0x)
$$
where the second equality comes from Lemma \JJJ.11. Clearly $V^\e$ is a tangent, and it satisfies
$S(V^\e, r) = S(U^\e,r)=\den^S K(r)$.  In particular, $V^\e$ is also maximal.
Furthermore, note that 
$$
\max\{U^\e(x), V^\e(x)\}\ \equiv\ \lim_{j\to\infty} r_j^{p-2} \max\{u^\e(r_j x), u^\e(r_j g_0x)\}
$$
is also a tangent and hence maximal. We have proved the following.

\Prop {\JJJ.13}  {\sl For all $g\in G$ and all $\e>0$ the function
$\max\{U^\e, g^*U^\e\}$  is  $F$-harmonic}.
\medskip

As in the first proof we now apply elliptic regularity and Lemma \JJJ.7 to conclude that 
 each function $\max\{U^\e, g^*U^\e\}$ is $C^1$, and therefore 
 that $U^\e$ is constant on each sphere.
Then by Corollary \PP.3   $U^\e$ is an increasing radial harmonic
and therefore a multiple of the Riesz kernel.
 Since $U^\e \to U$ in $\lloc$, we conclude that 
 $U = \den^S(u) K(|x|)$. This completes our second proof in the orthogonally invariant case.
 Arguments for the complex and quaternionic analogous proceed as above.\qed

\medskip
\noindent
{\bf Proof of Lemma \JJJ.11.}  Let  $U_j(x) \equiv  r_j^{p-2}u(r_j x)$, so that $U_j\to U$ in 
$\lloc(\rn-\{0\})$.  Set $A=\{r\leq |x|\leq R\}$. Then
$$
\eqalign
{
\left\|  U_j^\e -  U^\e \right\|_{L^1(A)} \ 
&=\ \int_A \left|  \int_G\left\{ U_j(gx) \chi(g) - U(gx)\chi(g)\right\}\,dg            \right| \, dx \cr
&\leq \ \int_G \int_A \left| U_j(gx) - U(gx)\right| \ dx \, \chi(g)\, dg    \cr
&=\ \int_G \left\| g^* U_j - g^* U \right\|_{L^1(A)}\, \chi(g)\,dg  \cr
&=\ \int_G \left\|  U_j -  U \right\|_{L^1(A)}\, \chi(g)\,dg  
\ =\ \left\|  U_j -  U \right\|_{L^1(A)}
}
$$
Thus $\lim_{j\to\infty} U^\e_j = \left\{\lim_{j\to\infty} U_j\right\}^\e$ as claimed.\qed

\medskip
\noindent
{\bf Proof of Lemma \JJJ.12.}   It is standard that the restriction of $U^\e$ to each sphere
$\{|x|=r\}$ is continuous (in fact, smooth). We see this as follows. Suppose $x_j \to x$
in $\{|x|=r\}$.  By transitivity we can write $x_j = g_jx$ where $g_j \to 1$ in $G$.  Then
$$
\eqalign
{
\left|  U^\e(x_j) - U^\e(x)    \right|    \ 
&=\   \left|  \int_G U(gx_j) \chi(g)\,dg  -    \int_G U(gx) \chi(g)\,dg  \right|   \cr
&=\ \left|  \int_G U(gg_jx) \chi(g)\,dg  -    \int_G U(gx) \chi(g)\,dg  \right|   \cr
&=\ \left|  \int_G U(hx) \chi(h g_j^{-1})\,dh  -    \int_G U(gx) \chi(g)\,dg  \right|   \cr
&=\ \left|  \int_G U(gx) \left\{ \chi(g g_j^{-1}) - \chi(g) \right\}\,dg  \right|   \cr
&\leq \ \int_G  \left|U(gx)\right|  \left|  \chi(g g_j^{-1}) - \chi(g)  \right| \,dg\cr
&\leq \ \left\{\int_{\{|x|=r\}}  \left|U(x)\right| \,dg \right\} \sup_{g\in G}\left|  \chi(g g_j^{-1}) - \chi(g)  \right| \ \to\ 0
}
$$
We also know that 
$U^\e$ is maximal, and in particular upper semi-continuous with $S(U^\e,t) \equiv \den K(t)$
for all $t$.

Now  for $|x|=r$,  $g_0\in G$, and any $r_1<r<r_2$, the calculation above  also shows that 
$$
\left|  U^\e(g_0x) - U^\e(x)    \right|   \ \leq \ \sup_{r_1\leq t \leq r_2}
\left\{\int_{\{|x|=t\}}    \left|U(x)\right| \,dg \right\} \sup_{g\in G}\left|  \chi(g g_0^{-1}) - \chi(g)  \right| 
$$
Now every  $y$ with $|y|=t$ and $|y-x|<\d$ can be written as $y=g_0x$
with $d(g_0, 1)<\e(\d)$ where $\e(\d) \to 0$ as $\d\to0$.  Thus we have
$$
\left|  U^\e(y) - U^\e(x)    \right|   \  \leq \ \sup_{r_1\leq t \leq r_2}
\left\{\int_{\{|x|=t\}}    \left|U(x)\right| \,dg \right\} \sup_{d(g_0, 1)<\e(\d) }  \sup_{g\in G}\left|  \chi(g g_0^{-1}) - \chi(g)  \right|    \ \leq C \vf(\d)   
$$
for all $|x|=t, |y|=t$, $|y-x| <\d$ and $r_1\leq t\leq r_2$.
This shows that the family of functions 
$$
V_t \  \equiv \ U^\e(tx) \ \ {\rm is\ uniformly\ equicontinuous\ on \ the\ sphere\ } S^{n-1}= \{|x|=1\}
$$
\medskip
\noindent
{\bf Claim:}  \hskip 1.2 in
$
\lim_{t\to t_0}\sup_{S^{n-1}}|V_t-V_{t_0}| \ =\ 0.
$
\pf
Let $t_j\to t_0$ be any sequence.  Then by the equicontinuity above, there is a subsequence
such that $V_{t_j}$ converges uniformly to a limit $\wt V$ on $S^{n-1}$. We are done if we show
that $\wt V = V_{t_0}$.

Now by the upper semi-continuity of $U^\e$ we have
$$
\wt V(x) \ =\ 
\lim_{j\to\infty} V_{t_j} (x) \ =   \ \lim_{j\to\infty}  U^\e(t_jx) \ \leq \ U^\e(t_0x).
$$
However, we also have that
$$
\int_{S^{n-1}} \wt V (x)\, dx \ =   \ 
\lim_{j\to\infty} \int_{S^{n-1}}V_{t_j} (x)\, dx \ =   \ 
\lim_{j\to\infty}  \int_{S^{n-1}}   U^\e(t_jx) \, dx \ = \   \int_{S^{n-1}}   U^\e(t_0x) \, dx.
$$
since the last two terms are just the averages 
$S(U^\e, t_j) = \den K(t_j) \to \den K(t_0)=S(U^\e, t_0) $.
By the inequality (2) we conclude that $\wt V (x)  =  U^\e(t_0 x)=V_{t_0}(x)$
for all $x\in S^{n-1}$. Thus we have shown that $U^\e$ is continuous for all $\e$.

Now it is a general fact that  $f^\e \to f$ in $\lloc$. The proof is easy and the
convergence is uniform when $f\in C^\infty_0$. The general case follows from
the fact that $C^\infty_0$ is dense in $L^1$ on compact domains.  This completes the proof
of Lemma \JJJ.12.
\qed
 
\Ex{\JJJ.14}  If one drops the convexity hypothesis in Theorem \JJJ.1, then in dimensions $n\geq 3$
there are orthogonally invariant subequations of every finite Riesz characteristic  for which strong uniqueness fails.  
To see this we consider the largest such subequation  of  characteristic $p$:
$$
\cp^{\rm min/max}_p \ \equiv \  \left\{A : \l_{\rm min}(A)  + (p-1) \l_{\rm max}(A) \ \geq\ 0 \right \}.
$$
(See Appendix A in Part II for a proof that there exists a largest and it is the one above.)
To see that strong uniqueness fails for $\cp^{\rm min/max}_p$ we consider the following 
functions.  Write $\rn = \bbr^m\times \bbr^{n-m}, m<n$ with coordinates $z=(x,y)$, and consider
the function 
$$
u(x,y) \ \equiv\ \overline{K}_p(|x|)
$$
where $\overline{K}_p$ is given by (\CC.8).  Then $D^2_zu = {1\over |x|^p} (P_{x^\perp} -(p-1)P_x)$ has
ordered eigenvalues
$$
-{(p-1) \over |x|^p}, \ 0,\  ...\ ,\  0,\  {1\over |x|^p},\  ... \ , \  {1\over |x|^p},
$$ 
from which it is clear that $u$ is $\cp^{\rm min/max}_p$-subharmonic on $\rn$ and, in fact,
 $\cp^{\rm min/max}_p$-harmonic for $x\neq 0$. Note that $u$ has Riesz homogeneity $p$
 and is therefore its own tangent at points of the form $(0,y)$. Hence strong uniqueness 
 fails for  $\cp^{\rm min/max}_p$.
 
Straightforward calculation shows, however, that these ``partial Riesz kernels'' are not
subharmonic for the largest {\sl convex} subequation of characteristic $p$ given in Proposition \JJJ.9
above.


\vfill\eject

\centerline{\headfont \MM.   The Structure of the Sets  E$_c$ where the Density is $\geq$ c.}
\medskip

 In this section we assume  the subequation $F$ on $\rn$ is convex with 
 finite Riesz characteristic $p\geq 2$.
Fix  $u\in F(X)$   where $X$ is an open subset in $\rn$. Let $\den = \den^V: X\to\bbr$ 
be the density function (for the volume function). 
For $c>0$ define 
$$
E_c(u) \ \equiv\ \{x\in X : \den(x)\geq c\}.
$$

For classical plurisubharmonic functions in $\bbc^n$ (where $F = \cp^\bbc$), these sets
have been of central importance.  A deep theorem, due to L. H\"ormander, E. Bombieri 
and  in its final form by Siu ([Ho$_1$], [B], [Siu]), 
states that in this case $E_c$ is a complex analytic subvariety.
One  straightforwardly deduces from this result  that for the subequation $\cp_2$ in $\bbr^{2n}$
the set $E_c$ is discrete, since $\cp^\bbc(J)\ss \cp_2$ for all parallel complex structures $J$
on $\bbr^{2n}$.  

This strong corollary has a quite general extension.

 \Theorem {\MM.1} {\sl Suppose strong uniqueness of tangents holds for $F$ (e.g., $F=\cp_p$).
 Then for any $F$-subharmonic function $u$ the set $E_c(u)$ is discrete.  
 }
 \medskip
 
 This result is essentially sharp.  See Remark \MM.2 below.
 
We will prove Theorem \MM.1 in the following equivalent form.
Consider an $F$-subharmonic function $u$ where $F$ has Riesz characteristic $p$ with $2 <p<\infty$.

  \Theorem { \MM.1$'$}  {\sl Suppose strong uniqueness of tangents holds for  $u$ at a point $x_0$,
  that is,  suppose that  the $p$-flow of $u$ has limit
  $$
  \lim_{r\downarrow 0}   u_r(x_0; x) \ =\ \den K(|x-x_0|) \quad {\rm in} \ \ \lloc(\rn), \ \ {\rm for\ some\ \ }
\den\geq0.
  \eqno{(\MM.1)}
  $$
  Then}
  $$
  \lim_{{\eqalign{&x \to x_0\cr 
  &x \neq x_0
  }}}      \den(u,x)\ =\ 0.
  $$
  \pf
  Suppose the conclusion fails.  Then there exists a sequence $x_j\to x_0$
   with $\den(u,x_j)\geq c >0$ for all $j$.  Assume $x_0=0$, and set $x_j = r_j \s_j$
   with $r_j = |x_j|$. Then $r_j\to 0$,  and passing to a subsequence we can assume
   that $\s_j\to \s \in S^{n-1}$.  The idea now is to apply the sequence of $r_j$-homotheties
   to $u$.  This will give a sequence $u_{r_j}$ of $F$-subharmonics with $\den(u_{r_j}, \s_j)\geq c$.
   With appropriate estimates from monotonicity, this will contradict (\MM.1).
   
   To begin pick $\rho>0$ small, and note that
   $$
   {V \left(u_{r_j}, \s_j, \rho \right)   \over   K(\rho)  }
   \ =\ 
    {V \left(u,  x_j,  r_j \rho \right)   \over   K(r_j \rho)  }
   \eqno{(\MM.2)}
   $$
  since 
  $$
V \left(u_{r_j}, \s_j, \rho \right) \ =\ \intave B u_{r_j} \left(\s_j + \rho x\right)\, dx
 \ =\  r_j^{p-2} \intave B u \left(x_j +  r_j \rho x\right)\, dx
  $$
 and
 $$
 {r_j^{p-2} \over K(\rho) } \ =\ {1\over K(r_j\rho)}.
 $$
  
 Next we show that for all $j$
 $$
  {V \left(u,  x_j,  r_j \rho \right)   \over   K(r_j \rho)  }  \ \geq\ {c\over 2}.
  \eqno{(\MM.3)}
   $$
 In fact, this uniform bound from below, on the convergence of 
  ${V \left(u,  x_j,  t \right)   \over   K(t)  } $ to $\den(u,x_j)$,
  independent of $x_j$, is obtained from the monotonicity property (Theorem \BB.4)
  as follows.  Set $\a \equiv 2^{{1\over p-2}}$.  Fix $x_j$ and abbreviate notation by setting
  $t=r_j\rho$ and $V(t) = V(u,x_j, t)= V(u,x_j,r_j\rho)$.  We now apply  the identity   
  $$
  {V(t)  \over  K(t)} \ \ =\ \ \left[  { V(\a t) - V(t)    \over  K(\a t) - K(t)   }    \right]
  {    
  \left( 1-  {K(\a t)   \over  K(t)} \right)      
  \over       
    \left( 1-  {V(\a t)   \over  V(t)} \right)    
  },
  \eqno{(\MM.4)}
   $$
   with the constant $\a>0$ chosen so 
 that ${K(\a t) \over K(t)} = \a^{-(p-2)} = \half$.
  We assume $u$ and hence $V(t)$ is
  $\leq 0$ which can be obtained by subtracting a constant, or noting that $\lim_{x\to 0} u(x) =-\infty$
  since $\den(u,0)\geq c$ by Theorem \KK.4.)   
  
  Then $V(t) \leq V(\a t) \leq0$ since $V(t)$ is increasing in $t$, which implies that the reciprocal 
  of $1- {V(\a t)\over V(t)}$ is $\geq 1$.

  By Theorem \BB.4 this proves that, as desired,
   $$
  {V(t)   \over   K(t)  }  \ \geq\ {c\over 2}.
  \eqno{(\MM.3)'}
   $$
  Combining (\MM.2) and (\MM.3) we have
  $$
   {V \left(u_{r_j},  \s_j,   \rho \right)   \over   K( \rho)  }  \ \geq\ {c\over 2}.
   \eqno{(\MM.5)}
   $$

\def\intavee{{\int \!\!\!\! \!\!-}}

By the hypothesis (\MM.1) we have
$$
\lim_{r_j\downarrow 0} V \left(u_{r_j},  \s_j,   \rho \right) \ =
\ \lim_{r_j\downarrow 0} \intavee_{B_\rho(\s_j)} u_{r_j} \ =\ 
 \den  \intavee_{B_\rho(\s)}  K(|y|)\, dy.
$$
Therefore, by (\MM.5) 
$$
- \rho^{p-2} \den  \intavee_{B_\rho(\s)}  K(|y|)\, dy \ \geq\ {c\over 2}.
$$
Since
$$
\lim_{\rho\to0}   \intavee_{B_\rho(\s)}  K(|y|)\, dy \ =\ K(1) \ =\ -1,
$$
this implies that $c=0$, a contradiction.\qed

\Remark{\MM.2}  For $F$   as above, any finite set can occur as the set $E_c$ for an
$F$-subharmonic function.  In fact, more is true.  In a separate paper [\HLADP] we construct
$F$-subharmonics with prescribed asymptotics at a finite set of points and prescribed boundary values.

 \Theorem {\MM.3. [\HLADP]}  {\sl Let $\O\ss\rn$ be a domain with smooth boundary
 $\bo$ which is strictly convex (or more generally strictly $F$-convex (cf.  [\HLDD]).
 Let $E=\{x_j\}_{j=1}^N \ss \O$ be a finite subset, and $\{\den_j\}_{j=1}^N $ any set
 of positive real numbers.  Then given any $\vf \in C(\bo)$, there exists a unique
 $u\in  \USC(\ob)$ such that:
 \medskip
 
 (1) \ \ $u$ is $F$-harmonic in $\O-E$,\medskip
 
 (2) \ \ $u\bigr|_{\bo} \ =\ \vf$, and\medskip
 
 (3) \ \ $\den(u, x_j) \ =\ \den_j$ for $j=1,..., N$.
 
}



\vskip.3in

\centerline{\headfont \NN.  Subequations with Riesz characteristic 1 $\leq $ p $< $ 2.  } \bigskip 

When the Riesz characteristic   satisfies $1\leq p <2$,  
the behavior and study of $F$-subharmonics  differs greatly from the case $p\geq 2$.
\bigskip
\centerline
{\bf
$C^{0,\a}$  Regularity of Subharmonics
}
\medskip

To begin,  all $F$-subharmonics (not just the $F$-harmonics)  are regular.

To be completely clear we formulate two hypotheses on a function $u$.
\medskip
\noindent
{\bf Hypothesis A:} \ $u\in F(X)$ where $F$ is a (not necessarily convex) \ST cone subequation with 
characteristic $p<\infty$.

\medskip
\noindent
{\bf Hypothesis B.}  \ $u\in \USC(X)$ satisfies the (MP) and {\bf $K_p$ double monotonicity}, 
that is, for all $y\in X$
$$
{M(u,y,t) - M(u,y,s) \over K_p(t) - K_p(s)} \ \ {\rm is\ non\, decreasing\ in\  } s\ {\rm and}\ t
\eqno{(\NN.1)}
$$
for all $0\leq s<t<\dist(y,\partial X)$.
\medskip

By Theorem \AA.7 and Theorem \BB.4
$$
{\rm  Hypothesis \ A} \qquad\Rightarrow \qquad {\rm  Hypothesis\  B}.
\eqno{(\NN.2)}
$$
Note that under  Hypothesis B the density $\T(u,y)$ exists with $0\leq \T(u,y)<\infty$ for each 
point $y\in X$.  For an arbitrary function $u$, we abbreviate the H\"older norm on a compact set $K$
(allowing the value $+\infty$) by
$$
\| u \|_\a(K) \ \equiv\ \| u\|_{C^{0,\a}(K)}.
\eqno{(\NN.3)}
$$

\Theorem {\NN.1} {\sl 
Assume Hypothesis B.
Then $u$ is locally H\"older continuous on $X$ with exponent $\a\equiv 2-p$.

More specifically,  if $B_{3\rho}(x_0) \ss X$, then
$$
\| u \|_\a \left(B_\rho(x_0)\right) \ \leq\ \left[ {  R^\a \over (R-\rho)^\a - \rho^\a   } \right]
{ M(u, x_0, R) -u(x_0)   \over  R^\a   }
\eqno{(\NN.4)}
$$
for all $0 <  3\rho \leq R <\dist(x_0, \partial X)$.
(In particular, $u(x_0)> -\infty$, i.e., $u$  is finite-valued at each point $x_0\in X$.)
}
\pf
  Assume $x,y\in B_\rho(x_0)$.  Note that $x\in  \partial B_{|x-y|}(y)$.
Hence,
$$
{ u(x)-u(y)   \over |x-y|^\a    } \ \leq\ {   M(u,y,|x-y|) -u(y)  \over   |x-y|^\a      }.
$$
Choose $R\geq 3\rho$.  Since $x,y\in B_{\rho}(x_0)$, we have $|x-y|\leq 2\rho$ and hence
$R\geq |x-y|+\rho$, or $R-\rho\geq |x-y|$. Therefore, by the monotonicity Hypothesis B
$$
{   M(u,y,|x-y|) -u(y)  \over   |x-y|^\a      } \ \leq\ {   M(u,y,R-\rho) - M(u,y,\rho) \over   (R-\rho)^\a - \rho^\a     }.
\eqno{(\NN.5)}
$$
Now $B_{R-\rho}(y) \ss B_R(x_0)$ since $y\in B_\rho(x_0)$.
This proves that
\medskip

(1)\ \ $M(u,y,R-\rho)\ \leq\ M(u, x_0, R)$.

\medskip
\noindent
Also  $x_0\in B_\rho(y)$ and hence $u(x_0)\leq M(u,y,\rho)$, or equivalently
\medskip

(2)\ \ $-M(u,y, \rho)\ \leq\ -u(x_0).$

\medskip
\noindent
Now (1) and (2) imply that $M(u,y, R-\rho) - M(u,y,\rho) \leq M(u, x_0, R) - u(x_0)$
and (\NN.4) follows from (\NN.5).\qed
\medskip

Define the {\bf infinitesimal H\"older norm of $u$ at $x_0$} to be
$$
\|u\|_\a (x_0) \ \equiv\ \lim_{\rho\to\infty} \|u\|_\a \left( B_\rho (x_0)  \right).
\eqno{(\NN.6)}
$$

\Prop{\NN.2} {\sl
Under Hypothesis B,
$$
\|u\|_\a (x_0) \ \leq\  {M(u, x_0, R) - u(x_0)   \over   R^\a }   \ \leq\ \|u\|_\a \left( B_R (x_0)  \right).
\eqno{(\NN.7)}
$$
for all $0< R <\dist(x_0, \partial X)$.
}
\pf
For the first inequality, let $\rho\to 0$ on both sides of the inequality (\NN.4)
in Theorem \NN.1.

By the (MP) there exists $y\in \partial B_R(x_0)$ such that $M(u, x_0, R) = u(y)$,
and hence
$$
 {M(u, x_0, R) - u(x_0)   \over   R^\a } \ =\ {u(y)-u(x_0)  \over  |y-x|^\a}  \ \leq\  \|u\|_\a \left( B_R (x_0)  \right).
$$
Now it is easy to prove  that the infinitesimal H\"older norm and the density are the same thing.

\Cor{\NN.3}
$$
 \|u\|_\a(x_0) \ =\ \T(u, x_0).
$$
\pf
Take the limit as $R\to 0$ in (\NN.7) and apply the definition of the density.\qed

\Remark{\NN.4. (Hypothesis A)}
Lemma  A.1 in part II states that $\cp^{\rm min/max}_p 
\equiv \{A : \l_{\rm min}(A) +(p-1) \l_{\rm max}(A) \geq 0\}$ 
is the maximal subequation of characteristic $p$ -- it contains every other subequation $F$  of characteristic $p$.
Thus the relevance of  Theorem \NN.1 for pure second-order subequations can be stated as follows.

Theorem \NN.1 holds under 
\medskip
\noindent
{\bf Hypothesis A$'$ ($0<\a\leq 1$):}
The function $u$ satisfies the subequation
 $$
 \l_{\rm min}(D^2 u) + (1-  \a)\l_{\rm max}(D^2 u) \ \geq\ 0\qquad{\rm on}  \ \ X
 $$  in the viscosity sense.  Said differently, Hypothesis A and Hypothesis A$'$ are the same.

\Remark {\NN.5}  The subequations $\cp^{\rm min/max}_p$ are never convex unless $p=1$.  In addition we have
$$
\cp^{\rm min/max}_p \ \ss\ \D
\qquad\iff\qquad
p\ \leq\ 1 + {1\over n-1}
\qquad\iff\qquad
{n-2\over n-1} \ \leq\ \a \ \leq\ 1.
$$
To see this note that 
$\l_1+(p-1)\l_n \geq 0 \Rightarrow \l_1+\cdots+\l_n \geq0$ if and only if $p-1 \leq {1\over n-1}$
since $ \l_1+\cdots+\l_n \geq (n-1)\l_1 + \l_n = (n-1)(\l_1+{1\over n-1}\l_n)$.

\vskip .3in

\centerline
{
\bf  Existence of  Tangents }

\bigskip

In the range $1\leq p <2$ the arguments for the existence and structure of tangents
have a different flavor from the case $p\geq2$.
Recall  that in this range the {\bf tangent flow} 
$$
u_r(x)\ = \ {1\over r^\a}(u(rx)-u(0)) \qquad {\rm where}\ \ \a=2-p,
$$ 
is defined by (\FF.1b).

Tangents to subharmonics have only been defined when $F$ is convex (see Definition \FF.3).
However, because of the H\"older continuity when $1\leq p<2$, the definition can be
extended to the more general cone case in Hypothesis A. In fact, Hypothesis B is enough.  Give $C(\rn)$ the topology of uniform convergence
on compact subsets.

\Def{\NN.6. (Tangents)}  Suppose that $u$ satisfies Hypothesis B
in a neighborhood of the origin in $\rn$. For each sequence $r_j \searrow 0$ such that
$$
U \ \equiv\ \lim_{j\to\infty} u_{r_j} \quad{\rm converges\ in\ } \ C(\rn),
\eqno{(\NN.8)}
$$
the limit function $U$ is called a {\bf tangent to $u$ at} 0, and $T_0(u)$ denotes the space
of all such tangents.
\medskip

The version of Theorem  \FFF.1 for $1\leq p<2$  is given as follows.

\Theorem {\NN.7. (Existence of Tangents)} {\sl Suppose $u$  satisfies Hypothesis B on a ball about the origin.
Then for each $\rho >0$ there exists a $\d>0$ such that  the family $\{u_r\}_{0<r\leq\d}$ is
bounded in norm in $C^{0,\a}(B_\rho)$. In fact,
$$
\limsup_{r\downarrow0} \|u_r\| \left (B_\rho\right)\ \leq\  \den^M(u,0)  \qquad\forall\, \rho>0.
\eqno{(\NN.9)}
$$
In particular, the set $\{u_r\}_{0<r\leq\d}$ is precompact in $C(\rn)$.}

\pf
Note that $u_r(0) = 0$ so that  Theorem \NN.1 states that the $\a$-H\"older norm of $u_r$ on $B_\rho$ satisfies
$$
\|u_r\|\left (B_\rho\right)  \ \leq\ { R^\a  \over  (R-\rho)^\a - \rho^\a   } { M(u_r, 0, R) \over R^\a   }
$$
if $rR$ is small and $0< 3\rho\leq R$.
Now by the definition of $u_r$
$$
M(u_r, 0, R)\ =\ {M(u, 0, rR) - u(0)   \over r^\a   },
$$
and therefore
$$
\|u_r\| \left (B_\rho\right)\ \leq\   { R^\a  \over  (R-\rho)^\a - \rho^\a   } { M(u, 0, rR) - u(0)\over (rR)^\a   }.
$$
Taking  the $\limsup$ as $r\downarrow 0$ yields
$$
\limsup_{r\downarrow0}  \|u_r\| \left (B_\rho\right)\ \leq\   { R^\a  \over  (R-\rho)^\a - \rho^\a   } \den^M(u,0).
$$
Finally we can let $R\to\infty$, proving (\NN.9).

By the standard compact embedding theorem this proves that (taking the topology of 
H\"older norms on compact subsets)
$$
\{u_r\}_{0<r\leq\d}\ \ \ {\rm is\ precompact\ in\ \ } C^{0,\b}(\rn)\ \ {\rm for\ each}\ \ 0\leq \b <\a,
\eqno{(\NN.10)}
$$
where  $C^{0,\b}(\rn)=C(\rn)$ when $\b=0$.
\qed
 
 \medskip
 \noindent
 {\bf Note. }  If $F$ is convex, then our previous $\lloc$ Definition \FF.3 of a tangent $U$ to $u$ at 0 is also applicable.  It agrees with Definition \NN.6 because of the precompactness.  
 \medskip

The analogue of Theorem \FFF.2 is the same except that $\lloc(\rn)$ is replaced by $C(\rn)$.

\Theorem{\NN.8}  {\sl
The tangent set $T_0(u)$ to an $F$-subharmonic function $u$ satisfies:\medskip

(1)\ \ $T_0(u)$ is non-empty.
\medskip

(2)\ \ $T_0(u)$ is a compact subset of $C(\rn)$.
\medskip

(3)\ \ $T_0(u)$ is invariant under the tangent flow $U\to U_r$.
\medskip

(4)\ \ $T_0(u)$ is a connected subset of $C(\rn)$.
 }
\medskip

The proof is similar to that of Theorem \FFF.2 and is omitted.

As a consequence of Theorem \NN.8 the H\'older norm of a tangent is finite
on all of $\rn$.

\Cor{\NN.9} {\sl
If $U\in T_0(u)$, then}
$$
\| U \|_\a (\rn) \ =\ \T(u,0)\ =\ \| u \|_\a (x_0).
$$

\bigskip
\centerline{
\bf
 Uniqueness, Strong Uniqueness, and Homogeneity of Tangents
}
\medskip

The three concepts are defined exactly as in Definition \JJ.1.
 For instance,  uniqueness of tangents holds for $u$ at 0 if $T_0(u) =\{U\}$
is   a singleton, or equivalently (cf. (\JJ.1))
$$
\lim_{r \to 0} u_r \quad{\rm exists\ in\ \ } C(\rn) \ \ {\rm and\ equals\ \ } U.
\eqno{(\NN.11)}
$$
Strong uniqueness holds for $u$ at 0 if this limit $U = \T K_p$ where $\T=\T^M(u,0)$.
In this setting strong uniqueness for $u$ is equivalent to the notion of asymptotic equivalence
$u \sim \T |y|^\a$ defined by (\NN.13) below.

\Lemma{\NN.10} {\sl
Strong uniqueness of tangents for $u$ at 0 holds, i.e.,
$$
\lim_{r\to 0} u_r \ =\ \T K_p \ =\ \T |x|^\a \quad{\rm in\ \ } C(\rn)\qquad {\rm with\ \ } \T\geq0
\eqno{(\NN.12)}
$$
if and only if $u(y) \sim \T |y|^\a$, i.e., }
$$
\lim_{y\to 0}  {u(y)  - u(0) \over  |y|^\a}  \ =\ \T \ \geq\ 0.
\eqno{(\NN.13)}
$$

\pf  Actually, the equivalence of (\NN.12) and (\NN.13) is an elementary fact which holds for any continuous function defined in a neighborhood of the origin.

We can assume $u(0)=0$.
 We first show that (\NN.13) $\Rightarrow$ (\NN.12).  The inequality
 $$
 \left|  {u(y)  \over |y|^\a} - \T   \right| \ \leq\ \e
 $$
 can be rewritten, with $y=rx$, as
 $$
\left | u_r(x) - \T |x|^\a \right|  \ \leq\ \e|x|^\a.
 $$
 If the first holds for $|y|\leq \d$, then the second holds for $|x|\leq R$
 and $r\leq \d/R$. Thus we have $\left | u_r(x) - \T |x|^\a \right|   \leq \e R^\a$,
 for all $ |x|\leq R$ and $r\leq \d/R$, which is enough to prove (\NN.13).
 
 For the converse we need only assume that $u_r \to \T K$ uniformly on some sphere $\partial B_R$.
 The inequality 
 $$
 \left | u_r(x) - \T |x|^\a \right|  \leq \e
 $$
 can be rewritten, with $y=rx$, as 
 $$
 \left|  {u(y)  \over |y|^\a} - \T   \right| \ \leq\ {\e  \over |x|^\a}.
 $$
 If the first holds for all $|x|=R$ and $r\leq \d$, then the second holds
 for all $|y|\leq \d R$ with the right-hand side replaced by $\e/R^\a$.  This is enough 
 to prove that $\lim_{y\to0} u(y)/|y|^\a = 0$.\qed

\Note{\NN.11} We say that {\sl strong uniqueness holds for a subequation $F$} if it holds for all
$F$-subharmonics at 0.
Recall that  by Theorem \JJJ.1 Strong Uniqueness of Tangents to subharmonics holds for 
every convex O$(n)$-invariant subequation $F$ with finite Riesz characteristic except $F=\cp$.
This section is only concerned with the cases $1\leq p<2$, or $1<p<2$ when $\cp$ is excluded.
  This includes the subequations: $\cp_p$\  ($1<p<2$), $\Sigma_k$ \  ($p \equiv {n\over k} <2$),
  $\cp(\d)$ \ ($\d < {n\over n-2}$), and  others.

\bigskip

\centerline
{
\bf  Harmonicity of Tangents  when $F$ is convex.}

\medskip

If $F$ is a convex cone \ST\ subequation with finite characteristic, then by
 Theorem \PP.2 every tangent to a subharmonic is maximal, and by Proposition \PP.5,  every
continuous maximal function is $F$-harmonic.  Thus the regularity
result Theorem \NN.1 implies the following for $1\leq p<2$..

\Theorem {\NN.12} {\sl
Let $F$  be as above. Then  for $u$ $F$-subharmonic in a neighborhood of 0,
every tangent $U\in T_0(u)$ is  $F$-harmonic in $\rn-\{0\}$. 
}

\bigskip

\centerline
{
\bf  Removable Point Singularities .}

\medskip

The next result should be compared with  Theorem 1.9 (the case $\a^*<0$) in  [ASS],
where $F$ is assumed to be uniformly elliptic.

\Theorem{\NN.13} 
{\sl
Suppose that $F$ is a cone subequation with a Riesz characteristic  $p$ and $1<p<2$.
Suppose Strong Uniqueness of Tangents holds for $F$ and $F+\cp_p\ss F$
(i.e., $F$ is $\cp_p$-monotone).  For each  function $H$ which is $F$-harmonic 
in a punctured neighborhood of $x_0$ and $F$-subharmonic across $x_0$, one has that
\medskip
\centerline
{
$H$ is $F$-harmonic across $x_0$ $\qquad\iff\qquad$ the density $\T^M(H, x_0)\ =\ 0$.
}
\medskip
}

\pf
Assume that $x_0=0$.  By Proposition A.5 in [\HLADP], the strong uniqueness hypothesis can 
be restated as an asymptotic equivalence  $\lim_{x\to 0} {(H(x)-H(0)) \over |x|^\a}=\T\geq0$,
which was denoted there as $H(x) \sim \den |x|^\a$, at $x_0=0$.

Suppose $\T=0$.  Then for all $\e>0$, $\exists \d>0$ such that $H(x)-H(0) \leq \e|x|^\a$ if
$|x|\leq \d$.  Set $V_\e(x) \equiv -(H(x)-H(0)) + 2\e |x|^\a$.  Then $\e|x|^\a  \leq  V_\e(x)$
on $|x|\leq\d$, which implies that $V_\e$ has no test functions at 0.  Since $\ft+\cp_p\ss \ft$,
the Addition Theorem  (cf. [\HLAE])   implies that $V_\e$ is $\ft$-subharmonic on $B_\d-\{0\}$.
Thus $V_\e$ is  $\ft$-subharmonic on $B_\d$.
Since $V_\e$ decreases to $-H(x)+H(0)$ as $\e\to0$, this proves that $-H$ is 
$\ft$-subharmonic on $B_\d$, and hence $H$ is $F$-harmonic.

Suppose $\T>0$.   Then for $0<\e<\T$ there exists $0<\d<1$ with 
$\e|x|^\a \leq H(x) -H(0)$ on $B_\d$.  Therefore,
$-(H(x)-H(0)) \leq -\e |x|^\a \leq  -\e |x|^2$ if $|x|\leq \d$,
which proves that  $ -\e |x|^2$ is a test function for $-H(x)$ at 0, and hence
$-H$ is not subaffine.  Finally, $0\in F \Rightarrow \cp\ss F \Rightarrow \ft \ss\cpt$, which proves that
$-H$ is not $\ft$-subharmonic.\qed

\vfill\eject


\centerline{\headfont Appendix A.  Subaffine Functions and a Dichotomy.  } \medskip 

For punctured radial subharmonics, i.e., a radial $F$-subharmonic function defined
on a ball, there is a useful dichotomy between those which are increasing 
and those which are decreasing, which we now discuss.  The subaffine equation
$\cpt = \{\l_{\rm max}\geq0\}$
is an important special case, since it contains every subequation $F$ (including itself) for  which
the maximum principle holds.  It is also  a special case in that the radial subequation
$R_{\cpt}$ on $(0,\infty)$ is constant coefficient.  Using the jet variables $(\l,a)$, we have
$$
R_{\cpt}\ =\ \wt {  \bbr_+ \times \bbr_+} \ \equiv\  \{(\l,a) : \ {\rm either} \ \l \geq 0\ {\rm or}\ a\geq 0\}.
\eqno{(A.1)}
$$
It is important to note that the maximum principle holds for this one-variable subequation.

This dual  subequation  $\wt {  \bbr_+ \times \bbr_+}$  is more restrictive than  one might guess.
The next result shows that near the left endpoint of $(a,b)$ there is a dichotomy for a subharmonic.
It is either increasing or it is convex and decreasing.

\Lemma{A.1. (Increasing/Decreasing)} 
{\sl
Suppose that $\psi$ is a general upper semi-continuous  $\wt {  \bbr_+ \times \bbr_+}$-subharmonic
function on an open interval $(a,b)$.  Then either
\medskip

(1)\ \ $\psi$ is increasing on $(a,b)$, or

\medskip

(2)\ \ $\psi$ is decreasing and convex on $(a,b)$, or

\medskip

(3)\ \   $\exists \,c \in (a,b)$ such that
$\psi$ is decreasing and convex on $(a,c)$ and increasing on $(c,b)$. 
}
\pf
Suppose that $\psi$ is not increasing on all of $(a,b)$, that is, $\psi(r) > \psi(s)$ for some 
$a<r<s<b$.  We claim  that $\psi$ is decreasing on $(a,r)$.
If not, there exist $r_1,r_2$ with $a<r_1<r_2<r$ and $\psi(r_1) < \psi(r_2)$.
If $\psi(r_2) < \psi(r)$, then  since $\psi(r) > \psi(s)$, $\psi$ has a strict maximum 
on $(r_2,s)$.
Thus $\psi(r_2) \geq \psi(r)>\psi(s)$, and since $\psi(r_1) < \psi(r_2)$, we must have a
strict maximum on $(r_1,s)$. 

Suppose further that $\psi$ is not decreasing on all of $(a,b)$, that is, $\psi(s) < \psi(t)$ for some 
$r<s<t<b$. The argument above shows that there exists  a maximal 
$c\in (s,t)$ so that $\psi$ is decreasing on $(a,c)$.
Now $\psi$ must be increasing on $(c,b)$ for if not, it would have a strict interior maximum
on that interval.


When $\psi$ is decreasing on $(a,c)$, it must be convex there.  To see this let
$\vf$ be a test function for $\psi$ at $t_0\in(a,c)$. 
Then $0\leq \psi(t)-\psi(t_0) \leq \vf(t)-\vf(t_0)$ for $t<t_0$. This implies that $\vf'(t_0) \leq 0$.
If   $\vf'(t_0) = 0$, then the same inequality implies that
$\vf''(t_0) \geq 0$.  On the other hand, if $\vf'(t_0)<0$, then $\vf''(t_0)\geq 0$  because
$\psi$ is  $\wt {  \bbr_+ \times \bbr_+}$-subharmonic.
\qed
\medskip

We say that the maximum principle (MP) holds for a subequation $F$ if it holds
for all $F$-subharmonic functions.

\Theorem{A.2} {\sl
The following conditions on a subequation $F\ss\Symn$ are equivalent.

\medskip
(1)\ \ The maximum principle holds for $F$.

\medskip
(2)\ \ $F\ss \cpt$  (i.e., the subequation $\cpt$ is universal for (MP)).

\medskip
(3)\ \ $0\notin \Int F$.

\medskip
(4)\ \ $R_F \subseteq \wt {  \bbr_+ \times \bbr_+}$.
}

\pf
Parts (1) -- (3)  were proved in  [\HLDD, Lemma 2.2 and Proposition 4.8].
For part (4) note that $F\ss\cpt \  \Rightarrow \  (R_F)_t \ss(R_{\cpt})_t =\wt{\bbr_+ \times \bbr_+}$.
If $F$ is not contained in $\cpt \equiv \{A : \l_{\rm max}(A)\geq0\}$, 
then there exists $B<0$  with $B\in F$.  By positivity $-\e I \in F$ for some $\e>0$, which implies that 
$(R_{\cpt})_t$ is not contained in $\wt{\bbr_+ \times \bbr_+}$.
\qed

\medskip

These two results can be combined as follows.
\Cor{A.3} {\sl
If the (MP) holds for $F$, then the conclusions (1), (2) and (3) of the Increasing/Decreasing Lemma A.1 hold for any  radial $F$-subharmonic function $u(x) = \psi(|x|)$ defined on an annulus. (In particular, if $u$ is $F$-subharmonic on
a ball, then $\psi(t)$ must be increasing.)
}

\pf  By Theorem \AA.4 and Theorem A.2, $\psi$ is $\wt{\bbr_+ \times \bbr_+}$-subharmonic, and hence Lemma A.1 applies to $\psi$.\qed

\vskip .3in
 

\centerline{\headfont Appendix B.    Uniform Ellipticity and $\cp(\d)$.  } \medskip 

The point of this section is to make clear that viscosity harmonics for the subequation
$$
\cp(\delta') \ =\ \left \{A \in\Symn : A+{\delta}\tr(A) \geq0\right\} \qquad \d={\d'\over n}
$$
are solutions to a uniformly elliptic equation $F(D^2u) =0$ as defined in 
[CC], [T], [CIL], etc.  We define the operator
$$
F:\Symn \ \arr\ \bbr \qquad{\rm by}\qquad  F(A) \ \equiv\ \l_{\rm min}(A) +{\d}\tr(A).
$$
It is straightforward to verify that for all $P\geq 0$ one has
$$
{\delta}\, \tr(P) \ \leq\ F(A+P) -F(A) \ \leq\ \left(1+{\delta}\right) \tr(P).
$$
which is one of the standard equivalent versions of uniform ellipticity for the operator $F$
appearing in the sources above. 

Now since 
$$
\cp(\d') = \{A : F(A)\geq0\}
\and
\Int\cp(\d')\ =\ \{A : F(A)>0\}
$$
 it is completely straightforward to verify that a continuous function $u$ is a viscosity
 solution of $F(D^2 u)=0$ if and only if (in our terminology) $u$ is $\cp(\d')$-harmonic.

\vfill\eject


\centerline{\bf REFERENCES}

\vskip .2in

\noindent
\item{[A$_1$]}   S. Alesker,  {\sl  Non-commutative linear algebra and  plurisubharmonic functions  of quaternionic variables}, Bull.  Sci.  Math., {\bf 127} (2003), 1-35. also ArXiv:math.CV/0104209.  

\smallskip

\noindent
\item{[A$_2$]}   \ \----------,   {\sl  Quaternionic Monge-Amp\`ere equations}, 
J. Geom. Anal., {\bf 13} (2003),  205-238.

 ArXiv:math.CV/0208805.  

\smallskip

\noindent
\item{[AV]}    S. Alesker and M. Verbitsky,  {\sl  Plurisubharmonic functions  on hypercomplex manifolds and HKT-geometry},  J. Geom. Anal. {\bf 16} (2006), No.\  3, 375Ð399.

\smallskip

\noindent
\item{[AS]}  S. N.   Armstrong, B.  Sirakov  and  C. K. Smart,  {\sl Fundamental
solutions of homogeneous fully nonlinear elliptic equations}, Comm.
Pure. Appl. Math., 64 (2011), No.\   6, 737-777.

\smallskip

\noindent
\item{[BT]}   E. Bedford and B. A. Taylor,  {The Dirichlet problem for a complex Monge-Amp\`ere equation}, 
Inventiones Math.{\bf 37} (1976), No.\  1, 1-44.

\smallskip

\noindent
\item{[B]}   
E. Bombieri,  {\sl
Algebraic values of meromorphic maps.}
Invent. Math.  {\bf  10}  (1970), 267--287. 

\smallskip

\noindent
 \item{[CC]}    L. Caffarelli and X. Cabre,  
Fully Nonlinear Elliptic Equations,  Colloquium Publications vol. 43, Amer. Math. Soc., Providence, RI, 1995.

\smallskip

\noindent
 \item{[CNS]}    L. Caffarelli, L. Nirenberg and J. Spruck,  {\sl
The Dirichlet problem for nonlinear second order elliptic equations, III: 
Functions of the eigenvalues of the Hessian},  Acta Math.
  {\bf 155} (1985),   261-301.

 \smallskip

\noindent
\item{[C]}   M. G. Crandall,  {\sl  Viscosity solutions: a primer},  
pp. 1-43 in ``Viscosity Solutions and Applications''  Ed.'s Dolcetta and Lions, 
SLNM {\bf 1660}, Springer Press, New York, 1997.

 \smallskip

\noindent
\item{[CE]}  
M. G. Crandall  and    L. C.  Evans, {\sl A remark on infinity harmonic functions}, 
Proceedings of the USA-Chile Workshop on Nonlinear Analysis (Vi–a del Mar-Valparaiso, 2000), 123Ð129, Electron. J. Differ. Equ. Conf., {\bf 6}, Southwest Texas State Univ., San Marcos, TX, 2001.

\smallskip

\noindent
  \item{[CEG]} M. G. Crandall,  L. C.  Evans,  and R. F.  Gariepy,  {\sl
  Optimal Lipschitz extensions and the infinity Laplacian.},  Calc. Var. Partial Differential Equations 13 (2001), No.\   2, 123Ð139.

 \smallskip

\noindent
\item{[CIL]}   M. G. Crandall, H. Ishii and P. L. Lions {\sl
User's guide to viscosity solutions of second order partial differential equations},  
Bull. Amer. Math. Soc. (N. S.) {\bf 27} (1992), 1-67.

 \smallskip

\noindent
\item{[D]} J.-P.  Demailly, {\sl Complex analytic and differential geometry}.  An e-book available at:
http://www-fourier.ujf-grenoble.fr/~demailly/documents.html.

\smallskip

\noindent
\item{[E]}  L. C. Evans  
 {\sl Regularity for fully nonlinear elliptic equations and motion by mean curvature}, 
pp. 98-133 in``Viscosity Solutions and Applications''  Ed.'s Dolcetta and Lions, 
SLNM {\bf 1660}, Springer Press, New York, 1997.

\smallskip

\noindent
\item{[ES]}  L. C. Evans   and C. K. Smart,
 {\sl Everywhere differentiability of infinity harmonic functions}, 
Calc. Var.  {\bf 42} (2011), 289Ð299.

\smallskip

\noindent
\item{[F]} H. Federer,   
Geometric measure theory. {\sl Die Grundlehren der mathematischen Wissenschaften}, Band 153,Springer-Verlag New York Inc., New York, 1969.

\smallskip

\noindent
\item{[G]}   L. G\aa rding, {\sl  An inequality for hyperbolic polynomials},
 J.  Math.  Mech. {\bf 8}   No.\   2 (1959),   957-965.

 \smallskip

\noindent
\item{[H]}    F. R. Harvey,   {\sl Removable singularities and structure theorems for positive currents}. Partial differential equations (Proc. Sympos. Pure Math., Vol. XXIII, Univ. California, Berkeley, Calif., 1971), pp. 129-133. Amer. Math. Soc., Providence, R.I., 1973.

\smallskip

 \noindent 
\item {[\HLCG]}   F. R. Harvey and H. B. Lawson, Jr,  {\sl Calibrated geometries}, Acta Mathematica 
{\bf 148} (1982), 47-157.

 \smallskip

\item {[\HLPCG]}  \ \----------, 
 {\sl  An introduction to potential theory in calibrated geometry}, Amer.\ J.\ Math.\  {\bf 131} No.\   4 (2009), 893-944.  ArXiv:math.0710.3920.

\smallskip

\item {[\HLPCGG]}  \ \----------,   {\sl  Duality of positive currents and plurisubharmonic functions in calibrated geometry},  Amer. J. Math.    {\bf 131} No.\   5 (2009), 1211-1240. ArXiv:math.0710.3921.

\smallskip

\item {[\HLDD]}   \ \----------,  {\sl  Dirichlet duality and the non-linear Dirichlet problem},    Comm. on Pure and Applied Math. {\bf 62} (2009), 396-443. ArXiv:math.0710.3991

\smallskip

\item {[\HLPUP]}  \ \----------,    {\sl  Plurisubharmonicity in a general geometric context},  Geometry and Analysis {\bf 1} (2010), 363-401. ArXiv:0804.1316.

\smallskip

\item {[\HLDDR]}  \ \----------,   {\sl  Dirichlet duality and the nonlinear Dirichlet problem
on Riemannian manifolds},  J. Diff. Geom. {\bf 88} (2011), 395-482.   ArXiv:0912.5220.
\smallskip

\item {[\HLHP]}  \ \----------,    {\sl  Hyperbolic polynomials and the Dirichlet problem},   ArXiv:0912.5220.
\smallskip

\item {[\HLHPP]}  \ \----------,   {\sl  G\aa rding's theory of hyperbolic polynomials},  
 Comm. Pure Appl. Math. {\bf 66} (2013), No.\ 7, 1102-1128.

\smallskip

\item {[\HLREST]}  \ \----------, {\sl  The restriction theorem for fully nonlinear subequations}, 
  Ann. Institut Fourier  {\bf 64} No.\ 1 (2014), p. 217-265.
ArXiv:1101.4850.

\smallskip


\item {[\HLPCON]}   \ \----------,  {\sl  p-convexity, p-plurisubharmonicity  and the Levi problem },
   Indiana Univ. Math. J.  {\bf 62} No.\ 1 (2014), 149-170.  ArXiv:1111.3895.

\smallskip

\item {[\HLSURVEY]}   \ \----------,  {\sl  Existence, uniqueness and removable singularities
for nonlinear partial differential equations in geometry},  pp. 102-156 in ``Surveys in Differential Geometry 2013'', vol. 18,  
H.-D. Cao and S.-T. Yau eds., International Press, Somerville, MA, 2013.
ArXiv:1303.1117.

 \smallskip


\item  {[\HLBP]} \ \----------, {\sl The equivalence of  viscosity and distributional
subsolutions for convex subequations -- the strong Bellman principle},
 Bulletin Brazilian Math.\ Soc.\  {\bf 44} No.\ 4 (2013),  621-652.  ArXiv:1301.4914.

\smallskip

\noindent
\item  {[\HLAE]} \ \----------, {\sl The AE Theorem and Addition Theorems for quasi-convex functions}, 
 ArXiv: 1309:1770.

\smallskip



\noindent
\item  {[\HLTangII]} \ \----------, {\sl Tangents to subsolutions -- existence and uniqueness, II}. 
 ArXiv:1408.5851.

\smallskip

\noindent
\item  {[\HLADP]} \ \----------, {\sl The Dirichlet Problem with Prescribed Interior Singularities}, 
Advances in Math. (to appear).
ArXiv:1508.02962.

\smallskip

\noindent
\item  {[\HLCLASSICAL]} \ \----------, {\sl Subequation characterizations of various classical functions}
(to appear).

\smallskip

 \noindent
\item{[Ho$_1$]}
L. H\"ormander,  An Introduction to Complex Analysis in Several Variables,  Third edition. North-Holland Mathematical Library, 7. North-Holland Publishing Co., Amsterdam, 1990. 
 
  \smallskip

 \noindent
\item{[Ho$_2$]}
 \ \----------,     Notions of Convexity. Progress in Mathematics, vol. 127. BirkhŠuser Boston, Inc., Boston, MA, 1994.  \smallskip

%

 \noindent
\item{[HS]}
L. H\"ormander,  and R. Sigurdsson,   {\sl
Limit sets of plurisubharmonic functions}, 
Math. Scand. {\bf 65} (1989), No.\   2, 308--320.
\smallskip

   \noindent
\item{[I]}    H. Ishii,    {\sl  On uniqueness and existence of viscosity solutions of fully nonlinear second-order elliptic pde's},    Comm. Pure and App. Math. {\bf 42} (1989), 14-45.

\smallskip

 \noindent
\item{[K$_1$]}
C. Kiselman, {\sl Tangents of plurisubharmonic functions}, 
 International Symposium in Memory of Hua Loo Keng, Vol. II (Beijing, 1988), 157Ð167, Springer, Berlin, 1991.

  \smallskip
 \noindent
\item{[K$_2$]}
 \ \----------,    {\sl
Plurisubharmonic functions and potential theory in several complex variables}. Development of mathematics 1950--2000, 655--714, BirkhŠuser, Basel, 2000.

\smallskip

\noindent
\item{[Kl]} 
M. Klimek, 
Pluripotential theory.
London Mathematical Society Monographs. New Series, 6. Oxford Science Publications. The Clarendon Press, Oxford University Press, New York, 1991.

\smallskip

   \noindent
\item{[Kr]}    N. V. Krylov,    {\sl  On the general notion of fully nonlinear second-order elliptic equations},    Trans. Amer.\ Math.\ Soc.\  (3)
 {\bf  347}  (1979), 30-34.

\smallskip

\item{[La$_1$]}  D.Labutin,
{\sl Isolated singularities for fully nonlinear elliptic equations},  J. Differential Equations  {\bf 177} (2001), No.\ 1, 49-76.

 \smallskip
 
\item{[La$_2$]}  \ \----------, {\sl Singularities of viscosity solutions of fully nonlinear elliptic equations}, 
Viscosity Solutions of Differential Equations and Related Topics, Ishii ed., RIMS K\^oky\^uroku
No.\   1287, Kyoto University, Kyoto (2002), 45-57

\smallskip

\item{[La$_3$]}   \ \----------, {\sl Potential estimates for a class of fully nonlinear elliptic equations}, 
Duke Math. J. {\bf 111} No.\ 1 (2002), 1-49.

\smallskip

\item {[L]}   N. S.  Landkof,   {Foundations of Modern Potential Theory},  Springer-Verlag, New York, 1972.

\smallskip

 \item {[R]} 
R, T.   Rockafellar, Convex analysis. Princeton Mathematical Series, No.\   28,  Princeton University Press, Princeton, N.J.,  1970.

 \smallskip

\item {[Si]} L. Simon, {\sl  Lectures on Geometric Measure Theory},
Proc. Centre for Math. Anal.  {\bf 3}, Australian Natl. Univ., 1983.

\smallskip

\item {[Sh]} J.-P. Sha, {\sl  $p$-convex riemannian manifolds},
Invent.  Math.  {\bf 83} (1986), 437-447.

\smallskip

   \noindent
\item{[S]}  
R. Sigurdsson,  {\sl
Growth properties of analytic and plurisubharmonic functions of finite order},
Math. Scand. {\bf 59} (1986), No.\   2, 235--304.
 
 \smallskip

   \noindent
\item{[Siu]}    Y.-T.    Siu,    {\sl  Analyticity of sets associated to Lelong numbers and the extension on closed positive currents},    Inventiones Math.,
 {\bf  27}  (1974), 53-156.

\smallskip

   \noindent
\item{[T]} 
N.  Trudinger,  {\sl H\"older gradient estimates for fully nonlinear equations},  Proc.
Roy. Soc. Edinburgh  Sec. A {\bf 108} (1988), 57-65.

\smallskip

   \noindent
\item{[TW$_1$]} 
N.  Trudinger and X-J.  Wang,  {\sl Hessian measures. I},  Dedicated to Olga Ladyzhenskaya. Topol. Methods Nonlinear Anal. {\bf 10} (1997), No.\   2, 225--239.

\smallskip

   \noindent
\item{[TW$_2$]} 
  \ \----------,   {\sl Hessian measures. II},  Ann. of Math. (2)  {\bf 150} (1999), No.\   2, 579--604.

\smallskip

   \noindent
\item{[TW$_3$]} 
  \ \----------,   {\sl Hessian measures. III},   J. Funct. Anal. {\bf  193} (2002), No.\   1, 1--23.

\smallskip

\item {[Wu]}   H. Wu,  {\sl  Manifolds of partially positive curvature},
Indiana Univ. Math. J. {\bf 36} No.\   3 (1987), 525-548.
 
  \smallskip

\vfill\eject

\end